\documentclass[reqno]{amsart}
\usepackage{eurosym}
\usepackage{amsmath}
\usepackage{amsfonts}
\usepackage{amssymb}
\usepackage{amsfonts}
\usepackage{epsfig}
\usepackage{caption}
\usepackage{subcaption}
\usepackage{multirow}
\usepackage{csquotes}
\usepackage{stackengine}
\usepackage{makecell}
\stackMath

\newtheorem{theorem}{Theorem}

\newtheorem{definition}[theorem]{Definition}

\newtheorem{lemma}[theorem]{Lemma}

\newtheorem{proposition}[theorem]{Proposition}
\newtheorem{remark}[theorem]{Remark}

\newtheorem{assumption}[theorem]{Assumption}

\begin{document}

\title[Asymptotic condition numbers] {Asymptotic condition numbers for linear ordinary differential equations: the generic real case}
\author[S. Maset] {S.\ Maset \\
	Dipartimento di Matematica, Informatica e Geoscienze \\
	Universit\`{a} di Trieste \\
	maset@units.it}

\begin{abstract}
	{The present paper studies, for a \emph{real} linear ordinary differential equation $y^\prime(t)=Ay(t)$ \emph{in a generic case}, the long-time propagation to the solution $y(t)$ of a perturbation of the initial value. By measuring the perturbations with normwise relative errors, the paper considers a directional pointwise condition number, defined for a specific initial value and for a specific direction of perturbation of this initial value, and a pointwise condition number, defined for a specific initial value and the worst-case scenario for the direction of perturbation. The asymptotic (long-time) behaviors of these two condition numbers are analyzed in depth.}
\end{abstract}

\maketitle



\section{Introduction}

Consider a linear  $n-$dimensional \emph{complex} Ordinary Differential Equation (ODE) 
\begin{equation}
	\left\{ 
	\begin{array}{l}
		y^{\prime }\left( t\right) =Ay\left( t\right) ,\ t\in\mathbb{R}, \\ 
		y\left( 0\right) =y_{0},
	\end{array}
	\right.  \label{ode}
\end{equation}
where $A\in \mathbb{C}^{n\times n}$ and $y_0,y(t)\in\mathbb{C}^n$. Let  $\Vert\ \cdot\ \Vert$ be an arbitrary norm on $\mathbb{C}^n$. Suppose that the initial value $y_{0}\neq 0$ in (\ref{ode}) is perturbed to $
\widetilde{y}_{0}$ with normwise relative error 
\begin{equation*}
\varepsilon :=\frac{\left\Vert \widetilde{y}_{0}-y_{0}\right\Vert }{
\left\Vert y_{0}\right\Vert }.
\end{equation*}
Due to the perturbation of the initial value, the solution $y(t)=\mathrm{e}^{tA}y_0$ is perturbed to $\widetilde{y}(t)=\mathrm{e}^{tA}\widetilde{y}_0$ with normwise relative error 
\begin{equation*}
\delta \left( t\right) :=\frac{\left\Vert \widetilde{y}\left( t\right)
-y\left( t\right) \right\Vert }{\left\Vert y\left( t\right) \right\Vert }.  \label{deltat}
\end{equation*}

The relation between the relative errors $\delta(t)$ and $\varepsilon$ is given by
\begin{equation*}
	\delta \left( t\right) =K\left( t,y_{0},\widehat{z}_{0}\right) \cdot  \label{magnification}
	\varepsilon ,
\end{equation*}
where
\begin{equation}
	K\left( t,y_{0},\widehat{z}_{0}\right) :=\frac{\left\Vert \mathrm{e}^{tA}
		\widehat{z}_{0}\right\Vert }{\left\Vert \mathrm{e}^{tA}\widehat{y}
		_{0}\right\Vert }.  \label{Ktayz0}
\end{equation}
Here, $\widehat{z}_0\in\mathbb{C}^n$ is the \emph{direction of perturbation}, i.e., the unit vector ($\Vert \widehat{z
}_0 \Vert=1$)  defined by writing the perturbed value $\widetilde{y}_0$ of $y_0$ as
\begin{equation*}
\widetilde{y}_0=y_0+\varepsilon \Vert y_0\Vert \widehat{z}_0,
\end{equation*}
and
$
\widehat{y}_{0}:=\frac{y_{0}}{\left\Vert y_{0}\right\Vert }
$
is the \emph{normalized initial value}. We call $K\left( t,y_{0},\widehat{z}_{0}\right)$ \emph{directional pointwise condition number} of the problem
\begin{equation}
y_0\mapsto \mathrm{e}^{tA}y_0. \label{due}
\end{equation}

Since the direction of perturbation $\widehat{z}_0$ is in general unknown, we introduce the \emph{pointwise condition number} of the
problem (\ref{due}) given by 
\begin{equation}
K\left( t,y_{0}\right) :=\max\limits_{\substack{ \widehat{z}_0\in \mathbb{C
}^{n}  \\ \Vert \widehat{z}_0\Vert =1}}K\left( t,y_{0},\widehat{z}
_0\right) =\frac{\left\Vert \mathrm{e}^{tA}\right\Vert }{\left\Vert \mathrm{e
}^{tA}\widehat{y}_{0}\right\Vert },  \label{KAy0}
\end{equation}
where $\left\Vert \mathrm{e}^{tA}\right\Vert $ is the matrix norm of $
\mathrm{e}^{tA}$ induced by the vector norm $\left\Vert \ \cdot \
\right\Vert $.  It is the worst $
K(t,y_0,\widehat{z}_0)$ by varying $\widehat{
z}_0$. For the definition of condition number for general problems, see \cite{Burgisser2013}.

In the paper \cite{M0}, the asymptotic forms $K_\infty\left( t,y_{0},\widehat{z}_{0}\right)$ and $K_\infty\left( t,y_{0}\right)$, as $t\rightarrow +\infty$, of the condition numbers (\ref{Ktayz0}) and  (\ref{KAy0}), respectively, were stated. Such asymptotic forms are called \emph{asymptotic condition numbers}. The aim of the present paper is to analyze the asymptotic condition numbers for a \emph{real} ODE (\ref{ode}), in a generic case. This analysis is particularly involved in case of a rightmost complex conjugate pair of eigenvalues of the matrix $A$.

Studying how the relative conditioning of the problem (\ref{due}) depends on $t$ is quite important. In fact, it is well known how the absolute error due to a perturbation of the initial value of (\ref{ode}) is propagated to the solution, but it is not equally known how the relative error propagates.

These dynamical aspects of relative conditioning have received adequate attention in the literature only in connection with the problem $A\mapsto \mathrm{e}^{tA}$ (see \cite{Levis1969}, \cite{van2006}, \cite{Kagstrom1977}, \cite{Van1977}, \cite{moler2003}, \cite{Mohy2008}, \cite{Zhu2008}).

\subsection{Plan of the paper}

Besides this introduction, the paper contains seven sections and eight appendices in the Supplementary Material. 

Section 2 analyzes the norms $\left\Vert Q_j(t)u\right\Vert$ and $\left\Vert Q_j(t)\right\Vert$ of the vectors $Q_j(t)u$, $u\in\mathbb{R}^n$, and matrices $Q_j(t)\in\mathbb{R}^{n\times n}$.  Such norms are involved in the asymptotic condition numbers and in establishing the dominance of the asymptotic forms at a finite time. For the case of a rightmost complex conjugate pair of eigenvalues of $A$, Section 3 particularizes the results of Section 2 in depth, for the Euclidean norm as vector norm. Section 4 analyzes the asymptotic condition numbers and shows, for the case of a rightmost complex conjugate pair of eigenvalues of $A$, that they can be factorized as product of a time-constant \emph{oscillating scale factor} and a time-periodic \emph{oscillating term}. Section 5 analyzes the oscillating scale factor and the oscillating term in depth, for the Euclidean norm. In particular, it is analyzed how the extreme values of the oscillating term vary. This is the hardest part of the paper and it requires the extensive material in the appendices. Section 5 also includes numerical examples. Section 6 analyzes the dominance of the asymptotic condition numbers at a finite time. Conclusions are in Section 7. 

The eight appendices contain highly technical material and should be read only by readers interested in the mathematical details. The main text refers to the first two appendices, which in turn refer to the remaining six.

The present paper is a sequel to \cite{M0} and is followed by \cite{Maset2022}. The sequel  \cite{Maset2022} illustrates the results of the present paper through many experimental tests, applications, and other issues, such as the relative conditioning for non-large $t$ and how long it takes for the asymptotic forms to set in.

\section{The norms $\Vert Q_j(t)u\Vert$ and $\Vert Q_j(t)\Vert$}

In the present paper, we consider the case of a real ODE (\ref{ode}), where we have $A\in\mathbb{R}^{n\times n}$
	and $y_0,y(t),\widetilde{y}_0,$ $\widetilde{y}(t)\in\mathbb{R}^n$. Moreover, the
	condition number (\ref{KAy0}) is defined by maximum over 
	$\widehat{z}_0\in \mathbb{R}^{n}$ with $\Vert \widehat{z}_0\Vert=1$. As a consequence, $\Vert \mathrm{e}^{tA}\Vert$ is the real induced norm of $\mathrm{e}^{tA}$.

\subsection{Real and complex induced norms}

For a real matrix $B \in \mathbb{R}^{n \times n}$, we use the real induced norm
\begin{equation}
	\Vert B \Vert := \max_{\substack{u \in \mathbb{R}^{n} \\ \left\Vert u \right\Vert = 1}} \left\Vert Bu \right\Vert. \label{rin}
\end{equation}
For a real row vector $b\in\mathbb{R}^{1\times n}$, we use the real induced norm
\begin{equation}
	\left\Vert b\right\Vert:=\max\limits_{\substack{ u\in \mathbb{R}^{n}  \\ \Vert u\Vert =1
	}}\left\vert bu\right\vert. \label{rinrv}
\end{equation}
For a complex row vector $b\in\mathbb{C}^{1\times n}$, we utilize both the real induced norm
\begin{equation*}
	\left\Vert b\right\Vert_\mathbb{R}:=\max\limits_{\substack{ u\in \mathbb{R}^{n}  \\ \Vert u\Vert =1
	}}\left\vert bu\right\vert  
\end{equation*}
and the complex induced norm
\begin{equation}
\left\Vert b\right\Vert:=\max_{\substack{u \in \mathbb{C}^{n} \\ \left\Vert u \right\Vert = 1}} \left\vert bu \right\vert. \label{cinrv}
\end{equation}
Observe that
\begin{equation}
	\left\Vert b\right\Vert_{\mathbb{R}}\leq \left\Vert b \right\Vert. \label{r<=c}
\end{equation}

\begin{remark} \label{22}
If the vector norm $\Vert \ \cdot \ \Vert $ is a $p$-norm, then, for the real induced norm $\Vert \ \cdot \ \Vert$ of a real row vector $b$, we have 
	$$
	\left\Vert b\right\Vert=\left\Vert
	b^T \right\Vert _{q},
	$$
	where $\Vert \ \cdot \ \Vert_q$ is the vector $q$-norm with $\frac{1}{p}+\frac{1}{q}=1$. The same is true for the complex induced norm of a complex row vector.
\end{remark}

\subsection{$\Lambda_j$ single (simple)  real or single (simple) complex}

We partition the spectrum of $A$ by real parts in the sets $\Lambda_1,\ldots,\Lambda_q$: for $j\in\{1,\ldots,q\}$, $\Lambda_j$ contains all eigenvalues with same real part $r_j$, and $r_1>r_2>\cdots >r_q$ holds. Observe that $\Lambda_1$ is the set of the rightmost eigenvalues of $A$.

\begin{definition}
Let $j \in \{1, \ldots, q\}$. We say that:
\begin{itemize}
	\item $\Lambda_j$  is \emph{single real} if $\Lambda_j$ consists of a real eigenvalue with a single Jordan mini-block of maximum order;
	\item $\Lambda_j$  is \emph{simple real} if $\Lambda_j$ consists of a real simple eigenvalue;
	\item $\Lambda_j$  is \emph{single complex} if $\Lambda_j$ consists of a complex conjugate pair of eigenvalues with a single Jordan mini-block of maximum order;
	\item $\Lambda_j$  is \emph{simple complex} if $\Lambda_j$ consists of a complex conjugate pair of simple eigenvalues.
\end{itemize}
\end{definition}
Of course, if $\Lambda_j$ is simple real, or simple complex, then it is single real, or single complex, respectively.

Observe that it is a generic case for the real matrix $A$ to have all sets $\Lambda_j$
, $j\in\{1,\ldots,q\}$, simple real or simple complex.

Let $j\in\{1,\ldots,q\}$ such that $\Lambda_j$ is single real or single complex. For $\Lambda_j$ single real, we denote by $\lambda_{j}$ the real eigenvalue in $\Lambda_j$. For $\Lambda_j$ single complex, we denote by $\lambda _{j}$ and $
\overline{\lambda _{j}}$ the complex conjugate pair in $\Lambda_j$, where $\lambda_j$ is the eigenvalue with positive imaginary part. This positive imaginary part is denoted by $\omega_j$.

Let $V$ be the matrix whose columns form a Jordan basis for the matrix $A$ and let $W: = V^{-1}$. Let
$$
v^{(j,1)},\ldots,v^{(j,M_j)}
$$
be the unique longest Jordan chain of length $M_j$ corresponding to the eigenvalue $\lambda_{j}$ (real eigenvalue or complex eigenvalue with positive imaginary part in the complex conjugate pair). The vectors of the chain are generalized eigenvectors of $A$ and they appear as columns in the matrix $V$.

Let $v^{(j)}:=v^{(j,1)}\in \mathbb{C}^{n}$ be the first generalized eigenvector of the chain, which is an eigenvector of $A$ corresponding to $\lambda_j$. Let $w^{(j)} \in \mathbb{C}^{1 \times n}$ be the row in the matrix $W$ corresponding to the last generalized eigenvector $v^{(j,M_j)}$ in the chain, i.e., the index of $v^{(j,M_j)}$ as column of $V$ matches the index of $w^{(j)}$ as row of $W$. Observe that $w^{(j)}$ is a left generalized eigenvector  of the matrix $A$ corresponding to the eigenvalue $\lambda_j$.

If $\Lambda_j$ is simple real or simple complex, then $M_j=1$ and $w^{(j)}$ is a left eigenvector  of the matrix $A$ corresponding to the eigenvalue $\lambda_j$.

We can assume (see \cite{M0}) to have a Jordan basis for the matrix $A$ such that: for any $j\in\{1,\ldots,q\}$,
\begin{itemize}
	\item if $
\Lambda_j$ is single real, then $v^{(j)}$
and $w^{(j)}$ are real vectors;
\item if $\Lambda_j$ is single complex, then
\begin{equation}
\overline{v^{(j,1)}},\ldots,\overline{v^{(j,M_j)}}\label{chain}
\end{equation}
is the unique longest Jordan chain of length $M_j$ corresponding to the eigenvalue $\overline{\lambda_{j}}$ and the row in the matrix $W$ corresponding to the last generalized eigenvector $\overline{v^{(j,M_j)}}$ in the chain (\ref{chain}) is $\overline{w^{(j)}}$ (the vectors of (\ref{chain}) appear as columns in the matrix $V$).
\end{itemize}

\subsection{The matrices $Q_j(t)$ and their significance}
We can summarize the results in \cite{M0} of our interest in the following two theorems, where we use the following notations.

For real scalar functions $a(t)$ and $b(t)$ of $t$ and $\epsilon\geq 0$:
\begin{itemize} 
	\item
	\begin{equation}
	a(t)\sim b(t),\ t\rightarrow +\infty, \label{af}
	\end{equation}
	means
	$$
	\lim\limits_{t\rightarrow +\infty}\frac{a(t)}{b(t)}=1;
	$$
	\item
	\begin{equation}
	a(t)\sim b(t)\ \text{with precision}\ \epsilon \label{daf}
	\end{equation}
	means
	$$
	\frac{a(t)}{b(t)}=1+\chi(t),
	$$
	with $\vert \chi(t)\vert\leq \epsilon$.
	\end{itemize}
If (\ref{af}) holds, we say that $b(t)$ is an \emph{asymptotic form} of $a(t)$. The relation (\ref{daf}) shows the dominance of the asymptotic form at the finite time $t$.

Moreover, for $j\in\{1,\ldots,q\}$ such that $\Lambda_j$ is single real or single complex, let
\begin{itemize}
	\item
\begin{equation}
	Q_{j}(t)=Q_{j}:=v^{(j)}w^{(j)}\in\mathbb{R}^{n\times n} \label{Qj0(t)ureal}
\end{equation}
for $\Lambda_j$ single real;
\item \begin{equation}
	Q_{j}(t):= \mathrm{e}^{\mathrm{i}\omega
		_{j}t}v^{(j)}w^{(j)}+ \mathrm{e}^{-\mathrm{i}\omega
		_{j}t}\overline{v^{(j)}}\ \overline{w^{(j)}}=2\mathrm{Re}\left( \mathrm{e}^{\mathrm{i}\omega
		_{j}t}v^{(j)}w^{(j)}\right)\in\mathbb{R}^{n\times n} \label{Qj0(t)ucomplex}
\end{equation}
for $\Lambda _{j}$ single complex.
\end{itemize}
 In (\ref{Qj0(t)ucomplex}), $\mathrm{i}$ denotes the imaginary unit. Moreover, we use the notation $\mathrm{Re}(Z)$, with $Z$ a matrix, for the matrix whose elements are the real parts of the elements of $Z$. The real matrices $Q_j(t)$ are the matrices denoted by $Q_{j0}(t)$ in \cite{M0}.
 
Next theorem states the asymptotic forms of the two condition numbers (\ref{Ktayz0}) and (\ref{KAy0}). Such asymptotic forms are called \emph{asymptotic condition numbers}.
\begin{theorem} \label{one}
	Suppose that $\Lambda_1$ is single real or single complex. If $w^{(1)}y_0\neq 0$ and $w^{(1)}\widehat{z}_0\neq 0$, then
	$$
	K\left( t,y_{0},\widehat{z}_{0}\right)\sim K_\infty\left( t,y_{0},\widehat{z}_{0}\right):=\frac{\left\Vert
		Q_{1}(t)\widehat{z}_{0}\right\Vert }{\left\Vert Q_{1}(t)\widehat{y}
		_{0}\right\Vert },\ t\rightarrow +\infty,
	$$
	and
	$$
		K\left( t,y_{0}\right)\sim K_\infty\left( t,y_{0}\right):=\frac{\left\Vert
		Q_{1}(t)\right\Vert }{\left\Vert Q_{1}(t)\widehat{y}
		_{0}\right\Vert },\ t\rightarrow +\infty.
	$$
\end{theorem}
Next theorem states the dominance of the asymptotic condition numbers at a finite time for the simple real or simple complex case.
\begin{theorem}\label{two}
	Suppose that, for any $j\in\{1,\ldots,q\}$, $\Lambda_j$ is simple real or simple  complex. If $w^{(1)}y_0\neq 0$ and $w^{(1)}\widehat{z}_0\neq 0$, then
	$$
	K\left( t,y_{0},\widehat{z}_{0}\right)\approx K_\infty\left( t,y_{0},\widehat{z}_{0}\right)\ \text{with precision}\ \frac{\epsilon(t,\widehat{z}_0)+\epsilon(t,\widehat{y}_0)}{1-\epsilon(t,\widehat{y}_0)}
	$$
	and
	$$
	K\left( t,y_{0}\right)\approx K_\infty\left( t,y_{0}\right)\ \text{with precision}\ \frac{\epsilon(t)+\epsilon(t,\widehat{y}_0)}{1-\epsilon(t,\widehat{y}_0)},
	$$
	whenever $\epsilon(t,\widehat{y}_0)<1$, where
	$$
	\epsilon(t,u):=\sum\limits_{j=2}^q\mathrm{e}^{(r_j-r_1)t}\frac{\Vert Q_j(t)u\Vert}{\Vert Q_1(t)u\Vert},\ u\in\mathbb{R}^n,
	$$
	and
	$$
	\epsilon(t):=\sum\limits_{j=2}^q\mathrm{e}^{(r_j-r_1)t}\frac{\Vert Q_j(t)\Vert}{\Vert Q_1(t)\Vert}.
	$$ 
\end{theorem}

The previous Theorems \ref{one} and \ref{two} are reformulations of theorems given in \cite{M0} and motivate our interest in the norms $\Vert Q_j(t)u\Vert$ and $\Vert Q_j(t)\Vert$. 

\begin{remark}\label{remof}
	Theorems \ref{one} and \ref{two}  consider a generic case for the real ODE (\ref{ode}). In fact, as previously noted, it is a generic case
	for a real matrix $A$ to have, for any $j\in\{1,\ldots,q\}$, $\Lambda_j$ simple real or simple complex. Moreover, it is a generic case for $y_0$ to satisfy $
	w^{(1)}y_{0}\neq 0$ and it is a generic case for $\widehat{z}_0$ to
	satisfy $w^{(1)}\widehat{z}_{0}\neq 0$.
\end{remark}

In the following, we consider an index $j\in\{1,\ldots,q\}$ such that $\Lambda_j$ is single real or single complex. Indeed, looking at Theorem \ref{two}, for $j>1$ we could consider only $\Lambda_j$ simple real or simple complex.

\subsection{The unit vectors $\widehat{v}^{\left( j\right) }$ and $\widehat{w}^{\left( j\right) }$ and the number $f_j$}
\label{unitvectors}
Let
\begin{equation}
\widehat{v}^{\left( j\right) }:=\frac{v^{\left( j\right) }}{\left\Vert
v^{\left( j\right) }\right\Vert },\ \widehat{w}^{\left(j\right) }:=\frac{
w^{\left( j\right) }}{\left\Vert w^{\left( j\right) }\right\Vert }\ \text{and
}\ f_{j}:=\left\Vert w^{\left( j\right) }\right\Vert \cdot \left\Vert
v^{\left( j\right) }\right\Vert.  \label{normalization}
\end{equation}

In (\ref{normalization}), for $\Lambda_j$ single real, $\left\Vert w^{\left( j\right) }\right\Vert$ denotes the real induced norm (\ref{rinrv}) of the real row vector $w^{(j)}$. On the other hand, for $\Lambda _{j}$ single complex, $\left\Vert w^{\left( j\right) }\right\Vert$ denotes the complex induced norm (\ref{cinrv}) of the complex row vector $w^{(j)}$. 
\subsection{The numbers $g_j(t,u)$ and $g_j(t)$, the vector $\widehat{\Theta}_j(t,u)$ and the matrix $\widehat{\Theta}_j(t)$.}
\label{thenumbers}
In the following, in both cases for $\Lambda_j$, single real or single complex, we write the norm of the real vector $Q_{j}(t)u$, $u\in
\mathbb{R}^n$, as 
\begin{equation}
\Vert Q_{j}(t)u\Vert =f_j\left\vert \widehat{w}^{(j)} u\right\vert g_j(t,u)
\label{normQJ0ureal}
\end{equation}
for a suitable number $g_j(t,u)$. Analogously, we write the real induced norm (\ref{rin}) of the real 
matrix $Q_{j}(t)$ as 
\begin{equation}
\Vert Q_{j}(t)\Vert =f_jg_j(t)  \label{normQJ0}
\end{equation}
for a suitable number $g_j(t)$.

Next two lemmas state what $g_j(t,u)$ and $g_j(t)$ are.

\begin{lemma} \label{Lemmareal}
Let $j\in\{1,\ldots,q\}$. If $\Lambda_j$ is single real, then 
\begin{equation*}
g_j(t,u)=1,\ u\in\mathbb{R}^n,\text{\ and\ \ }g_j(t)=1. 
\end{equation*}
\end{lemma}

\begin{proof}
For $u\in\mathbb{R}^n$, by (\ref{Qj0(t)ureal}) we have
\begin{equation*}
Q_{j}(t)u =v^{(j)}w^{(j)}u=\left(w^{(j)}u\right)v^{(j)}=f_j\left(\widehat{w}^{(j)}u\right)\widehat{v}^{(j)}
\end{equation*}
and then 
\begin{equation}
\left\Vert Q_{j}(t)u\right\Vert =\left\vert w
^{(j)}u\right\vert \left\Vert v^{(j)}\right\Vert=f_j\left\vert \widehat{w}
^{(j)}u\right\vert.  \label{ggjj}
\end{equation}
By the second equality in (\ref{ggjj}), we obtain $g_j(t,u)=1$. By the first equality, we obtain 
\begin{equation*}
\left\Vert Q_{j}(t)\right\Vert=\left\Vert w
^{(j)}\right\Vert\left\Vert v^{(j)}\right\Vert=f_j
\end{equation*}
and then $g_j(t)=1$.
\end{proof}

\begin{lemma}\label{Lemmacomplex}
\label{gj2} Let $j\in\{1,\ldots,q\}$. If $\Lambda _{j}$ is single complex, then 
\begin{equation*}
g_j(t,u)=2\left\Vert \widehat{\Theta}_j(t,u)\right\Vert,\ u\in\mathbb{R}^n,\text{\ and\ \ }g_j(t)=2\left\Vert \widehat{\Theta}_j(t)\right\Vert, 
\end{equation*}
where 
\begin{eqnarray*}
&& \widehat{\Theta }_j\left( t,u\right) =\left( \left\vert \widehat{v}
_{k}^{\left( j\right) }\right\vert \cos \left( \omega _{j}t+\alpha
_{jk}+\gamma_j\left( u\right) \right) \right) _{k=1,\ldots
,n}\in \mathbb{R}^{n}\\
&& \widehat{\Theta }_j\left( t\right) =\left[\left\vert \widehat{v}
_{k}^{\left( j\right) }\right\vert \left\vert \widehat{w}_{l}^{\left(
j\right) }\right\vert \cos \left( \omega _{j}t+\alpha _{jk}+\beta _{jl}\right) \right]_{k,l=1,\ldots ,n}\in \mathbb{
R}^{n\times n} \\
\end{eqnarray*}
with $\alpha_{jk}$, $\beta_{jl}$ and $\gamma_j(u)$ defined by 
\begin{eqnarray*}
&&\widehat{v}_k^{\left( j\right) }=\left\vert \widehat{v}_{k}^{\left(
j\right) }\right\vert \mathrm{e}^{\mathrm{i}\alpha _{jk}},\ k\in\{1,\ldots ,n\},\\
&&\widehat{w}_l^{\left( j\right) }=\left\vert \widehat{w}_{l}^{\left(
j\right) }\right\vert \mathrm{e}^{\mathrm{i}\beta _{jl}},\ l\in\{1,\ldots ,n\}, \\
&&\widehat{w}^{\left( j\right) }u=\left\vert\widehat{w}^{\left( j\right)
}u\right\vert\mathrm{e}^{\mathrm{i}\gamma _j\left( u\right) },
\end{eqnarray*}
i.e., they are, respectively, the angles in the polar forms of the components of the complex vector $\widehat{
v}^{(j)}$, the angles in the polar forms of the components of the complex vector $\widehat{w}^{(j)}$ and  the angle in the polar form of the complex number $
\widehat{w}^{\left( j\right) }u$. 

Moreover, we have
\begin{equation} \label{A==}
\left\Vert \widehat{\Theta }_j\left( t\right) \right\Vert =\max\limits_{\substack{ u\in \mathbb{R}^{n} 
\\ \Vert u\Vert =1}} \left\vert \widehat{w}^{\left( j\right) }u\right\vert\left\Vert \widehat{\Theta }_j \left( t,u\right)
\right\Vert.
\end{equation}
\end{lemma}
In this lemma, $\left\Vert \widehat{\Theta}_j(t)\right\Vert$  is the real induced norm (\ref{rin}) of the real matrix  $\widehat{\Theta}_j(t)$.

\begin{proof}
By (\ref{Qj0(t)ucomplex}), we have
	\begin{equation}
		Q_{j}(t)=2\mathrm{Re}\left( \mathrm{e}^{\mathrm{i}\omega
			_{j}t}v^{(j)}w^{(j)}\right)=2f_j\mathrm{Re}\left( \mathrm{e}^{\mathrm{i}\omega _{j}t}
		\widehat{v}^{(j)}\widehat{w}^{(j)}\right).  \label{Q2}
	\end{equation}
	
Thus, for $u\in\mathbb{R}^n$, we have
\begin{eqnarray*}
&&Q_{j}(t)u=2f_j\mathrm{Re}\left( \mathrm{e}^{\mathrm{i}\omega _{j}t}
\widehat{v}^{(j)}\widehat{w}^{(j)}\right)u=2f_j\mathrm{Re}\left( \mathrm{e}^{
\mathrm{i}\omega _{j}t}\widehat{v}^{(j)}\widehat{w}^{(j)}u\right) \\
&&=2f_j\mathrm{Re}\left( \mathrm{e}^{
	\mathrm{i}\omega _{j}t}\left( \widehat{w}^{\left( j\right) }u\right)\widehat{v}^{(j)}\right)=2f_j\left\vert \widehat{w}^{\left( j\right) }u\right\vert \mathrm{Re}
\left( \mathrm{e}^{\mathrm{i}\omega _{j}t}\mathrm{e}^{\mathrm{i}\gamma_j\left( u\right) }\widehat{v}^{(j)}\right)
\end{eqnarray*}
and 
\begin{equation*}
\mathrm{Re}\left( \mathrm{e}^{\mathrm{i}\omega _{j}t}\mathrm{
e}^{\mathrm{i}\gamma_j\left( u\right) }\widehat{v}^{(j)}\right)=\widehat{\Theta }_j\left(
t,u\right). 
\end{equation*}
We obtain
\begin{equation}
	Q_{j}(t)u=2f_j\left\vert \widehat{w}^{\left( j\right) }u\right\vert\widehat{\Theta }_j\left(
	t,u\right) \label{Qj0u(t)}
\end{equation}
and then $g_j(t,u)=2\left\Vert \widehat{\Theta}_j(t,u)\right\Vert$.

Moreover, in (\ref{Q2}) we have 
\begin{equation*}
\mathrm{Re}\left( \mathrm{e}^{\mathrm{i}\omega _{j}t}\widehat{v}^{(j)}
\widehat{w}^{(j)}\right)=\widehat{\Theta}_j(t).
\end{equation*}
Therefore, we obtain
\begin{equation}
	Q_{j}(t)=2f_j\widehat{\Theta }_j\left(
	t\right) \label{Qj0(t)}
\end{equation}
and then $g_j(t)=2\left\Vert \widehat{\Theta}_j(t)\right\Vert$.

Finally, by (\ref{Qj0(t)}) and (\ref{Qj0u(t)}), we have
\begin{equation*} 
\left\Vert \widehat{\Theta }_j\left( t\right) \right\Vert =\frac{\left\Vert
Q_{j}(t)\right\Vert }{2f_{j}}=\frac{1}{2f_j}\max\limits_{\substack{ u\in \mathbb{R}^{n} 
\\ \Vert u\Vert =1}}\left\Vert Q_{j}(t)u\right\Vert=\max\limits_{\substack{ u\in \mathbb{R}^{n} 
\\ \Vert u\Vert =1}}\left\vert \widehat{w}^{\left( j\right) }u\right\vert\left\Vert \widehat{\Theta }_j \left( t,u\right)
\right\Vert.
\end{equation*}
\end{proof}

\begin{remark}
\label{remTheta} \quad

\begin{itemize}
\item[1.] The functions $t\mapsto \widehat{\Theta}_j(t,u)$ and $t\mapsto 
\widehat{\Theta}_j(t)$ are periodic functions of period $\frac{2\pi}{
\omega_j}$ and the functions $t\mapsto \Vert \widehat{\Theta}_j(t,u)\Vert $ and $t\mapsto 
\Vert \widehat{\Theta}_j(t)\Vert $ are periodic functions of period $\frac{\pi}{
\omega_j}$. 

\item[2.] When the vector norm is a $p$-norm, we have 
\begin{equation*}
\left\Vert \widehat{\Theta }_j\left( t,u\right) \right\Vert
_{p}=\left(\sum\limits_{k=1}^{n}\left\vert \widehat{v}_{k}^{\left( j\right)
}\right\vert^p \left\vert \cos \left( \omega _{j}t+\alpha _{jk}+\gamma_j\left( u\right) \right)
\right\vert^p\right)^{\frac{1}{p}} \leq \left\Vert \widehat{v}^{\left(
j\right) }\right\Vert_p =1.
\end{equation*}
\item[3.] When the vector norm is a $p$-norm, we have $\left\Vert \widehat{\Theta }_j\left( t\right)
\right\Vert_p\leq 1$. This follows by (\ref{A==}),
$$
\max\limits_{\substack{ u\in \mathbb{R}^{n} 
\\ \Vert u\Vert =1}}\left\vert \widehat{w}^{\left( j\right) }u\right\vert=\left\Vert \widehat{w}^{\left( j\right) }\right\Vert_{\mathbb{R}}\leq \left\Vert \widehat{w}^{\left( j\right) }\right\Vert=1\text{\ \ \ (recall(\ref{r<=c}))}
$$
and the previous point 2.
\end{itemize}
\end{remark}

\section{The norms $\left\Vert Q_j(t)u \right\Vert$ and $\left\Vert Q_j(t) \right\Vert$ for the Euclidean norm}

In this section, we consider the case where $\Lambda_{j}$ is single complex and the vector norm $\Vert \ \cdot \ \Vert$ is the Euclidean norm. Consequently, the vectors $\widehat{v}^{(j)}$ and $\widehat{w}^{(j)}$ in (\ref{normalization}) are defined using the Euclidean norm. Our aim is to derive expressions for the norms $\left\Vert \widehat{\Theta}_j(t,u) \right\Vert_2$ and $\left\Vert \widehat{\Theta}_j(t) \right\Vert_2$ appearing in
\begin{equation}
	\left\Vert Q_{j}(t)u\right\Vert_2=2f_j\left\vert \widehat{w}^{(j)}u\right\vert \left\Vert \widehat{\Theta}_{j}(t,u)\right\Vert_2  \label{normQJ0u}
\end{equation}
and
\begin{equation*}
	\left\Vert Q_{j}(t)\right\Vert_2=2f_j\left\Vert \widehat{\Theta}_{j}(t)\right\Vert_2
\end{equation*}
(see (\ref{normQJ0ureal}) and (\ref{normQJ0}) and Lemma \ref{Lemmacomplex}).

\subsection{Notations}\label{notations}

In the following, for a complex row or column vector $z$, $z^T$ denotes the transpose of $z$, representing only transposition without conjugation. The conjugate transpose of $z$ is denoted by $z^H$.

Before continuing with our analysis, we introduce some important quantities. They are:

\begin{itemize}
	
	\item $\gamma_j(u)$, the angle in the polar form of the complex number $\widehat{w}^{(j)}u$, $u\in\mathbb{R}^n$, already introduced in Lemma \ref{gj2};
	\item $V_j$, the modulus of the complex number $\left(\widehat{v}^{\left( j\right) }\right) ^{T}\widehat{v}^{\left( j\right)}$;
		\item  $\delta_j$, the angle in the polar form of the complex number $\left(\widehat{v}^{\left( j\right) }\right) ^{T}\widehat{v}^{\left( j\right)}$;
	\item $W_j$, the modulus of the complex number $\widehat{w}^{(j)}\left(\widehat{w}^{(j)}\right)^T$. 
\end{itemize}
\begin{remark}
	\label{VW}
	\quad
	\begin{itemize}
		\item [1.]  Observe that  $\left(\widehat{v}^{\left( j\right) }\right) ^{H}
		\widehat{v}^{\left( j\right)}=1$ and $\widehat{w}^{(j)}\left(\widehat{w}^{(j)}\right)^H=1$. On the other hand, $\left(\widehat{v}^{\left( j\right) }\right) ^{T}
		\widehat{v}^{\left( j\right)}$ and $\widehat{w}^{(j)}\left(\widehat{w}^{(j)}\right)^T$ are complex numbers and $V_j$ and $W_j$ are their moduli, respectively.
		\item [2.] By the Cauchy-Schwarz inequality we have 
		\begin{equation*}
			V_j=\left\vert \left( \widehat{v}^{\left( j\right) }\right) ^{T}\widehat{v}
			^{\left( j\right) }\right\vert = \left\vert \left( \overline{\widehat{v}
				^{\left( j\right) }}\right) ^{H}\widehat{v}^{\left( j\right)
			}\right\vert\leq 1. 
		\end{equation*}
		Moreover, $V_j<1$ holds, otherwise $\widehat{v}^{\left( j\right) }$ and $
		\overline{\widehat{v}^{\left( j\right) }}$ were linearly dependent, which is
		impossible since $v^{\left( j\right) }$ and $\overline{v
			^{\left( j\right) }}$ are columns of the non-singular matrix $V$, whose columns constitute the Jordan basis for $A$. Therefore, we have $V_j\in[0,1)$.
		\item [3.]  By the Cauchy-Schwarz inequality and Remark \ref{22}, or by the definition of induced norm $\left\Vert w^{(j)}\right\Vert$ and Remark \ref{22}, we have 
		\begin{equation*}
			W_j=\left\vert \widehat{w}^{\left( j\right) } \left(\widehat{w}^{\left(
				j\right) }\right)^T\right\vert = \left\vert \left(\left(\overline{\widehat{w}^{\left( j\right) }}\right)^T\right)^H \left(\widehat{w}^{\left(
				j\right) }\right)^T
			\right\vert\leq 1. 
		\end{equation*}
		Moreover, $W_j<1$ holds, otherwise $\widehat{w}^{\left( j\right) }$ and $
		\overline{\widehat{w}^{\left( j\right) }}$ were linearly dependent, which is
		impossible since $w^{\left( j\right) }$ and $\overline{w
			^{\left( j\right) }}$ are rows of the non-singular matrix $W=V^{-1}$. Therefore, we have $W_j\in[0,1)$.
	\end{itemize}
\end{remark}

\subsection{The set $\mathcal{W}_j$ and the ellipse $\mathcal{E}_j$.}

The complex number $ \widehat{w}^{(j)}u $, $u\in\mathbb{R}^n$, plays a crucial role in the vector $\widehat{\Theta}_j(t,u)$ and the matrix $\widehat{\Theta}_j(t)$: $ \widehat{\Theta}_j(t,u) $ is defined in Lemma \ref{Lemmacomplex} by means of  the angle $ \gamma_j(u) $ in the polar form of $ \widehat{w}^{(j)}u $ and the modulus $\vert  \widehat{w}^{(j)}u\vert$ appears in (\ref{A==}). 

For this reason, now we study the set 
\begin{equation*}
\mathcal{W}_j:=\{\widehat{w}^{(j)}u:u\in\mathbb{R}^n\text{\ and\ }\left\Vert
u\right\Vert_2=1\}  \label{set}
\end{equation*}
of the complex numbers $\widehat{w}^{(j)}u$, by varying the unit real vector $u$.

Next proposition describes the set $\mathcal{W}_j$ and states the real induced norm 
\begin{equation*}
	\left\Vert \widehat{w}^{\left( j\right) }\right\Vert_\mathbb{R}=\max\limits_{\substack{ u\in \mathbb{R}^{n}  \\ \Vert u\Vert_2 =1
			}}\left\vert w^{\left( j\right) }u\right\vert  
\end{equation*}
of the complex row vector $\widehat{w}^{(j)}$, whose complex induced norm is equal to $1$.

\begin{proposition}
\label{euclidean} Let $j\in\{1,\ldots,q\}$ such
that $\Lambda _{j}$ is single complex. Consider, in the complex plane, the ellipse $\mathcal{E}_j$ centered at the origin with semi-axes the vectors $\sigma_j
\alpha^{(j)}$ and $\mu_j\beta^{(j)}$, where $\sigma_j\geq \mu_j\geq 0$ are
the singular values of the $2\times n$ real matrix 
\begin{equation*}
R_j:=\left[
\begin{array}{c}
\mathrm{Re}(\widehat{w}^{(j)}) \\ 
\mathrm{Im}(\widehat{w}^{(j)})
\end{array}
\right]  \label{matrix}
\end{equation*}
and $\alpha^{(j)},\beta^{(j)}\in\mathbb{R}^2$ are the corresponding left
singular vectors. We have
\begin{equation}
\sigma_j=\sqrt{\frac{1+W_j}{2}}\text{\ \ and\ \ }
\mu_j=\sqrt{\frac{1-W_j
}{2}}.  \label{wR}
\end{equation}

If $n=2$, the set $\mathcal{W}_j$ is the ellipse $\mathcal{E}_j$. If $n>2$, the set $\mathcal{W}_j$ is the union of the ellipse $\mathcal{E}_j$ and its interior. Moreover, we have
\begin{equation*}
\left\Vert \widehat{w}^{\left( j\right) }\right\Vert_\mathbb{R}=\sigma_j=\sqrt{\frac{1+W_j}{2}}. \label{maxnormR}
\end{equation*}

\end{proposition}

\begin{proof}
We begin by proving (\ref{wR}). The singular values $\sigma_j$ and $\mu_j$ of the matrix $R_j$ are the square roots of the largest
eigenvalue and the smallest eigenvalue, respectively, of the matrix 
\begin{equation*}
R_jR_j^T=\left[
\begin{array}{c}
\mathrm{Re}(\widehat{w}^{(j)}) \\ 
\mathrm{Im}(\widehat{w}^{(j)})
\end{array}
\right]\left[\mathrm{Re}(\widehat{w}^{(j)})^T\ \mathrm{Im}\widehat{w}w^{(j)})^T\right]=\left[ 
\begin{array}{cc}
a_j^2 & p_j \\ 
p_j & b_j^2
\end{array}
\right], 
\end{equation*}
where 
\begin{equation*}
a_j=\left\Vert \mathrm{Re}\left( \widehat{w}^{(j)}\right) ^{T}\right\Vert
_{2},\ b_j=\left\Vert \mathrm{Im}\left( \widehat{w}^{(j)}\right)
^{T}\right\Vert _{2}\text{\ \ and\ \ }p_j=\mathrm{Re}\left( \widehat{w}
^{(j)}\right) \mathrm{Im}\left( \widehat{w}^{(j)}\right) ^{T}.
\end{equation*}
Therefore, we obtain 
\begin{equation*}
\sigma_j^2=\frac{1+\sqrt{\left(a_j^2-b_j^2\right)^2+4p_j^2}}{2}\text{\ and\ }
\mu_j^2=\frac{1-\sqrt{\left(a_j^2-b_j^2\right)^2+4p_j^2}}{2}.
\end{equation*}
by observing that $a_j^2+b_j^2=\left\Vert \left(\widehat{w}^{(j)}\right)^T\right\Vert_2^2=\left\Vert \widehat{w}^{(j)}\right\Vert_2^2=1$
. Since 
\begin{equation*}
\widehat{w}^{(j)}\left(\widehat{w}^{(j)}\right)^T=a_j^2-b_j^2+
\mathrm{i}\cdot 2p_j, 
\end{equation*}
we obtain
\begin{equation*}
W_j=\sqrt{\left(a_j^2-b_j^2\right)^2+4p_j^2} 
\end{equation*}
and (\ref{wR}) follows.

Now, we pass to describing the set $\mathcal{W}_j$ in terms of the ellipse $\mathcal{E}_j$. By looking at a complex number as the pair constituted by the real and
imaginary parts, the complex number $\widehat{w}^{(j)}u$, $u\in\mathbb{R}^n$, is $R_ju$. By
using the Singular Value Decomposition of $R_j$, we see that
\begin{equation*}
R_ju=\sigma_jc_j(u)\alpha^{(j)}+\mu_jd_j(u)\beta^{(j)}
\end{equation*}
where $c_j(u)$ and $d_j(u)$ are the components of $u$ along the first and second right singular vectors of $R_j$, i.e., the right singular vectors corresponding to the singular values $\sigma_j$ and $\mu_j$, respectively. Thus,
\begin{equation}
\mathcal{W}_j=\{\sigma_jc\alpha^{(j)}+\mu_jd\beta^{(j)}: (c,d)\in\mathbb{R}^2\text{\ and\ }
c^2+d^2=1\}  \label{setn=2}
\end{equation}
when $n=2$ and
\begin{equation}
\mathcal{W}_j=\{\sigma_jc\alpha^{(j)}+\mu_jd\beta^{(j)}: (c,d)\in\mathbb{R}^2\text{\ and\ }
c^2+d^2\leq 1\}   \label{setn>2}
\end{equation}
when $n>2$.

In the orthonormal basis of $\mathbb{R}^2$ given by the vector $\alpha^{(j)}$
and $\beta^{(j)}$, the coordinates $x=\sigma_jc$ and $y=\mu_jd$ of the
elements of the set (\ref{setn=2}) satisfy the equation 
\begin{equation*}
\frac{x^2}{\sigma_j^2}+\frac{y^2}{\mu_j^2}=1, 
\end{equation*}
which is the equation of an ellipse centered at the origin with lengths of
the semi-axes $\sigma_j$ and $\mu_j$. Thus, we have $\mathcal{W}_j=\mathcal{E}_j$. On the other hand, the coordinates $x=\sigma_jc$ and $y=\mu_jd$ of the
elements of the set (\ref{setn>2}) satisfy the inequality 
\begin{equation*}
\frac{x^2}{\sigma_j^2}+\frac{y^2}{\mu_j^2}\leq 1. 
\end{equation*}
Thus, $\mathcal{W}_j$ is the union of the ellipse $\mathcal{E}_j$ and its interior.

Finally, we have
\begin{eqnarray*}
\left\Vert \widehat{w}^{\left( j\right) }\right\Vert_\mathbb{R}=\max\limits_{\substack{ u\in \mathbb{R}^{n} 
\\ \Vert u\Vert_2 =1}}\left\vert \widehat{w}^{\left( j\right) }u\right\vert=\max\limits_{\substack{ u\in \mathbb{R}^{n} 
\\ \Vert u\Vert_2 =1}}\left\Vert R_ju\right\Vert_2=\left\Vert R_j\right\Vert_2=\sigma_j.
\end{eqnarray*}
\end{proof}

\begin{remark}\label{remsigma}
Since $W_j<1$ (see point 3 in Remark \ref{VW}), we have 
\begin{eqnarray*}
\left\Vert \widehat{w}^{\left( j\right) }\right\Vert_\mathbb{R}=\sqrt{\frac{1+W_j}{2}}<1
\end{eqnarray*}
for the real induced norm of the row vector $\widehat{w}^{(j)}$. On the other hand, $\left\Vert \widehat{w}^{(j)} \right\Vert_2=1$
for the complex induced norm of $\widehat{w}^{(j)}$. This confirms the inequality (\ref{r<=c}). 
\end{remark}

\subsection{$w^{(j)}u$ and the projection onto the eigenspaces of $\lambda_j$ and $\overline{\lambda}_j$} \label{Pe}

The projection of $u\in\mathbb{R}^n$ onto the sum of the eigenspaces corresponding to $\lambda_j$ and $\overline{\lambda}_j$ is
$$
(w^{(j)}u)v^{(j)}+(\overline{w^{(j)}}u)\ \overline{v^{(j)}},
$$
where the first and the second terms are the projections onto the eigenspace corresponding to $\lambda_j$ and $\overline{\lambda_j}$, respectively: observe that 
$w^{(j)}u$ and $\overline{w^{(j)}}u$ are the components of $u$ along $v^{(j)}$ and $\overline{v^{(j)}}$, respectively, in the Jordan basis (since $w^{(j)}$ and $\overline{w^{(j)}}$ are the rows of $W=V^{-1}$ corresponding to the columns $v^{(j)}$ and $\overline{v^{(j)}}$, respectively, of $V$).

Since $\widehat{v}^{(j)}$ and $\overline{\widehat{v}^{(j)}}$ are linearly independent, the projection onto the sum of the eigenspaces is zero if and only if the projection onto the eigenspace corresponding to $\lambda_j$ is zero, i.e., $\widehat{w}^{(j)}u=0$.

The Euclidean norm of the projection of $u$ onto the eigenspace corresponding to $\lambda_j$ is $f_j\left\vert \widehat{w}^{(j)}u \right\vert$, equal to the norm of the projection onto the eigenspace corresponding to $\overline{\lambda_j}$. The Euclidean norm of the projection of $u$ onto the sum of the eigenspaces is given in Remark \ref{project2} below.

\subsection{$w^{(j)}u$ and the right singular vectors of the matrix $R_j$}
Vectors in $\mathbb{R}^n$ can be  expressed in the orthonormal basis of $\mathbb{R}^n$ of the right singular vectors of the matrix $R_j$ of Proposition \ref{euclidean}.

As we saw in the proof of Proposition \ref{euclidean}, for a vector $u\in\mathbb{R}^n$  we have 
\begin{equation}
\widehat{w}^{(j)}u=a+\mathrm{i}b \label{whatu}
\end{equation}
with
\begin{equation}
	(a,b)=R_ju=\sigma_jc_j(u)\alpha^{(j)}+\mu_jd_j(u)\beta^{(j)}, \label{whatu2}
\end{equation}
where: $\sigma_j$ and $\mu_j$ are the two singular values of $R_j$; $c_j(u)$ and $d_j(u)$ are the components of $u$ along the first and second, respectively, right singular vectors of $R_j$; and $\alpha^{(j)}$ and $\beta^{(j)}$ are the left singular vectors of $R_j$. 

Therefore, $u$ has projection zero onto the sum of the eigenspaces corresponding to $\lambda_j$ and $\overline{\lambda}_j$, i.e. $\widehat{w}^{(j)}u=0$ by recalling the previous Subsection \ref{Pe}, if and only if  $u$ is spanned by the last $n-2$ right singular vectors of $R_j$.

By (\ref{whatu})-(\ref{whatu2}), we have
\begin{equation}
\left\vert \widehat{w}^{(j)}u\right\vert=\left\Vert R_ju\right\Vert_2=\sqrt{\sigma_j^2c_j(u)^2+\mu_j^2d_j(u)^2}. \label{size}
\end{equation}

\subsection{Further notations}\label{notation2}

In addition to the quantities introduced in Subsection \ref{notations}, here we introduce a few others. They are:
\begin{itemize}

	\item $\theta_j$, the angle in the complex plane between the real axis and the major axis $\sigma_j\alpha^{(j)}$ of the ellipse $\mathcal{E}_j$ of Proposition \ref{euclidean}; 
	\item
	\begin{equation}
			\Omega_j(t):=2(\omega_{j}t+\theta_j)+\delta_{j}, \label{xt}
		\end{equation}
		where $\delta_j$ is introduced in Subsection \ref{notations} and (recall) $\omega_j$ is the imaginary part of the eigenvalue $\lambda_j$, which is positive;
	\item for $u\in\mathbb{R}^n$, 
	\begin{equation}
			\Delta_j(u):=2(\gamma_{j}(u)-\theta_{j}), \label{xt2}
		\end{equation}
	where $\gamma_j(u)$ is introduced in Subsection \ref{notations}.
\end{itemize}

\subsection{The norm $\left\Vert \widehat{\Theta }_{j}\left(
t,u\right) \right\Vert _{2}$}

Next proposition provides an expression for  the norm $\left\Vert \widehat{\Theta }_{j}\left(
t,u\right) \right\Vert _{2}$.

\begin{proposition}
\label{oscillation} Let $j\in\{1,\ldots,q\}$
such that $\Lambda _{j}$ is single complex. For $u\in\mathbb{R}^n$, we have 
\begin{equation}
\left\Vert \widehat{\Theta }_{j}\left( t,u\right) \right\Vert _{2}=\sqrt{\frac{1+V_{j}\cos \left( \Omega_j(t)+\Delta_j(u)\right)}{2}}, \label{Theta1}
\end{equation}
where $\Omega_j(t)$ and $\Delta_j(u)$ are defined in (\ref{xt}) and (\ref{xt2}).
\end{proposition}

\begin{proof}
Let 
\begin{equation*}
\widehat{w}^{(j)}u=\vert \widehat{w}^{(j)}u\vert \mathrm{e}^{\mathrm{i}\gamma_j(u)}\text{\ \ and\ \ }\left( \widehat{v}^{\left( j\right) }\right) ^{T}\widehat{v}^{\left(
j\right) }=V_{j}\mathrm{e}^{\mathrm{i}\delta _{j}}
\end{equation*}
be the polar forms of  $\widehat{w}^{(j)}u$ and $\left( \widehat{v}^{\left( j\right)
}\right) ^{T}\widehat{v}^{\left( j\right) }$ (recall the notations introduced in Subsection \ref{notations}). By (\ref{Qj0(t)ucomplex}), we have 
\begin{eqnarray}
&&\left\Vert Q_{j}(t)u\right\Vert_{2}^{2}=\left\Vert \mathrm{e}^{\mathrm{i}
\omega_j t}v^{(j)}w^{(j)}u+\overline{\mathrm{e}^{\mathrm{i}\omega_j
t}v^{(j)}w^{(j)}u}\right\Vert _{2}^{2} \notag \\
&&=\left\Vert \mathrm{e}^{\mathrm{i}\omega_j t}v^{(j)}w^{(j)}u\right\Vert
_{2}^{2}+2\mathrm{Re}\left( \left( \overline{\mathrm{e}^{\mathrm{i}\omega_j
t}v^{(j)}w^{(j)}u}\right) ^{H}\mathrm{e}^{\mathrm{i}\omega
_{j}}v^{(j)}w^{(j)}u\right) \notag \\
&&\quad +\left\Vert \overline{\mathrm{e}^{\mathrm{i}\omega_j t}v^{(j)}w^{(j)}u}\right\Vert
_{2}^{2}\notag \\
&&=2\left\Vert \mathrm{e}^{\mathrm{i}\omega_j t}v^{(j)}w^{(j)}u\right\Vert
_{2}^{2}+2\mathrm{Re}\left( \left( \mathrm{e}^{\mathrm{i}\omega_j
t}v^{(j)}w^{(j)}u\right) ^{T}\mathrm{e}^{\mathrm{i}\omega_j
t}v^{(j)}w^{(j)}u\right)\notag \\
&&=2\left\Vert v^{(j)}\right\Vert _{2}^{2}\left\vert
w^{(j)}u\right\vert ^{2}+2\mathrm{Re}\left( \mathrm{e}^{2\mathrm{i}\omega_j
t}\left(w^{(j)}u\right)
^{2}\left(v^{(j)}\right)^{T}v^{(j)} \right) \notag \\
&&=2f_j^2\vert \widehat{w}^{(j)}u\vert^2+2\mathrm{Re}\left(f_j^2 \mathrm{e}^{2\mathrm{i}\omega_j
	t}\left(\widehat{w}^{(j)}u\right)
^{2}\left(\widehat{v}^{(j)}\right)^{T}\widehat{v}^{(j)} \right)\notag \\
&&=2f_j^2\vert \widehat{w}^{(j)}u\vert^2+2\mathrm{Re}\left(f_j^2\vert \widehat{w}^{(j)}u\vert^2V_j\mathrm{e}^{2\mathrm{i}\omega_jt}\mathrm{e}^{2\mathrm{i}\gamma_j(u)}\mathrm{e}^{\mathrm{i}\delta_j}\right)\notag \\
&&=2f_j^{2}\left\vert \widehat{w}^{(j)}u\right\vert ^{2}\left( 1+V_j\cos
\left( 2\omega_j t+2\gamma_j (u)+\delta_j \right) \right)\label{expQ} \\
&&=2f_j^{2}\left\vert \widehat{w}^{(j)}u\right\vert ^{2}\left( 1+V_j\cos \left( \Omega_j(t)+\Delta_j(u)\right)\right).  \notag
\end{eqnarray}
Now, (\ref{Theta1}) follows by (\ref{normQJ0u}).
\end{proof}

Although in (\ref{Theta1}) the argument of $\cos$, namely
\begin{equation}
\Omega_j(t)+\Delta_j(u)=2(\omega_jt+\gamma_j(u))+\delta_j, \label{Omegajplus}
\end{equation}
does not involve $\theta_j$, it is convenient for the subsequent developments in the proof of Proposition \ref{oscillation2} to split this argument as the sum of $\Omega_j(t)$ and $\Delta_j(u)$, with both of them involving $\theta_j$.

The following theorem is a direct consequence of the previous Proposition \ref{oscillation}.

\begin{theorem}
The $\frac{\pi}{\omega_j}$-periodic function $\left\Vert \widehat{\Theta }_{j}\left( t,u\right) \right\Vert _{2}$ of $t$ oscillates monotonically between the extreme values
\begin{eqnarray}
&&\max\limits_{t\in\mathbb{R}} \left\Vert \widehat{\Theta }_{j}\left( t,u\right)
\right\Vert _{2}=\sqrt{\frac{1+V_{j}}{2}}\label{maxtheta}
\end{eqnarray}
and
\begin{eqnarray}
	&&\min\limits_{t\in\mathbb{R}} \left\Vert \widehat{\Theta }_{j}\left(
	t,u\right) \right\Vert _{2}= \sqrt{\frac{1-V_{j}}{2}}. \label{mintheta}
\end{eqnarray}
The maximum value is achieved if $\Omega_j(t)+\Delta_j(u)$ is an even multiple of $\pi$ and the minimum value is achieved if $\Omega_j(t)+\Delta_j(u)$ is an odd multiple of $\pi$.
\end{theorem}

\begin{proof}
	The theorem follows from the fact that the argument (\ref{Omegajplus}) of $\cos$ in (\ref{Theta1}), considered as a function of $t$, is a strictly increasing continuous bijection $\mathbb{R}\rightarrow \mathbb{R}$, i.e., an orientation preseving homeomorphism of $\mathbb{R}$. 
\end{proof}

Observe that, when $V_j=0$, we have
$$
\left\Vert \widehat{\Theta }_{j}\left( t,u\right) \right\Vert _{2}=\sqrt{\frac{1}{2}},
$$
a constant function of $t$.

\begin{remark} \label{project2}
Recall Subsection \ref{Pe}, in particular the last paragraph. By using $t=0$ in (\ref{expQ}) in the proof of Proposition \ref{oscillation}, we obtain 
$$
\left\Vert Q_{j}(0)u\right\Vert_2=f_j\left\vert \widehat{w}^{(j)}u\right\vert\sqrt{2\left(1+V_j\cos \left(\delta_j +2\gamma_j (u)\right)\right)}.
$$
This is the norm of the projection of $u\in\mathbb{R}^n$ onto the sum of the eigenspaces corresponding to  $\lambda_j$ and $\overline{\lambda_j}$: $Q_j(0)$ is the projection onto the sum of the eigenspaces corresponding to  $\lambda_j$ and $\overline{\lambda_j}$, see (\ref{Qj0(t)ucomplex}). The norm of the projection on each eigenspace is $f_j\left\vert \widehat{w}^{(j)}u\right\vert$.
\end{remark}

\subsection{The functions $f_{V_jW_j}$, $f_{V_jW_j}^{\max}$ and $f_{V_jW_j}^{\min}$} \label{fVWsect}

For the  analysis of  the norm $\left\Vert \widehat{\Theta }_{j}\left(
t\right) \right\Vert _{2}$, as well as for the analysis of the oscillating terms in Section \ref{OSFOT}, we introduce the functions $f_{V_jW_j}$, $f_{V_jW_j}^{\max}$ and $f_{V_jW_j}^{\min}$.

The function $f_{V_jW_j}:\mathbb{R}^2\rightarrow \mathbb{R}$ is given by
\begin{equation}
	f_{V_jW_j}\left(\alpha,x\right):=\frac{1+V_j\cos\left( x+\alpha\right)}{1-W_j\cos\alpha}
	,\ \left(\alpha,x\right)\in\mathbb{R}^2, \label{fVW0}
\end{equation}
and the $2\pi$-periodic functions $f_{V_jW_j}^{\max},f_{V_jW_j}^{\min}:\mathbb{R}\rightarrow \mathbb{R}$ are given by
\begin{equation}
	f_{V_jW_j}^{\max}(x)=\max\limits_{\alpha\in\mathbb{R}} f_{V_jW_j}\left(\alpha,x\right),\ x\in\mathbb{R},  \label{fmax0}
\end{equation}
and
\begin{equation}
	f_{V_jW_j}^{\min}(x)=\min\limits_{\alpha\in\mathbb{R}} f_{V_jW_j}\left(\alpha,x\right),\ x\in\mathbb{R}.  \label{fmin0}
\end{equation}

The functions $f_{V_jW_j}$, $f_{V_jW_j}^{\max}$ and $f_{V_jW_j}^{\min}$ are studied in Appendix \ref{AA}.

\subsection{The norm $\left\Vert \widehat{\Theta }_{j}\left(
	t\right) \right\Vert _{2}$}

Next proposition provides an expression for the norm $\left\Vert \widehat{\Theta }_{j}\left(
t\right) \right\Vert _{2}$.

\begin{proposition}
\label{oscillation2} Let $j\in\{1,\ldots,q\}$
such that $\Lambda _{j}$ is single complex. We have
\begin{equation}
\left\Vert \widehat{\Theta }_{j}\left( t\right) \right\Vert _{2}=\sqrt{\frac{1-W_j^2}{4}f_{V_jW_j}^{\max}\left(\Omega_j(t)\right)},  \label{Theta2}
\end{equation}
where $f_{V_jW_j}^{\max}$ is defined in (\ref{fmax0}) and $\Omega_j(t)$ is defined in (\ref{xt}).
\end{proposition}

\begin{proof}
By (\ref{A==}), Proposition \ref{euclidean} and Proposition \ref{oscillation}, we obtain
\begin{equation}
\left\Vert \widehat{\Theta }_j\left( t\right) \right\Vert_2 =\max\limits_{\rho\mathrm{e}^{\mathrm{i}\theta}\in \mathcal{W}_j}\rho\sqrt{\frac{1+V_{j}\cos \left(\Omega_{j}(t)+2(\theta-\theta_j)\right) }{2}},
\label{normtheta}
\end{equation}
where $\mathcal{W}_j$ is the ellipse $\mathcal{E}_j$ of Proposition \ref{euclidean} when $n=2$, and $\mathcal{W}_j$ is the union of $\mathcal{E}_j$ and its interior when $n>2$.
Moreover, $\rho\mathrm{e}^{\mathrm{i}\theta}$ is the polar form of the complex numbers in $\mathcal{W}_j$, over which we take the maximum. We now show that (\ref{Theta2}) holds.

Suppose $n>2$. By using, in the complex plane, cartesian coordinates with axes coinciding with those of the ellipse $\mathcal{E}_j$, the set $\mathcal{W}_j$ is the set of point $(x,y)$ such that
$$
\frac{x^2}{\sigma_j^2}+\frac{y^2}{\mu_j^2}\leq 1,
$$
where $\sigma_j$ and $\mu_j$ are the lengths of the semi-axes of the ellipse. Thus, since the cartesian coordinate of a point of polar form $\rho\mathrm{e}^{\mathrm{i}\theta}$ are
$$
x=\rho\cos(\theta-\theta_j)\text{\ \ and\ \ }y=\rho\sin(\theta-\theta_j),
$$
we obtain
\begin{eqnarray*}
\rho\mathrm{e}^{\mathrm{i}\theta}\in\mathcal{W}_j &\Leftrightarrow& \frac{\rho^2\cos^2(\theta-\theta_j)}{\sigma_j^2}+\frac{\rho^2\sin^2(\theta-\theta_j)}{\mu_j^2}\leq 1\\
&\Leftrightarrow& \rho^2\leq R_j\left(\theta\right)^2,
\end{eqnarray*}
where
\begin{equation}
R_j\left(\theta\right):=\sqrt{\frac{\sigma_j^2\mu_j^2}{\mu^2_j\cos^2(\theta-\theta_j)+\sigma^2_j\sin^2(\theta-\theta_j)}}.  \label{Rtheta}
\end{equation}
Thus, by (\ref{normtheta}) we have
\begin{eqnarray*}
\left\Vert \widehat{\Theta }_j\left( t\right) \right\Vert_2&=&\max\limits_{\theta\in[0,2\pi]}\max\limits_{\rho\in\left[0,R_j\left(\theta\right)\right]}
\rho\sqrt{\frac{1+V_{j}\cos\left(\Omega_{j}(t)+2(\theta-\theta_j)\right) }{2}}\notag\\
&=&\max\limits_{\theta\in[0,2\pi]}R_j\left(\theta\right)\sqrt{\frac{1+V_{j}\cos\left( \Omega_{j}(t)+2(\theta-\theta_j)\right) }{2}}. \label{cos}
\end{eqnarray*}

In the case $n=2$, where $\mathcal{W}_j$ is the set of points $(x,y)$ such that
$$
\frac{x^2}{\sigma_j^2}+\frac{y^2}{\mu_j^2}=1,
$$
we obtain
$$
\rho\mathrm{e}^{\mathrm{i}\theta}\in\mathcal{W}_j\Leftrightarrow  \rho^2= R_j\left(\theta\right)^2
$$
and then
$$
\left\Vert \widehat{\Theta }_j\left( t\right) \right\Vert_2=\max\limits_{\theta\in[0,2\pi]}R_j\left(\theta\right)\sqrt{\frac{1+V_{j}\cos\left(\Omega_{j}(t)+2(\theta-\theta_j)\right) }{2}}
$$
in this case as well.

Recall (\ref{wR}) in Proposition \ref{euclidean}. In (\ref{Rtheta}), we have
$$
\sigma_j^2\mu_j^2=\frac{1-W^2_j}{4}
$$
and, by introducing $\alpha=2(\theta-\theta_j)$, 
\begin{eqnarray*}
&&\mu^2_j\cos^2(\theta-\theta_j)+\sigma^2_j\sin^2(\theta-\theta_j)=\mu^2_j\cos^2\frac{\alpha}{2}+\sigma^2_j\sin^2\frac{\alpha}{2}\\
&&=\mu^2_j\frac{1+\cos\alpha}{2}+\sigma_j^2\frac{1-\cos\alpha}{2}=\frac{\mu^2_j+\sigma^2_j+\left(\mu^2_j-\sigma^2_j\right)\cos\alpha}{2}\\
&&=\frac{1-W_j\cos\alpha}{2}.
\end{eqnarray*}
Therefore, we obtain
$$
R_j\left(\theta\right)=\sqrt{\frac{1-W_j^2}{2\left(1-W_j\cos\alpha\right)}}.
$$
Now, we can write
\begin{eqnarray*}
&&\left\Vert \widehat{\Theta }_j\left( t\right) \right\Vert_2=\max\limits_{\theta\in[0,2\pi]}R_j\left(\theta\right)\sqrt{\frac{1+V_{j}\cos\left(\Omega_{j}(t)+2(\theta-\theta_j)\right)}{2}}\\
&&=\max\limits_{\alpha\in[-2\theta_j,2(2\pi-\theta_j)]}\sqrt{\frac{1-W_j^2}{2\left(1-W_j\cos\alpha\right)}}\cdot\sqrt{\frac{1+V_{j}\cos\left(\Omega_j(t)+\alpha\right)}{2}}\\
&&=\max\limits_{\alpha\in[0,2\pi]}\sqrt{\frac{1-W_j^2}{2\left(1-W_j\cos\alpha\right)}}\cdot\sqrt{\frac{1+V_{j}\cos\left(\Omega_j(t)+\alpha\right)}{2}},
\end{eqnarray*}
i.e. (\ref{Theta2}) by recalling the definition (\ref{fmax0})-(\ref{fVW0}) of $f^{\max}_{V_jW_j}$. 
\end{proof}

The following theorem is a direct consequence of Proposition \ref{oscillation2}.

\begin{theorem}
The $\frac{\pi}{\omega_j}$-periodic function $\left\Vert \widehat{\Theta }_{j}\left( t\right) \right\Vert _{2}$ of $t$ oscillates monotonically between the extreme values
\begin{eqnarray}
	&&\max\limits_{t\in\mathbb{R}} \left\Vert \widehat{\Theta }_{j}\left( t\right)
	\right\Vert _{2}=\sqrt{\frac{(1+V_j)(1+W_{j})}{4}} \label{maxTheta}
\end{eqnarray}
and
\begin{eqnarray}
	&&\min\limits_{t\in\mathbb{R}} \left\Vert \widehat{\Theta }_{j}\left(
	t\right) \right\Vert _{2}= \sqrt{\frac{a_{V_jW_j}}{4}}, \label{minTheta}
\end{eqnarray}
where
\begin{equation}
	a_{V_jW_j}:=\left\{
	\begin{array}{l}
		(1-V_j)(1+W_j)\text{\ if\ }V_j\leq W_j\\
		\\
		(1+V_j)(1-W_j)\text{\ if\ }V_j\geq W_j.
	\end{array}
	\right. \label{aVjWj}
\end{equation}
The maximum value is achieved if $\Omega_j(t)$ is an even multiple of $\pi$ and the minimum value  is achieved if $\Omega_j(t)$ is an odd multiple of $\pi$.
\end{theorem}

\begin{proof}
The function $f_{V_jW_j}^{\max}$ in (\ref{fmax0}) oscillates monotonically between the extreme values
\begin{eqnarray*}
	&&\max\limits_{x\in\mathbb{R}} f_{V_jW_j}^{\max}\left(x\right)
	=\frac{1+V_j}{1-W_j}
\end{eqnarray*}
and
\begin{eqnarray*}
	&&\min\limits_{x\in\mathbb{R}} f_{V_jW_j}^{\max}\left(x\right)
	=\left\{
	\begin{array}{l}
		\frac{1-V_j}{1-W_j}\text{\ if\ }V_j\leq W_j\\
		\\
		\frac{1+V_j}{1+W_j}\text{\ if\ }V_j\geq W_j,
	\end{array}
	\right.\label{minPhi}
\end{eqnarray*}
with the maximum value achieved if $x$ is an even multiple of $\pi$ and the minimum value achieved if $x$ is an odd multiple of $\pi$ (see Appendix \ref{AA}). The theorem now follows from the fact that the argument $\Omega_j(t)$  of $f^{\max}_{V_jW_j}$ in (\ref{Theta2}), considered as a function of $t$, is an orientation preserving homeomorphism of $\mathbb{R}$. 
\end{proof}

Observe that, when $V_j=0$, we have $f^{\max}_{V_jW_j}$ constant of value $\frac{1}{1-W_j}$ and then 
$$
\left\Vert \widehat{\Theta }_{j}\left( t\right) \right\Vert _{2}=\sqrt{\frac{1+W_j}{4}}
$$
is a constant function of $t$. Moreover, when $W_j=0$, we have $f^{\max}_{V_jW_j}$ constant of value $1+V_j$ and then
$$
\left\Vert \widehat{\Theta }_{j}\left( t\right) \right\Vert _{2}=\sqrt{\frac{1+V_j}{4}}
$$
is a constant function of $t$.

\section{Asymptotic condition numbers}

In this section, we describe in more detail the asymptotic condition numbers $K_{\infty}(t,y_0,\widehat{z}_0)$ and $K_{\infty}(t,y_0)$ of Theorem \ref{one}. This description is given in the next theorem.
\begin{theorem}
\label{Ar1} Suppose that $\Lambda_1$ is single real or single complex. If $w^{(1)}y_{0}\neq 0$ and $w^{(1)}\widehat{z}
_{0}\neq 0$, then 
\begin{equation}
K_{\infty}(t,y_{0},\widehat{z}_{0})=\frac{\left\vert \widehat{w}^{(1)}
\widehat{z}_{0}\right\vert}{\left\vert \widehat{w}^{(1)}\widehat{y}
_{0}\right\vert}\cdot\frac{g_1\left(t,\widehat{z}_0\right)}{g_1\left(t,
\widehat{y}_0\right)}  \label{first}
\end{equation}
and 
\begin{equation}
K_\infty(t,y_{0})=\frac{1}{\left\vert \widehat{w}^{(1)}\widehat{y}
_{0}\right\vert}\cdot\frac{g_1\left(t\right)}{g_1\left(t,\widehat{y}_0\right)
}  \label{second}
\end{equation}
where $g_1(t,\widehat{z}_0)$, $g_1(t,\widehat{y}_0)$ and $g_1(t)$ are defined in Subsection \ref{thenumbers}.
\end{theorem}

\begin{proof}
By Theorem \ref{one}, we have
\begin{equation*}
K_\infty\left( t,y_{0},\widehat{z}_{0}\right)=\frac{\left\Vert
Q_{1}(t)\widehat{z}_{0}\right\Vert }{\left\Vert Q_{1}(t)
\widehat{y}_{0}\right\Vert }\text{\ and\ \ }K_\infty\left( t,y_{0}\right) = \frac{\left\Vert Q_{1}(t)\right\Vert }{\left\Vert Q_{1}(t)\widehat{y}_{0}\right\Vert }.
\end{equation*}
Now, (\ref{first}) and (\ref{second}) follow by (\ref{normQJ0ureal}) and (\ref{normQJ0}).
\end{proof}

Separating the two cases $\Lambda_1$ single real and $\Lambda_1$ single complex, and using Lemmas \ref{Lemmareal} and \ref{Lemmacomplex}, we can restate the previous theorem as follows.

\begin{itemize}
\item If $\Lambda_1$ is single real,
then 
\begin{equation*}
K_\infty(t,y_{0},\widehat{z}_{0})=\frac{\left\vert \widehat{w}^{(1)}
\widehat{z}_{0}\right\vert }{\left\vert \widehat{w}^{(1)}\widehat{y}
_{0}\right\vert }
\end{equation*}
and
\begin{equation*}
K_\infty(t,y_{0})=\frac{1}{\left\vert \widehat{w}^{(1)}\widehat{y}
		_{0}\right\vert }.\label{K2}
\end{equation*}
\item If $\Lambda_1$ is single complex, then 
\begin{equation}
K_\infty(t,y_{0},\widehat{z}_{0})=\frac{\left\vert \widehat{w}^{(1)}
\widehat{z}_{0}\right\vert }{\left\vert \widehat{w}^{(1)}\widehat{y}
_{0}\right\vert }\cdot \frac{\left\Vert \widehat{\Theta}_1
\left( t,\widehat{z}_{0}\right) \right\Vert }{\left\Vert \widehat{\Theta}_1
\left( t,\widehat{y}_{0}\right) \right\Vert }  \label{KK1}
\end{equation}
and
\begin{equation}
	K_\infty(t,y_{0})= \frac{1}{\left\vert \widehat{w}^{(1)}\widehat{y}
		_{0}\right\vert }\cdot \frac{\left\Vert \widehat{\Theta}_1
		\left( t\right) \right\Vert }{\left\Vert \widehat{\Theta}_1
		\left( t,\widehat{y}_{0}\right) \right\Vert }. \label{KK2}
\end{equation}
\end{itemize}

\subsection{The oscillation scale factors and the oscillating terms}\label{sectionOSFOT}

Consider the case where $\Lambda_1$ is single complex. We can rewrite (\ref{KK1}) and (\ref{KK2}) as
\begin{eqnarray*}
&&K_\infty(t,y_{0},\widehat{z}_{0})=\mathrm{OSF}(y_0,\widehat{z}_0)\cdot\mathrm{OT}(t,y_0,\widehat{z}_0)\notag \\
&&K_\infty(t,y_{0})=\mathrm{OSF}(y_0)\cdot\mathrm{OT}(t,y_0), \label{KK22}
\end{eqnarray*}
where the constants 
\begin{eqnarray}
&&\mathrm{OSF}(y_0,\widehat{z}_0):=\frac{\left\vert \widehat{w}^{(1)}\widehat{z}_{0}\right\vert }{\left\vert 
\widehat{w}^{(1)}\widehat{y}_{0}\right\vert }\label{OSFz}\\
&&\mathrm{OSF}(y_0):=\frac{1}{
\left\vert \widehat{w}^{(1)}\widehat{y}_{0}\right\vert }\label{OSF}
\end{eqnarray}
are called the \emph{oscillation scale factors}
and the functions
\begin{eqnarray}
&&t\mapsto \mathrm{OT}(t,y_0,\widehat{z}_0):=\frac{\left\Vert \widehat{\Theta}_1 \left( t,\widehat{z}
_{0}\right) \right\Vert }{\left\Vert \widehat{\Theta}_1
\left( t,\widehat{y}_{0}\right) \right\Vert },\ t\in\mathbb{R}, \label{OTz}\\
&&t\mapsto \mathrm{OT}(t,y_0):=\frac{\left\Vert \widehat{\Theta}_1 \left( t\right) \right\Vert }{
\left\Vert \widehat{\Theta}_1 \left( t,\widehat{y}
_{0}\right) \right\Vert },\ t\in\mathbb{R}, \label{OT}
\end{eqnarray}
which are  periodic functions of period $\frac{\pi}{\omega_1}$, are called the \emph{oscillating terms}.

\begin{remark} \label{OSFOTRemark}
	
	The oscillation scale factor $\mathrm{OSF}(y_0, \widehat{z}_0)$ and the oscillating term $\mathrm{OT}\left(t, y_0, \widehat{z}_0\right)$ separate the contributions of moduli and angles in the polar forms of $\widehat{w}^{(1)}\widehat{y}_{0}$ and $\widehat{w}^{(1)}\widehat{z}_{0}$.
	
	Specifically, $\mathrm{OSF}(y_0, \widehat{z}_0)$ depends on $y_0$ and $\widehat{z}_0$ through the moduli of the complex numbers $\widehat{w}^{(1)}\widehat{y}_{0}$ and $\widehat{w}^{(1)}\widehat{z}_{0}$, and $\mathrm{OT}\left(t, y_0, \widehat{z}_0\right)$ depends on $y_0$ and $\widehat{z}_0$ through the angles $\gamma_1(\widehat{y}_0)$ and $\gamma_1(\widehat{z}_0)$ in these polar forms (refer to the definition of $\widehat{\Theta}(t, u)$ in Lemma \ref{Lemmacomplex}).
	
	Similarly, $\mathrm{OSF}(y_0)$ depends on $y_0$ through the modulus of the complex number $\widehat{w}^{(1)}\widehat{y}_{0}$ and $\mathrm{OT}\left(t, y_0\right)$ depends on $y_0$ through the angle $\gamma_1(\widehat{y}_0)$ in this polar form.
\end{remark}

\section{Oscillation scale factors and oscillating terms for the euclidean norm}\label{Section5}

\label{OSFOT}

In this section, we study the oscillation scale factors and the oscillating terms when the Euclidean norm is used as the vector norm. In particular, we derive expressions for the oscillation scale factors in terms of the components $c_1(u)$ and $d_1(u)$ of $u=\widehat{y}_0,\widehat{z}_0$ along the first and second right singular vectors, respectively, of the matrix $R_1$ of Proposition \ref{euclidean}. Moreover, for the oscillating terms, we analyze the extreme values
$$
\max\limits_{t\in\mathbb{R}}\mathrm{OT}(t,y_0,\widehat{z}_0),\ \min\limits_{t\in\mathbb{R}}\mathrm{OT}(t,y_0,\widehat{z}_0),\ \max\limits_{t\in\mathbb{R}}\mathrm{OT}(t,y_0)\text{\ \ and\ \ }\min\limits_{t\in\mathbb{R}}\mathrm{OT}(t,y_0)
$$
as $y_0$ and $\widehat{z}_0$ vary, always in terms of the components $c_1(u)$ and $d_1(u)$ of $u=\widehat{y}_0,\widehat{z}_0$.

\subsection{The oscillation scale factor $\mathrm{OSF}(y_0,\widehat{z}_0)$}
The oscillation scale factor\\ $\mathrm{OSF}(y_0,\widehat{z}_0)$ in (\ref{OSFz}) is the ratio of the moduli of the non-zero complex numbers $\widehat{w}^{(1)}\widehat{y}_{0}$ and $\widehat{w}^{(1)}\widehat{z}_{0}$. In case of the Euclidean norm as vector norm, such complex numbers can be arbitrary non-zero points in the set $\mathcal{W}_1$ of Proposition \ref{euclidean}.

Thus, as $y_0$ and $\widehat{z}_0$ vary, $\mathrm{OSF}(y_0,\widehat{z}_0)$ can be any number in
$(0,+\infty)$ when $n>2$ and any number in
$$
\left[\frac{\mu_1}{\sigma_1},\frac{\sigma_1}{\mu_1}\right]=\left[\sqrt{\frac{1-W_1}{1+W_1}},\sqrt{\frac{1+W_1}{1-W_1}}\right],
$$
where $\sigma_1$ and $\mu_1$ are the lengths of the semi-axes of the ellipse $\mathcal{E}_1$ of Proposition \ref{euclidean}, when $n=2$.

\begin{remark} 
In the case $n=2$, the eigenvalues of the complex conjugate pair in $\Lambda_1$ are all the eigenvalues of $A$. Therefore, we have $K(t,y_0,\widehat{z}_0)=K_\infty(t,y_0,\widehat{z}_0)$ and $K(t,y_0)=K_\infty(t,y_0)$ (we have $\epsilon(t,u)=0$, $u=\widehat{y}_0,\widehat{z}_0$, and $\epsilon(t)=0$ in Theorem \ref{two}).
\end{remark}

By (\ref{size}), we can express $\mathrm{OSF}(y_0,\widehat{z}_0)$ in terms of the components $c_1(u)$ and $d_1(u)$ of $u=\widehat{y}_0,\widehat{z}_0$  along the first two right singular vectors of the matrix $R_1$ of Proposition \ref{euclidean}. We have
\begin{equation}
	\mathrm{OSF}(y_0,\widehat{z}_0)=\sqrt{\frac{(1+W_1)c_1(\widehat{z}_0)^2+(1-W_1)d_1(\widehat{z}_0)^2}{(1+W_1)c_1(\widehat{y}_0)^2+(1-W_1)d_1(\widehat{y}_0)^2}}. \label{OSFz1}
\end{equation}

\subsection{The oscillating term $\mathrm{OT}\left(t,y_0,\widehat{z}_0\right)$}
We are interested in the extreme values 
\begin{equation}
	\max\limits_{t\in\mathbb{R}} \mathrm{OT}\left(t,y_0,\widehat{z}_0\right)\text{\ \ and\ \ }\min\limits_{t\in\mathbb{R}} \mathrm{OT}\left(t,y_0,\widehat{z}_0\right) \label{extreme}
\end{equation}
of the oscillating term, as $t$ varies.

Let us introduce the functions $G^{\max}_{V_1},G^{\min}_{V_1}:\mathbb{R}\rightarrow \mathbb{R}$ given by
\begin{equation}
G^{\max}_{V_1}(\gamma)=\sqrt{f_{V_1V_1}^{\max}(2\gamma-\pi)}\text{\ \ and\ \ }G^{\min}_{V_1}(\gamma)=\sqrt{f_{V_1V_1}^{\min}(2\gamma-\pi)},\ \gamma\in\mathbb{R},\label{GmaxminV1}
\end{equation}
where the functions $f_{V_1V_1}^{\max}$ and $f_{V_1V_1}^{\min}$ are defined in (\ref{fmax0}) and (\ref{fmin0}).

Next theorem specifies the extreme values (\ref{extreme}).

\begin{theorem}\label{ThmOTyz}
We have
	\begin{equation}
		\max\limits_{t\in\mathbb{R}} \mathrm{OT}\left(t,y_0,\widehat{z}_0\right)=G^{\max}_{V_1}\left(\gamma
		_{1}(\widehat{z}_0)-\gamma
		_{1}(\widehat{y}_0)\right) \label{alphamax0}
\end{equation}
and
\begin{equation}
		\min\limits_{t\in\mathbb{R}} \mathrm{OT}\left(t,y_0,\widehat{z}_0\right)=G^{\min}_{V_1}\left(\gamma
		_{1}(\widehat{z}_0)-\gamma
		_{1}(\widehat{y}_0)\right), \label{alphamin0}
	\end{equation}
	where $\gamma_1\left(\widehat{z}_0\right)$ and $\gamma_1\left(\widehat{y}_0\right)$ are defined in Subsection \ref{notations}.
\end{theorem} 

\begin{proof}
By recalling Proposition \ref{oscillation}, we can write the oscillating term (\ref{OTz}) as
\begin{equation}
	\mathrm{OT}\left(t,y_0,\widehat{z}_0\right)
	=\sqrt{\frac{1+V_1\cos(\Omega_1(t)+\Delta_1(\widehat{z}_0))}{1+V_1\cos(\Omega_1(t)+\Delta_1(\widehat{y}_0))}}=\sqrt{f_{V_1V_1}\left(\alpha,x\right)},  \label{OTyz2}
\end{equation}
where $f_{V_1V_1}$ is the function defined in (\ref{fVW0}), 
\begin{equation}
\alpha=\Omega_1(t)+\Delta_1(\widehat{y}_0)+\pi=2(\omega_{1}t+\gamma_1(\widehat{y}_0))+\delta_1+\pi \label{alpha}
\end{equation}
and
\begin{equation}
x=\Delta_1(\widehat{z}_0)-\Delta_1(\widehat{y}_0)-\pi=2\left(\gamma
_{1}(\widehat{z}_0)-\gamma
_{1}(\widehat{y}_0)\right)-\pi. \label{x=}
\end{equation}
with $\Omega_1(t)$, $\Delta_1(\widehat{y}_0)$ and $\Delta_1(\widehat{z}_0)$ defined in (\ref{xt}) and (\ref{xt2}). Observe that $\alpha$ and $x$ are chosen so that
$$
\frac{1+V_1\cos(\Omega_1(t)+\Delta_1(\widehat{z}_0))}{1+V_1\cos(\Omega_1(t)+\Delta_1(\widehat{y}_0))}=\frac{1+V_1\cos(x+\alpha)}{1-V_1\cos\alpha}=f_{V_1V_1}\left(\alpha,x\right)
$$
in (\ref{OTyz2}).

Since $\alpha$ in (\ref{alpha}) varies over all $\mathbb{R}$, as $t$ varies over all $\mathbb{R}$, we obtain
$$
\max\limits_{t\in\mathbb{R}} \mathrm{OT}\left(t,y_0,\widehat{z}_0\right)=\sqrt{f_{V_1V_1}^{\max}\left(x\right)}\text{\ \ and\ \ }\min\limits_{t\in\mathbb{R}} \mathrm{OT}\left(t,y_0,\widehat{z}_0\right)=\sqrt{f_{V_1V_1}^{\min}\left(x\right)},
$$
i.e., (\ref{alphamax0}) and (\ref{alphamin0}).
\end{proof}

The forms of the functions  $G^{\max}_{V_1}$ and $G^{\min}_{V_1}$ are stated in the next Theorem.

\begin{theorem}\label{ThOTz}
The function $G^{\max}_{V_1}$ is $\pi$-periodic, even and it oscillates monotonically between the extreme values
\begin{equation*}
\max\limits_{\gamma\in \mathbb{R}}G^{\max}_{V_1}\left(\gamma\right)=\sqrt{\frac{1+V_1}{1-V_1}}\text{\ \ and\ \ }\min\limits_{\gamma\in \mathbb{R}}G^{\max}_{V_1}\left(\gamma\right)=1.
\end{equation*}
The maximum value is achieved if $\gamma$ is an odd multiple of $\frac{\pi}{2}$ and the minimum value is achieved if $\gamma$ is a multiple of $\pi$.

The function $G^{\min}_{V_1}$ is $\pi$-periodic, even and it oscillates monotonically between the extreme values
\begin{equation*}
\max\limits_{\gamma\in \mathbb{R}}G^{\min}_{V_1}\left(\gamma\right)=1\text{\ \ and\ \ }\min\limits_{\gamma\in \mathbb{R}}G^{\min}_{V_1}\left(\gamma\right)=\sqrt{\frac{1-V_1}{1+V_1}}.
\end{equation*}
The maximum value is achieved if $\gamma$ is a multiple of $\pi$ and  the minimum value is achieved if $\gamma$ is an odd multiple of $\frac{\pi}{2}$.
\end{theorem}

\begin{proof}
The function $f_{V_1V_1}^{\max}$ oscillates monotonically between the extreme values
\begin{equation*}
	\max\limits_{x\in\mathbb{R}} f_{V_1V_1}^{\max}\left(x\right)
	=\frac{1+V_1}{1-V_1}
\end{equation*}
and
\begin{equation*}
	\min\limits_{x\in\mathbb{R}} f_{V_1V_1}^{\max}\left(x\right)
	=1
\end{equation*}
with the maximum value achieved if $x$ is an even multiple of $\pi$ and the minimum value achieved if $x$ is an odd multiple of $\pi$. Moreover, the function $f_{V_1V_1}^{\min}$ oscillates monotonically between the extreme values
\begin{equation*}
	\max\limits_{x\in\mathbb{R}} f_{V_1V_1}^{\min}\left(x\right)
	=1
\end{equation*}
and
\begin{equation*}
	\min\limits_{x\in\mathbb{R}} f_{V_1V_1}^{\min}\left(x\right)
	=\frac{1-V_1}{1+V_1}
\end{equation*}
with the maximum value achieved if $x$ is an odd multiple of $\pi$ and the minimum value achieved if $x$ is an even multiple of $\pi$. See Appendix \ref{AA}.

The theorem now follows from the fact that the argument $2\gamma-\pi$  of $f^{\max}_{V_1V_1}$ and $f^{\min}_{V_1V_1}$ in (\ref{GmaxminV1}), considered as a function of $\gamma$, is an orientation preserving homeomorphism of $\mathbb{R}$.

That $G^{\max}_{V_1}$ and $G^{\min}_{V_1}$ are $\pi$-periodic and even follows from that fact that  $f^{\max}_{V_1V_1}$ and $f^{\min}_{V_1V_1}$ are $2\pi$-periodic and even (see Appendix \ref{AA}).   
\end{proof}

\subsection{How the extreme values $	\max\limits_{t\in\mathbb{R}} \mathrm{OT}\left(t,y_0,\widehat{z}_0\right)$ and $\min\limits_{t\in\mathbb{R}} \mathrm{OT}\left(t,y_0,\widehat{z}_0\right)$ depend on $y_0$ and $\widehat{z}_0$}\label{How1}

As $y_0$ and $\widehat{z}_0$ vary in $\mathbb{R}^n$ with $\Vert \widehat{z}_0\Vert_2=1$, $w^{(1)}y_0\neq 0$ and $w^{(1)}\widehat{z}_0\neq 0$, the complex numbers $\widehat{w}^{(1)}y_0$ and $\widehat{w}^{(1)}\widehat{z}_0$ can be arbitrary non-zero points in the set $\mathcal{W}_1$ of Proposition \ref{euclidean}. Then, the angle $\gamma=\gamma_1\left(\widehat{z}_0\right)-\gamma_1\left(\widehat{y}_0\right)$ of the vector $\widehat{w}^{(1)}\widehat{z}_0$ with respect to the vector $\widehat{w}^{(1)}y_0$ varies in $[0,2\pi]$.

Therefore, by Theorem \ref{ThOTz}, the extreme values $\max\limits_{t\in\mathbb{R}}\mathrm{OT}(t,y_0,\widehat{z}_0)=G^{\max}_{V_1}(\gamma)$ and $\min\limits_{t\in\mathbb{R}}\mathrm{OT}(t,y_0,\widehat{z}_0)=G^{\min}_{V_1}(\gamma)$ (see Theorem \ref{ThmOTyz}) behave as follows, as the angle $\gamma$ varies.
\begin{itemize}
	\item  [a)] When $\gamma$ increases from $0$ to $\frac{\pi}{2}$, the extreme value $\max\limits_{t\in\mathbb{R}} \mathrm{OT}\left(t,y_0,\widehat{z}_0\right)$ increases from its minimum value
	$$
	\min\limits_{\substack{y_0,\widehat{z}_0\in\mathbb{R}^n\\ \left\Vert \widehat{z}_0\right\Vert_2=1\\ w^{(1)}y_0\neq 0,w^{(1)}\widehat{z}_0\neq 0}}\max\limits_{t\in\mathbb{R}} \mathrm{OT}\left(t,y_0,\widehat{z}_0\right)=G^{\max}_{V_1}(0)=1
	$$
	to its maximum value
	$$
	\max\limits_{\substack{y_0,\widehat{z}_0\in\mathbb{R}^n\\ \left\Vert \widehat{z}_0\right\Vert_2=1\\ w^{(1)}y_0\neq 0,w^{(1)}\widehat{z}_0\neq 0}}\max\limits_{t\in\mathbb{R}} \mathrm{OT}\left(t,y_0,\widehat{z}_0\right)=G^{\max}_{V_1}\left(\frac{\pi}{2}\right)=\sqrt{\frac{1+V_1}{1-V_1}}
	$$
	and the extreme value $\min\limits_{t\in\mathbb{R}} \mathrm{OT}\left(t,y_0,\widehat{z}_0\right)$ decrease from its maximum value
	$$
	\max\limits_{\substack{y_0,\widehat{z}_0\in\mathbb{R}^n\\ \left\Vert \widehat{z}_0\right\Vert_2=1\\ w^{(1)}y_0\neq 0,w^{(1)}\widehat{z}_0\neq 0}}\min\limits_{t\in\mathbb{R}} \mathrm{OT}\left(t,y_0,\widehat{z}_0\right)=G^{\min}_{V_1}(0)=1
	$$
	to its minimum value
	$$
	\min\limits_{\substack{y_0,\widehat{z}_0\in\mathbb{R}^n\\ \left\Vert \widehat{z}_0\right\Vert_2=1\\ w^{(1)}y_0\neq 0,w^{(1)}\widehat{z}_0\neq 0}}\min\limits_{t\in\mathbb{R}} \mathrm{OT}\left(t,y_0,\widehat{z}_0\right)=G^{\min}_{V_1}\left(\frac{\pi}{2}\right)=\sqrt{\frac{1-V_1}{1+V_1}}.
	$$
	\item [b)] When $\gamma$ increases from $\frac{\pi}{2}$ to $\pi$, the extreme value $\max\limits_{t\in\mathbb{R}} \mathrm{OT}\left(t,y_0,\widehat{z}_0\right)$ decreases from its maximum value $G^{\max}_{V_1}\left(\frac{\pi}{2}\right)$ to its minimum value $G^{\max}_{V_1}(0)$ and the extreme value $\min\limits_{t\in\mathbb{R}} \mathrm{OT}\left(t,y_0,\widehat{z}_0\right)$ increases from its minimum value $G^{\max}_{V_1}\left(\frac{\pi}{2}\right)$ to its maximum value $G^{\max}_{V_1}(0)$. This happens in a symmetric way with respect to the process in a): in fact, since $G^{\max}_{V_1}$ and $G^{\min}_{V_1}$ are $\pi$-periodic and even, we have, for $\gamma\in\left[0,\frac{\pi}{2}\right]$,
	\begin{equation*}
	\quad\quad\quad  G^{\max}_{V_1}\left(\frac{\pi}{2}-\gamma\right)=G^{\max}_{V_1}\left(\frac{\pi}{2}+\gamma\right)\text{\ \ and\ \ }G^{\min}_{V_1}\left(\frac{\pi}{2}-\gamma\right)=G^{\min}_{V_1}\left(\frac{\pi}{2}+\gamma\right).
	\end{equation*}
	\item [c)] When $\gamma$ increases from $\pi$ to $\frac{3\pi}{2}$, the process in a) repeats.
	\item [d)] When $\gamma$ increases from $\frac{3\pi}{2}$ to $2\pi$, the process in b) repeats.
	\end{itemize}

The following subsections contain important considerations related to the previous description of the extreme values of the oscillating term.

\subsubsection {$\gamma_{1}(\widehat{z}_0)-\gamma_{1}(\widehat{y}_0)$ odd multiple of $\frac{\pi}{2}$}\label{oddmultiplepi/2}

	The condition $\gamma_{1}(\widehat{z}_0)-\gamma_{1}(\widehat{y}_0)$ is an odd multiple of $\frac{\pi}{2}$ is equivalent to the vectors $\widehat{w}^{(1)}\widehat{y}_{0}$ and $\widehat{w}^{(1)}\widehat{z}_{0}$ being orthogonal in the complex plane, i.e.,
\begin{equation}
	(1+W_1)c_1(\widehat{y}_0)c_1(\widehat{z}_0)+(1-W_1)d_1(\widehat{y}_0)d_1(\widehat{z}_0)=0
	\label{ocnd}
\end{equation}
by reminding (\ref{whatu})-(\ref{whatu2}).

	When the orthogonality condition (\ref{ocnd}) holds, we have
	\begin{equation}
		\max\limits_{t\in\mathbb{R}} \mathrm{OT}\left(t,y_0,\widehat{z}_0\right)=\sqrt{\frac{1+V_1}{1-V_1}}\text{\ \ and\ \ }
		\min\limits_{t\in\mathbb{R}} \mathrm{OT}\left(t,y_0,\widehat{z}_0\right)=\sqrt{\frac{1-V_1}{1+V_1}} \label{maxmaxz2}
	\end{equation}
	with the $\max$  at its maximum value and the $\min$ at its minimum value, as $y_0$ and $\widehat{z}_0$ vary.

\subsubsection{$\gamma_{1}(\widehat{z}_0)-\gamma_{1}(\widehat{y}_0)$ multiple of $\pi$} \label{multipleofpi}

The condition $\gamma_{1}(\widehat{z}_0)-\gamma_{1}(\widehat{y}_0)$ is a multiple of $\pi$ is equivalent to the vectors $\widehat{w}^{(1)}\widehat{y}_{0}$ and $\widehat{w}^{(1)}\widehat{z}_{0}$ being parallel in the complex plane, i.e.,
\begin{equation}
	c_1(\widehat{y}_0)d_1(\widehat{z}_0)-d_1(\widehat{y}_0)c_1(\widehat{z}_0)=0.
	\label{pcnd}
\end{equation}	

When the parallelism condition (\ref{pcnd}) holds, we have
\begin{equation*}
\max\limits_{t\in\mathbb{R}} \mathrm{OT}\left(t,y_0,\widehat{z}_0\right)=1=\min\limits_{t\in\mathbb{R}} \mathrm{OT}\left(t,y_0,\widehat{z}_0\right) \label{maxminz}
\end{equation*}
with the $\max$  at its minimum value and $\min$ at its maximum value, as $y_0$ and $\widehat{z}_0$ vary. In this case, 
$$
\mathrm{OT}\left(t,y_0,\widehat{z}_0\right)=1
$$
holds, as it can be directly verified by (\ref{OTyz2}) and the fact that $\Delta_1(\widehat{z}_0)-\Delta_1(\widehat{y}_0)=2(\gamma_1(\widehat{z}_0)-\gamma_1(\widehat{y}_0))$ is a multiple of $2\pi$.

\subsubsection{$V_1$ not close to $1$}\label{V1nc1}

If $V_1$ is not close to $1$, then, for any $y_0$ and $\widehat{z}_0$, $\mathrm{OT}\left(t,y_0,\widehat{z}_0\right)$ does not assume large or small  values, as $t$ varies. In fact, we have
$$
\sqrt{\frac{1-V_1}{1+V_1}}\leq\min\limits_{t\in\mathbb{R}} \mathrm{OT}\left(t,y_0,\widehat{z}_0\right)\text{\ \ and\ \ } \max\limits_{t\in\mathbb{R}} \mathrm{OT}\left(t,y_0,\widehat{z}_0\right)\leq\sqrt{\frac{1+V_1}{1-V_1}}.
$$

As a consequence, if $V_1$ is not close to $1$, then the asymptotic condition number
$$
K_\infty(t,y_0,\widehat{z}_0)=\mathrm{OSF}(y_0,\widehat{z}_0)\cdot \mathrm{OT}(t,y_0,\widehat{z}_0),
$$
has the same order of magnitude as the oscillation scale factor $\mathrm{OSF}(y_0,\widehat{z}_0)$.

\subsubsection{Critical case}\label{critcase} There is a critical case: $V_1$ close to $1$.

If $V_1$ is close to $1$, then, for $y_0$ and $\widehat{z}_0$ satisfying the orthogonality condition (\ref{ocnd}), $\mathrm{OT}\left(t,y_0,\widehat{z}_0\right)$ assumes both large and small values, as $t$ varies (see (\ref{maxmaxz2})). In other words, in such a critical case, $\mathrm{OT}\left(t,y_0,\widehat{z}_0\right)$ exhibits large oscillations, from small to large values.

Suppose $n>2$. Since the oscillating scale factor $\mathrm{OSF}(y_0,\widehat{z}_0)$ (ratio of the moduli of $\widehat{w}^{(1)}\widehat{y}_{0}$ and $\widehat{w}^{(1)}\widehat{z}_{0}$) can be any number in $(0,+\infty)$ even if $\widehat{w}^{(1)}\widehat{y}_{0}$ and $\widehat{w}^{(1)}\widehat{z}_{0}$ are orthogonal point of the set $\mathcal{W}_1$ of Proposition \ref{euclidean}, nothing can be said about the asymptotic condition number $K_\infty(t,y_0,\widehat{z}_0)$ when $V$ is close to $1$: it can attain values much larger than the large values of $\mathrm{OT}\left(t,y_0,\widehat{z}_0\right)$  when $\mathrm{OSF}\left(y_0,\widehat{z}_0\right)$ is large, but it can avoid to assume large values when a small $\mathrm{OSF}\left(y_0,\widehat{z}_0\right)$ compensates for large values of $\mathrm{OT}\left(t,y_0,\widehat{z}_0\right)$.

\subsubsection{Specific cases}\label{690}
Some specific cases are worth considering.
\begin{itemize}
	\item [a)] Consider $y_0$ orthogonal to the first right singular vector of the matrix $R_1$ of Proposition \ref{euclidean} and $\widehat{z}_0$ orthogonal to the second right singular vector (i.e., $c_1(\widehat{y}_0)=0$ and $d_1(\widehat{z}_0)=0$).
	
	We have
	\begin{equation*}
		\mathrm{OSF}(y_0,\widehat{z}_0)=\sqrt{\frac{1+W_1}{1-W_1}}\cdot \frac{\vert c_1(\widehat{z}_0)\vert}{\vert d_1(\widehat{y}_0)\vert}.
	\end{equation*}
	by (\ref{OSFz1}). Moreover, since $\widehat{y}_0$ and $\widehat{z}_0$ satisfy the orthogonality condition (\ref{ocnd}), we have (\ref{maxmaxz2}), i.e., 
	\begin{equation*}
		\max\limits_{t\in\mathbb{R}} \mathrm{OT}\left(t,y_0,\widehat{z}_0\right)=\sqrt{\frac{1+V_1}{1-V_1}} \text{\ and\ }
		\min\limits_{t\in\mathbb{R}} \mathrm{OT}\left(t,y_0,\widehat{z}_0\right)=\sqrt{\frac{1-V_1}{1+V_1}}.
	\end{equation*}
	Regarding the asymptotic condition number $K_\infty(t,y_0,\widehat{z}_0)$, we have
	\begin{equation}
		\max\limits_{t\in  \mathbb{R}}K_\infty(t,y_0,\widehat{z}_0)=\sqrt{\frac{(1+V_1)(1+W_1)}{(1-V_1)(1-W_1)}}\cdot \frac{\vert c_1(\widehat{z}_0)\vert}{\vert d_1(\widehat{y}_0)\vert} \label{Kyz1}
	\end{equation}
	and
	\begin{equation}
		\min\limits_{t\in  \mathbb{R}}K_\infty(t,y_0,\widehat{z}_0)=\sqrt{\frac{(1-V_1)(1+W_1)}{(1+V_1)(1-W_1)}}\cdot \frac{\vert c_1(\widehat{z}_0)\vert}{\vert d_1(\widehat{y}_0)\vert}. \label{Kyz11}
	\end{equation}
	We can conclude the following.
	\begin{itemize}
		\item It is expected that $\mathrm{OSF}(y_0,\widehat{z}_0)$ is large if and only if $W_1$ is close to $1$ (since the ratio $\frac{\vert c_1(\widehat{z}_0)\vert}{\vert d_1(\widehat{y}_0)\vert}$ is expected to be neither large nor small).
		\item $\mathrm{OT}(t,y_0,\widehat{z}_0)$ assumes large values, as $t$ varies, if and only if $V_1$ is close to $1$.
		\item $\mathrm{OT}(t,y_0,\widehat{z}_0)$ assumes small values, as $t$ varies, if and only if $V_1$ is close to $1$.
		\item It is expected that $K_\infty(t,y_0,\widehat{z}_0) $ assumes  large values, as $t$ varies, if and only if $V_1$ is close to $1$ or $W_1$ is close to $1$.
		\item It is expected that $K_\infty(t,y_0,\widehat{z}_0) $ assumes  small values, as $t$ varies, if and only if the ratio $\frac{1-V_1}{1-W_1}$ is small: in this case, even a large $\mathrm{OSF}(y_0,\widehat{z}_0)$ does not compensate for small values $\mathrm{OT}(t,y_0,\widehat{z}_0)$. 
	\end{itemize}
	
	\item [b)] Consider $y_0$ orthogonal to the second right singular vector of $R_1$ and $\widehat{z}_0$ orthogonal to the first right singular vector (i.e., $d_1(\widehat{y}_0)=0$ and $c_1(\widehat{z}_0)=0$).
	
	We have  
	\begin{equation*}
		\mathrm{OSF}(y_0,\widehat{z}_0)=\sqrt{\frac{1-W_1}{1+W_1}}\cdot \frac{\vert d_1(\widehat{z}_0)\vert}{\vert c_1(\widehat{y}_0)\vert}
	\end{equation*}
	and 
	\begin{equation*}
		\max\limits_{t\in\mathbb{R}} \mathrm{OT}\left(t,y_0,\widehat{z}_0\right)=\sqrt{\frac{1+V_1}{1-V_1}} \text{\ and\ }
		\min\limits_{t\in\mathbb{R}} \mathrm{OT}\left(t,y_0,\widehat{z}_0\right)=\sqrt{\frac{1-V_1}{1+V_1}}.
	\end{equation*}
	again. 
		Regarding the asymptotic condition number $K_\infty(t,y_0,\widehat{z}_0)$, we have
		\begin{equation}
		\max\limits_{t\in  \mathbb{R}}K_\infty(t,y_0,\widehat{z}_0)=\sqrt{\frac{(1+V_1)(1-W_1)}{(1-V_1)(1+W_1)}}\cdot \frac{\vert d_1(\widehat{z}_0)\vert}{\vert c_1(\widehat{y}_0)\vert} \label{Kyz2}
		\end{equation}
		and
		\begin{equation}
			\min\limits_{t\in  \mathbb{R}}K_\infty(t,y_0,\widehat{z}_0)=\sqrt{\frac{(1-V_1)(1-W_1)}{(1+V_1)(1+W_1)}}\cdot \frac{\vert d_1(\widehat{z}_0)\vert}{\vert c_1(\widehat{y}_0)\vert}. \label{Kyz22}
		\end{equation}
		We can conclude the following.
		\begin{itemize}
			\item It is expected that $\mathrm{OSF}(y_0,\widehat{z}_0)$ is small if and only if $W_1$ is close to $1$;
			\item $\mathrm{OT}(t,y_0,\widehat{z}_0)$ assumes large values, as $t$ varies, if only if $V_1$ is close to $1$;
			\item $\mathrm{OT}(t,y_0,\widehat{z}_0)$ assumes small values, as $t$ varies, if only if $V_1$ is close to $1$;
			\item It is expected that $K_\infty(t,y_0,\widehat{z}_0) $ assumes  large values, as $t$ varies, if and only if  the ratio $\frac{1-V_1}{1-W_1}$ is small: in this case, even a small $\mathrm{OSF}(y_0,\widehat{z}_0)$ does not compensate for large values $\mathrm{OT}(t,y_0,\widehat{z}_0)$.
			\item It is expected that $K_\infty(t,y_0,\widehat{z}_0) $ assumes  small values, as $t$ varies, if and only if $V_1$ is close to $1$ or $W_1$ is close to $1$.
		\end{itemize}

	\item [c)] Consider $y_0$ and $\widehat{z}_0$ orthogonal to the first right singular vector of $R_1$ (i.e. $c_1(\widehat{y}_0)=0$ and $c_1(\widehat{z}_0)=0$).
	
	We have $\mathrm{OT}\left(t,y_0,\widehat{z}_0\right)=1$ since $\widehat{y}_0$ and $\widehat{z}_0$ satisfy the parallelism condition (\ref{pcnd}). 	Thus, we have
	\begin{equation*}
		K_\infty(t,y_0,\widehat{z}_0)=\mathrm{OSF}(y_0,\widehat{z}_0)=\frac{\vert d_1(\widehat{z}_0)\vert}{\vert d_1(\widehat{y}_0)\vert}
	\end{equation*}
	by (\ref{OSFz1}).

	\item [d)] Consider $y_0$ and $\widehat{z}_0$ orthogonal to the second right singular vector of $R_1$ (i.e. $d_1(\widehat{y}_0)=0$ and $d_1(\widehat{z}_0)=0$).
	
	We have $\mathrm{OT}\left(t,y_0,\widehat{z}_0\right)=1$ again and
	\begin{equation*}
		K_\infty(t,y_0,\widehat{z}_0)=\mathrm{OSF}(y_0,\widehat{z}_0)=\frac{\vert c_1(\widehat{z}_0)\vert}{\vert c_1(\widehat{y}_0)\vert}.
	\end{equation*}
\end{itemize}

\subsection{The oscillation scale factor $\mathrm{OSF}(y_0)$}\label{SectionOSFy0}

The oscillation scale factor $\mathrm{OSF}(y_0)$ in (\ref{OSF}) is the reciprocal of the modulus  of the complex number $\widehat{w}^{(1)}\widehat{y}_{0}$. In case of the Euclidean norm as vector norm, such complex number can be an arbitrary non-zero point in the set $\mathcal{W}_1$ of Proposition \ref{euclidean}.

Thus, by varying $y_0$, $\mathrm{OSF}(y_0)$ can be any number in
$$
\left[\frac{1}{\sigma_1},+\infty\right)=\left[\sqrt{\frac{2}{1+W_1}},+\infty\right)
$$
when $n>2$ and any number in
$$
\left[\frac{1}{\sigma_1},\frac{1}{\mu_1}\right]=\left[\sqrt{\frac{2}{1+W_1}},\sqrt{\frac{2}{1-W_1}}\right].
$$
when $n=2$.

Observe that, unlike $\mathrm{OSF}(t,y_0,\widehat{z}_0)$,  $\mathrm{OSF}(t,y_0)$ is never small.

By (\ref{size}), we can express $\mathrm{OSF}(y_0)$ in terms of the components $c_1(\widehat{y}_0)$ and $d_1(\widehat{y}_0)$ of $\widehat{y}_0$ along the first two right singular vectors of the matrix $R_1$ of Proposition \ref{euclidean}. We have
	\begin{equation}
	\mathrm{OSF}(y_0)=\sqrt{\frac{2}{(1+W_1)c_1(\widehat{y}_0)^2+(1-W_1)d_1(\widehat{y}_0)^2}}. \label{OSF2}
\end{equation}

\subsection{The oscillating term $\mathrm{OT}(t,y_0)$} \label{OTy0}

We are interested in the extreme values 
\begin{equation}
	\max\limits_{t\in\mathbb{R}} \mathrm{OT}\left(t,y_0\right)\text{\ \ and\ \ }\min\limits_{t\in\mathbb{R}} \mathrm{OT}\left(t,y_0\right) \label{extreme1}
\end{equation}
of the oscillating term, as $t$ varies.

\begin{remark}\label{V1W1zero}
	In the case $W_1= 0$, we have (see (\ref{OTy})) 
	\begin{equation}
		\mathrm{OT}\left(t,y_0\right)=\sqrt{\frac{1}{2}\cdot\frac{1+V_1}{1+V_1\cos(\Omega_1(t)+\Delta_1(\widehat{y}_0))}}, \label{W10OT}
	\end{equation}
	since
	\begin{eqnarray*}
		f^{\max}_{V_1W_1}(x)=1+V_1,\ x\in\mathbb{R},
	\end{eqnarray*}
	(see  Appendix \ref{AA}). Note that, $\mathrm{OT}\left(t,y_0\right)$, as function of $t$, oscillates monotonically between the extreme values
	\begin{eqnarray*}
		\max\limits_{t\in\mathbb{R}}\mathrm{OT}\left(t,y_0\right)=\sqrt{\frac{1+V_1}{2(1-V_1)}}\text{\ \ and\ \ }\min\limits_{t\in\mathbb{R}}\mathrm{OT}\left(t,y_0\right)=\sqrt{\frac{1}{2}},
	\end{eqnarray*}
	Moreover, note that, as $y_0$ varies, the functions $\mathrm{OT}(t,y_0)$ of $t$ are all time shifts of the same function: for the the argument of $\cos$ in (\ref{W10OT}), we have, for $u,v\in\mathbb{R}^n$, 
	$$
	\Omega_1(t)+\Delta_1(v)=\Omega_1(t-\overline{t})+\Delta_1(u),\text{\ \ where $\overline{t}=\frac{\gamma_1(u)-\gamma_1(v)}{\omega_1}$}. 
	$$
\end{remark}

Let us introduce the functions $G^{\max}_{V_1W_1},G^{\min}_{V_1W_1}:\mathbb{R}\rightarrow \mathbb{R}$ given by
\begin{equation}
	G^{\max}_{V_1W_1}(\gamma)=\sqrt{\frac{1-W_1^2}{2}H_{V_1W_1}^{\max}(2\gamma)}\text{\ \ and\ \ }G^{\min}_{V_1W_1}(\gamma)=\sqrt{\frac{1-W_1^2}{2}H_{V_1W_1}^{\min}(2\gamma)},\ \gamma\in\mathbb{R},\label{GmaxminV1W1}
\end{equation}
where the functions $H_{V_1W_1}^{\max},H_{V_1W_1}^{\min}:\mathbb{R}\rightarrow \mathbb{R}$ are given by
\begin{equation*}
	H_{V_1W_1}^{\max}(\beta):=\max\limits_{x\in\mathbb{R}}H_{V_1W_1}(x,\beta)\text{\ \ and\ \ }H_{V_1W_1}^{\min}(\beta):=\min\limits_{x\in\mathbb{R}}H_{V_1W_1}(x,\beta),\ \beta\in\mathbb{R},
\end{equation*}
with the function  $H_{V_1W_1}:\mathbb{R}^2\rightarrow \mathbb{R}$ given by
\begin{equation*}
	H_{V_1W_1}(x,\beta)=\frac{f^{\max}_{V_1W_1}(x)}{1+V_1\cos(x+\beta)},\ (x,\beta)\in\mathbb{R}^2,
\end{equation*}
(see Appendix \ref{AB}).

Next theorem specifies the extreme values (\ref{extreme1}).

\begin{theorem}\label{formOTyz0}
	We have
	\begin{equation}
		\max\limits_{t\in\mathbb{R}} \mathrm{OT}\left(t,y_0\right)=G^{\max}_{V_1W_1}\left(\gamma
		_{1}(\widehat{y}_0)-\theta_1\right) \label{alphamax1}
	\end{equation}
	and
	\begin{equation}
		\min\limits_{t\in\mathbb{R}} \mathrm{OT}\left(t,y_0\right)=G^{\min}_{V_1W_1}\left(\gamma
		_{1}(\widehat{y}_0)-\theta_{1}\right), \label{alphamin1}
	\end{equation}
	where $\gamma_1\left(\widehat{y}_0\right)$ and $\theta_1$ are defined in Subsections \ref{notations} and \ref{notation2}, respectively.
\end{theorem}

\begin{proof}

By recalling Propositions \ref{oscillation} and \ref{oscillation2}, we can write the oscillating term (\ref{OT}) as
\begin{eqnarray}
\mathrm{OT}\left(t,y_0\right)&=&\sqrt{\frac{1-W_1^2}{2}\cdot\frac{f_{V_1W_1}^{\max}(\Omega_1(t))}{1+V_1\cos(\Omega_1(t)+\Delta_1(\widehat{y}_0))}}\notag \\
&=&\sqrt{\frac{1-W_1^2}{2}H_{V_1W_1}(\Omega_1(t),\Delta_1(\widehat{y}_0))}, \notag \\
\label{OTy}
\end{eqnarray}
where $\Omega_1(t)$ and $\Delta_1(\widehat{y}_0)$ are defined in (\ref{xt}) and (\ref{xt2}).

Since
$$
\Omega_1(t)=2(\omega_1t+\theta_1)+\delta_1
$$
varies over all $\mathbb{R}$, as $t$ varies, we obtain 
\begin{equation*}
	\max\limits_{t\in\mathbb{R}} \mathrm{OT}\left(t,y_0\right)
	=\sqrt{\frac{1-W_1^2}{2}H_{V_1W_1}^{\max}(\Delta_1(\widehat{y}_0))} \label{maxy}
\end{equation*}
and
\begin{equation*}
	\min\limits_{t\in\mathbb{R}}\mathrm{OT}\left(t,y_0\right)
	=\sqrt{\frac{1-W_1^2}{2}H_{V_1W_1}^{\min}(\Delta_1(\widehat{y}_0))}, \label{miny}
\end{equation*} 
i.e., (\ref{alphamax1}) and (\ref{alphamin1}) (recall $\Delta_1(\widehat{y}_0)=2\left(\gamma_1(\widehat{y}_0)-\theta_1\right)$).
\end{proof}

The forms of the functions  $G^{\max}_{V_1W_1}$ and $G^{\min}_{V_1W_1}$ are stated in the next Theorem.

\begin{theorem} \label{formOTyz}
The function $G_{V_1W_1}^{\max}$ is $\pi$-periodic, even and it oscillates monotonically between the values
\begin{equation*}
	\max\limits_{\gamma\in \mathbb{R}}G^{\max}_{V_1W_1}\left(\gamma\right)=\sqrt{\frac{(1+V_1)(1+W_1)}{2(1-V_1)}}
\end{equation*}
and
\begin{equation*}
\min\limits_{\overline{\gamma}\in \mathbb{R}}G^{\max}_{V_1W_1}\left(\gamma\right)=\left\{
	\begin{array}{l}
		\sqrt{\frac{(1-V_1^2)(1+W_1)}{2(1-Q_1V_1)}}\text{\ if\ }Q_1\leq 1\\
		\\
		\sqrt{\frac{(1+V_1)(1-W_1)}{2(1-V_1)}}\text{\ if\ }Q_1\geq 1,
	\end{array}
	\right.,
\end{equation*}
where
\begin{equation}
	Q_1:=\frac{V_1(1+W_1)}{2W_1}. \label{Quno}
\end{equation}
In the case $V_1\neq 0$ and $W_1=0$, we set $Q_1=+\infty$. In the case $V_1=0$ and $W_1=0$, $Q_1$ can be chosen arbitrarily.

The maximum value is achieved if $\gamma$ is an odd multiple of $\frac{\pi}{2}$ and the minimum value is achieved if $\gamma$ is a multiple of $\pi$.

The function $G^{\min}_{V_1W_1}$ is $\pi$-periodic, even and it oscillates monotonically between the values
\begin{equation*}
	\max\limits_{\gamma\in \mathbb{R}}G^{\min}_{V_1W_1}\left(\gamma\right)=\sqrt{\frac{1+W_1}{2}} \text{\ \ and\ \ }\min\limits_{\gamma\in \mathbb{R}}G^{\min}_{V_1W_1}\left(\gamma\right)=\left\{
	\begin{array}{l}
		\sqrt{\frac{(1-V_1)(1+W_1)}{2(1+V_1)}}\text{\ if\ }V_1\leq W_1\\
		\\
		\sqrt{\frac{1-W_1}{2}}\text{\ if\ }V_1\geq W_1.
	\end{array}
	\right.
\end{equation*}
The maximum value is achieved if $\gamma$ is a multiple of $\pi$ and  the minimum value is achieved if  $\gamma$ is an odd multiple of $\frac{\pi}{2}$.
\end{theorem}
\begin{proof}
	The function $H_{V_1W_1}^{\max}$ oscillates monotonically between the extreme values
	\begin{equation*}
		\max\limits_{\beta\in\mathbb{R}} H_{V_1W_1}^{\max}\left(\beta\right)
		=\frac{1+V_1}{(1-V_1)(1-W_1)}
	\end{equation*}
	and
	\begin{equation*}
		\min\limits_{\beta\in\mathbb{R}} H_{V_1W_1}^{\max}\left(\beta\right)
		=\left\{
		\begin{array}{l}
			\frac{1-V_1^2}{(1-Q_1V_1)(1-W_1)}\text{\ if\ }Q_1\leq 1\\
			\\
			\frac{1+V_1}{(1-V_1)(1+W_1)}\text{\ if\ }Q_1\geq 1,
		\end{array}
		\right.
	\end{equation*}
	with the maximum value achieved if $\beta$ is an odd multiple of $\pi$ and the minimum value achieved if $\beta$ is an even multiple of $\pi$. Moreover, the function $H_{V_1W_1}^{\min}$ oscillates monotonically between the extreme values
	\begin{equation*}
		\max\limits_{\beta\in\mathbb{R}} H_{V_1W_1}^{\min}\left(\beta\right)
		=\frac{1}{1-W_1}
	\end{equation*}
	and
	\begin{equation*}
		\min\limits_{\beta\in\mathbb{R}} H_{V_1W_1}^{\min}\left(\beta\right)
		=\left\{
		\begin{array}{l}
			\frac{1-V_1}{(1+V_1)(1-W_1)}\text{\ if\ }V_1\leq W_1\\
			\\
			\frac{1}{1+W_1}\text{\ if\ }V_1\geq W_1,
		\end{array}
		\right.
	\end{equation*}
	with the maximum value achieved if $\beta$ is an even multiple of $\pi$ and the minimum value achieved if $\beta$ is an odd multiple of $\pi$. See Appendix \ref{AB}.
	
	The theorem now follows from the fact that the argument $2\gamma$  of $H^{\max}_{V_1W_1}$ and $H^{\min}_{V_1W_1}$ in (\ref{GmaxminV1W1}), considered as a function of $\gamma$, is an orientation preserving homeomorphism of $\mathbb{R}$.
	
	That $G^{\max}_{V_11W_1}$ and $G^{\min}_{V_1W_1}$ are $\pi$-periodic and even follows from that fact that  $H^{\max}_{V_1W_1}$ and $H^{\min}_{V_1W_1}$ are $2\pi$-periodic and even (see Appendix \ref{AB}).   
\end{proof}

\subsection{How the extreme values $	\max\limits_{t\in\mathbb{R}} \mathrm{OT}\left(t,y_0\right)$ and $\min\limits_{t\in\mathbb{R}} \mathrm{OT}\left(t,y_0\right)$ depend on $y_0$}\label{How2}

As $y_0$ vary in $\mathbb{R}^n$ with $w^{(1)}y_0\neq 0$, the complex numbers $\widehat{w}^{(1)}y_0$ can be an arbitrary non-zero point in the set $\mathcal{W}_1$ of Proposition \ref{euclidean}. Then, the angle $\gamma=\gamma_1\left(\widehat{y}_0\right)-\theta_1$ of the vector $\widehat{w}^{(1)}\widehat{y}_0$ with respect to the major semi-axis of the ellipse $\mathcal{E}_1$ of of Proposition \ref{euclidean} varies in $[0,2\pi]$.

Therefore, by Theorem \ref{formOTyz}, the extreme values $\max\limits_{t\in\mathbb{R}}\mathrm{OT}(t,y_0)=G^{\max}_{V_1W_1}(\gamma)$ and $\min\limits_{t\in\mathbb{R}}\mathrm{OT}(t,y_0)=G^{\min}_{V_1W_1}(\gamma)$ behave as follows, as the angle $\gamma$ varies. We refer to Table \ref{taba}.

\begin{table}[ht]
	\centering
	\begin{tabular}{|l|}
		\hline
		$a^{\max\max}_{V_1W_1}:=G_{V_1,W_1}^{\max}\left(\frac{\pi}{2}\right)=\sqrt{\frac{(1+V_1)(1+W_1)}{2(1-V_1)}}$  \\
		\hline
		 $a^{\min\max}_{V_1W_1}:=G_{V_1,W_1}^{\max}\left( 0\right)=\left\{
		 \begin{array}{l}
		 	\sqrt{\frac{(1-V_1^2)(1+W_1)}{2(1-Q_1V_1)}}\text{\ if\ }Q_1\leq 1\\
		 	\\
		 	\sqrt{\frac{(1+V_1)(1-W_1)}{2(1-V_1)}}\text{\ if\ }Q_1\geq 1,
		 \end{array}
		 \right.$  \\
		\hline
		 $a^{\max\min}_{V_1W_1}:=G_{V_1,W_1}^{\min}\left(0\right)=\sqrt{\frac{1+W_1}{2}}$   \\
		\hline
		 $a^{\min\min}_{V_1W_1}:=G_{V_1,W_1}^{\min}\left(\frac{\pi}{2}\right)=\left\{
		 \begin{array}{l}
		 	\sqrt{\frac{(1-V_1)(1+W_1)}{2(1+V_1)}}\text{\ if\ }V_1\leq W_1\\
		 	\\
		 	\sqrt{\frac{1-W_1}{2}}\text{\ if\ }V_1\geq W_1.
		 \end{array}
		 \right.$   \\
		\hline
	\end{tabular}
	\caption{$a^{\max\max}_{V_1W_1}$, $a^{\min\max}_{V_1W_1}$, $a^{\max\min}_{V_1W_1}$ and $a^{\min\min}_{V_1W_1}$.}
	\label{taba}
\end{table}
\begin{itemize}
	\item  [a)] When $\gamma$ increases from $0$ to $\frac{\pi}{2}$, the extreme value $\max\limits_{t\in\mathbb{R}} \mathrm{OT}\left(t,y_0\right)$ increases from its minimum value
	$$
	\min\limits_{\substack{y_0\in\mathbb{R}^n\\ w^{(1)}y_0\neq 0}}\max\limits_{t\in\mathbb{R}} \mathrm{OT}\left(t,y_0\right)=a^{\min\max}_{V_1W_1}
	$$
	to its maximum value
	$$
	\max\limits_{\substack{y_0\in\mathbb{R}^n\\ w^{(1)}y_0\neq 0}}\max\limits_{t\in\mathbb{R}} \mathrm{OT}\left(t,y_0\right)=a^{\max\max}_{V_1W_1}
	$$
	and the extreme value $\min\limits_{t\in\mathbb{R}} \mathrm{OT}\left(t,y_0\right)$ decrease from its maximum value
	$$
	\max\limits_{\substack{y_0\in\mathbb{R}^n\\ w^{(1)}y_0\neq 0\neq 0}}\min\limits_{t\in\mathbb{R}} \mathrm{OT}\left(t,y_0,\widehat{z}_0\right)=a^{\max\min}_{V_1W_1}
	$$
	to its minimum value
	$$
	\min\limits_{\substack{y_0\in\mathbb{R}^n\\\ w^{(1)}y_0\neq 0\neq 0}}\min\limits_{t\in\mathbb{R}} \mathrm{OT}\left(t,y_0,\widehat{z}_0\right)=a^{\min\min}_{V_1W_1}.
	$$
	\item [b)] When $\gamma$ increases from $\frac{\pi}{2}$ to $\pi$, the extreme value $\max\limits_{t\in\mathbb{R}} \mathrm{OT}\left(t,y_0\right)$ decreases from its maximum value $a^{\max\max}_{V_1W_1}$ to its minimum value $a^{\min\max}_{V_1W_1}$ and the extreme value $\min\limits_{t\in\mathbb{R}} \mathrm{OT}\left(t,y_0\right)$ increases from its minimum value $a^{\min\min}_{V_1W_1}$ to its maximum value $a^{\max\min}_{V_1W_1}$. This happens in a symmetric way with respect to the process in a): in fact, since $G^{\max}_{V_1W_1}$ and $G^{\min}_{V_1W_1}$ are $\pi$-periodic and even, we have, for $\gamma\in\left[0,\frac{\pi}{2}\right]$,
	\begin{equation*}
		\quad\quad\quad\quad  G^{\max}_{V_1W_1}\left(\frac{\pi}{2}-\gamma\right)=G^{\max}_{V_1W_1}\left(\frac{\pi}{2}+\gamma\right)\text{\ \ and\ \ }G^{\min}_{V_1W_1}\left(\frac{\pi}{2}-\gamma\right)=G^{\min}_{V_1W_1}\left(\frac{\pi}{2}+\gamma\right)
	\end{equation*}
	\item [c)] When $\gamma$ increases from $\pi$ to $\frac{3\pi}{2}$, the process in a) repeats.
	\item [d)] When $\gamma$ increases from $\frac{3\pi}{2}$ to $2\pi$, the process in b) repeats.
\end{itemize}

The following subsections contain important considerations related to the previous description of the oscillating term.

\subsubsection{$\gamma_{1}(\widehat{y}_0)-\theta_{1}$ odd multiple of $\frac{\pi}{2}$} \label{ompi/2}

The condition $\gamma_{1}(\widehat{y}_0)-\theta_{1}$ is an odd multiple of $\frac{\pi}{2}$ is equivalent to the vector $\widehat{w}^{(1)}\widehat{y}_{0}$ being orthogonal  in the complex plane to the major semi-axis of the ellipse $\mathcal{E}_1$ of Proposition \ref{euclidean}, i.e., $c_1(\widehat{y}_0)=0$ by (\ref{whatu})-(\ref{whatu2}), namely $y_0$ is orthogonal to the first right singular vector of the matrix $R_1$ of Proposition \ref{euclidean}.

When $y_0$ is orthogonal to the first singular vector of $R_1$, we have
\begin{equation*}
\max\limits_{t\in\mathbb{R}} \mathrm{OT}\left(t,y_0\right)=a^{\max\max}_{V_1W_1}\text{\ and\ }
\min\limits_{t\in\mathbb{R}} \mathrm{OT}\left(t,y_0\right)=a^{\min\min}_{V_1W_1}  \label{minmin}
\end{equation*}
with $\max\limits_{t\in\mathbb{R}} \mathrm{OT}\left(t,y_0\right)$ at its maximum value and $\min\limits_{t\in\mathbb{R}} \mathrm{OT}\left(t,y_0\right)$ at its minimum value, as $y_0$ varies.

\subsubsection{$\gamma_{1}(\widehat{y}_0)-\theta_{1}$ is a multiple of $\pi$} \label{sssection}
The condition $\gamma_{1}(\widehat{y}_0)-\theta_{1}$ is a multiple of $\pi$ is equivalent to the vector $\widehat{w}^{(1)}\widehat{y}_{0}$ being parallel in the complex plane to the major semi-axis of the ellipse $\mathcal{E}_1$ of Proposition \ref{euclidean}, i.e., $d_1(\widehat{y}_0)=0$ by (\ref{whatu})-(\ref{whatu2}), namely $y_0$ is orthogonal to the second right singular vector of the matrix $R_1$ of Proposition \ref{euclidean}.

When $y_0$ is orthogonal to the second singular vector of $R_1$, we have
\begin{equation*}
\max\limits_{t\in\mathbb{R}}\mathrm{OT}\left(t,y_0\right)=a^{\min\max}_{V_1W_1} \text{ and }
\min\limits_{t\in\mathbb{R}}\mathrm{OT}\left(t,y_0\right)
=a^{\max\min}_{V_1W_1} \label{maxmin}
\end{equation*}
with $\max\limits_{t\in\mathbb{R}} \mathrm{OT}\left(t,y_0\right)$ at its minimum value and $\min\limits_{t\in\mathbb{R}} \mathrm{OT}\left(t,y_0\right)$ at its maximum value, as $y_0$ varies.

\subsubsection{$V_1$ not close to $1$}\label{V1notclose}
If $V_1$ is not close to $1$, then, for any $y_0$, $\mathrm{OT}(t,y_0)$  does not assume large values,  as $t$ varies. In fact, we have
$$
\max\limits_{t\in\mathbb{R}}\mathrm{OT}(t,y_0)\leq a^{\max\max}_{V_1W_1}=\sqrt{\frac{(1+V_1)(1+W_1)}{2(1-V_1)}}.
$$
If $V_1$ is not close to $1$ or $W_1$ is not close to $1$, then, for any $y_0$, $\mathrm{OT}(t,y_0)$ does not assume small values, as $t$ varies. In fact, we have
$$
\min\limits_{t\in\mathbb{R}}\mathrm{OT}(t,y_0)\geq a^{\min\min}_{V_1W_1}=\left\{
\begin{array}{l}
	\sqrt{\frac{(1-V_1)(1+W_1)}{2(1+V_1)}}\text{\ if\ }V_1\leq W_1\\
	\\
	\sqrt{\frac{1-W_1}{2}}\text{\ if\ }V_1\geq W_1.
\end{array}
\right. 
$$

As a consequence, if $V_1$ is not close to $1$, then the asymptotic condition number
$$
K_\infty(t,y_0)=\mathrm{OSF}(y_0)\cdot \mathrm{OT}(t,y_0),
$$
has the same order of magnitude as the oscillation scale factor $\mathrm{OSF}(y_0)$.

\subsubsection{Critical case} \label{critcase2}

There is a critical case: $V_1$ close to one.

If $V_1$ is close to $1$, then, there exists an initial value $y_0$ such that $\mathrm{OT}(t,y_0)$  assumes large values,  as $t$ varies, and, as a consequence, the asymptotic condition number $K_\infty(t,y_0)$  assumes large values,  as $t$ varies (recall that  $\mathrm{OSF}(y_0)$ is never small, see Subsection \ref{SectionOSFy0}). 

In fact, for $y_0$ orthogonal to the first right singular vector of the matrix $R_1$ we have (see Subsection \ref{ompi/2})
$$
\max\limits_{t\in\mathbb{R}}\mathrm{OT}\left(t,y_0\right)=a^{\max\max}_{V_1W_1}=\sqrt{\frac{(1+V_1)(1+W_1)}{2(1-V_1)}}.
$$
For such $y_0$, see Subsection \ref{69} a) below.

\subsubsection{Supercritical case}\label{scc}

There is a super-critical case: $\frac{1-V_1}{1-W_1}$ small.

If the ratio $\frac{1-V_1}{1-W_1}$ is small, then, for any initial value $y_0$, $\mathrm{OT}(t,y_0)$  assumes large values,  as $t$ varies, and, as a consequence, the asymptotic condition number $K_\infty(t,y_0)$  assumes large values,  as $t$ varies.

In fact, for any $y_0$, we have 
$$
\max\limits_{t\in\mathbb{R}}\mathrm{OT}\left(t,y_0\right)\geq a^{\min\max}_{V_1W_1}=\left\{
\begin{array}{l}
	\sqrt{\frac{(1-V_1^2)(1+W_1)}{2(1-Q_1V_1)}}\text{\ if\ }Q_1\leq 1\\
	\\
	\sqrt{\frac{(1+V_1)(1-W_1)}{2(1-V_1)}}\text{\ if\ }Q_1\geq 1,
\end{array}
\right.
$$
and, for
$$
\frac{1-V_1}{1-W_1}\leq \frac{1}{1+W_1},
$$
the case $Q_1\geq 1$ holds.

\subsubsection{The critical, but not supercritical, case} \label{cnscc}
Observe that $\frac{1-V_1}{1-W_1}$ small implies $V_1$ close to $1$. Therefore the supercritical case is a subcase of the critical case.

Suppose the critical, but not supercritical, case, i.e., $V_1$ is close to $1$, but $\frac{1-V_1}{1-W_1}$ is not small. In this case, along with existence of $y_0$ such that $\mathrm{OT}(t,y_0)$  assume large values,  as $t$ varies, there also exists $y_0$ such that $\mathrm{OT}(t,y_0)$  does not assume large values,  as $t$ varies.

In fact, for $y_0$ orthogonal to the second right singular vector of the matrix $R_1$ we have (see Subsection \ref{sssection})
$$
\max\limits_{t\in\mathbb{R}}\mathrm{OT}\left(t,y_0\right)= a^{\min\max}_{V_1W_1}=\left\{
\begin{array}{l}
	\sqrt{\frac{(1-V_1^2)(1+W_1)}{2(1-Q_1V_1)}}\text{\ if\ }Q_1\leq 1\\
	\\
	\sqrt{\frac{(1+V_1)(1-W_1)}{2(1-V_1)}}\text{\ if\ }Q_1\geq 1.
\end{array}
\right.
$$
In the case $Q_1\leq 1$, we have
$$
\sqrt{\frac{(1-V_1^2)(1+W_1)}{2(1-Q_1V_1)}}\leq\sqrt{\frac{(1+V_1)(1+W_1)}{2}}
$$
and then $a^{\min\max}_{V_1W_1}$ is not large. In the case $Q_1\geq 1$, $a^{\min\max}_{V_1W_1}$ is not large if the ratio  $\frac{1-V_1}{1-W_1}$ is not small.

For such $y_0$, see Subsection \ref{69} b) below.

\subsubsection{Specific cases} \label{69}
Two specific cases are worth considering.
\begin{itemize}
	\item [a)] Consider $y_0$ orthogonal to the first right singular vector of the matrix $R_1$ of Proposition \ref{euclidean} (i.e., $c_1(\widehat{y}_0)=0$).
	
	We have
	\begin{equation*}
		\mathrm{OSF}(y_0)=\sqrt{\frac{2}{1-W_1}}\cdot \frac{1}{\vert d_1(\widehat{y}_0)\vert}
	\end{equation*}
	by (\ref{OSF2}). Moreover, by recalling Subsection \ref{ompi/2}, we have
	\begin{eqnarray*}
		&&\max\limits_{t\in\mathbb{R}}\mathrm{OT}\left(t,y_0\right)=a^{\max\max}_{V_1W_1}=\sqrt{\frac{(1+V_1)(1+W_1)}{2(1-V_1)}}\\
		&&\min\limits_{t\in\mathbb{R}}\mathrm{OT}\left(t,y_0\right)
		=a^{\min\min}_{V_1W_1}=\left\{
		\begin{array}{l}
			\sqrt{\frac{(1-V_1)(1+W_1)}{2(1+V_1)}}\text{\ if\ }V_1\leq W_1\\
			\\
			\sqrt{\frac{1-W_1}{2}}\text{\ if\ }V_1\geq W_1.
		\end{array}
		\right.  
	\end{eqnarray*}
Regarding the asymptotic condition number
	$
	K_\infty(t,y_0),
	$
	we have
	$$
	\max\limits_{t\in  \mathbb{R}}K_\infty(t,y_0)=\sqrt{\frac{(1+V_1)(1+W_1)}{(1-V_1)(1-W_1)}}\cdot \frac{1}{\vert d_1(\widehat{y}_0)\vert}
	$$
	and
	$$
	\min\limits_{t\in  \mathbb{R}}K_\infty(t,y_0)=\left\{
	\begin{array}{l}
		\sqrt{\frac{(1-V_1)(1+W_1)}{(1+V_1)(1-W_1)}}\cdot \frac{1}{\vert d_1(\widehat{y}_0)\vert}\text{\ if\ }V_1\leq W_1\\
		\\
		\frac{1}{\vert d_1(\widehat{y}_0)\vert}\text{\ if\ }V_1\geq W_1.
	\end{array}
	\right.\ 
	$$
	Compare to (\ref{Kyz1})-(\ref{Kyz11}). 
	
		We can conclude the following.
		\begin{itemize}
			\item It is expected that $\mathrm{OSF}(y_0)$ is large if and only if $W_1$ is close to $1$ (since $\frac{1}{\vert d_1(\widehat{y}_0)\vert}$ is expected to be not large).
			\item $\mathrm{OT}(t,y_0)$ assumes large values, as $t$ varies, if and only if $V_1$ is close to $1$.
			\item $\mathrm{OT}(t,y_0)$ assumes small values, as $t$ varies, if and only if $V_1$ is close to $1$ and $W_1$ is close to $1$. 
			\item $K_\infty(t,y_0)$ assumes large values, as $t$ varies, if $V_1$ is close to $1$ or $W_1$ is close to $1$ (we have $\vert d_1(\widehat{y}_0)\vert\leq 1$).
			\item It is expected that $K_\infty(t,y_0)$ assumes large values, as $t$ varies, only if $V_1$ is close to $1$ or $W_1$ is close to $1$ (since $\frac{1}{\vert d_1(\widehat{y}_0)\vert}$ is expected to be not large).
			\item $\min\limits_{t\in\mathbb{R}}K_\infty\left(t,y_0\right)\geq \frac{1}{\vert d_1(\widehat{y}_0)\vert}$.
			\item $\min\limits_{t\in\mathbb{R}}K_\infty\left(t,y_0\right)$ is large, i.e. $K_\infty(t,y_0)$ always assumes large values, as $t$ varies, if the ratio $\frac{1-V_1}{1-W_1}$ is large.
		\end{itemize}
	\item [b)] Consider $y_0$ orthogonal to the second right singular vector of $R_1$ (i.e., $d_1(\widehat{y}_0)=0$).
	
	We have
	\begin{equation*}
		\mathrm{OSF}(y_0)=\sqrt{\frac{2}{1+W_1}}\cdot \frac{1}{\vert c_1(\widehat{y}_0)\vert}\
	\end{equation*}
	by (\ref{OSF2}). Moreover, by recalling Subsection \ref{sssection}, we have
	\begin{eqnarray*}
		&&\max\limits_{t\in\mathbb{R}}\mathrm{OT}\left(t,y_0\right)=a^{\min\max}_{V_1W_1}=\left\{
		\begin{array}{l}
			\sqrt{\frac{(1-V_1^2)(1+W_1)}{2(1-Q_1V_1)}}\text{\ if\ }Q_1\leq 1\\
			\\
			\sqrt{\frac{(1+V_1)(1-W_1)}{2(1-V_1)}}\text{\ if\ }Q_1\geq 1,
		\end{array}
		\right. \\
		&&\min\limits_{t\in\mathbb{R}}\mathrm{OT}\left(t,y_0\right)
		=a^{\max\min}_{V_1,W_1}=\sqrt{\frac{1+W_1}{2}}  
	\end{eqnarray*}
Regarding the asymptotic condition number $K_\infty(t,y_0)$, we have
	$$
	\quad\quad\quad \quad\quad\quad \max\limits_{t\in  \mathbb{R}}K_\infty(t,y_0)=\left\{
	\begin{array}{l}
		\sqrt{\frac{1-V_1^2}{1-Q_1V_1}}\cdot \frac{1}{\vert c_1(\widehat{y}_0)\vert}\text{\ if $Q_1\leq 1$}\\
		\\
		\sqrt{\frac{(1+V_1)(1-W_1)}{(1-V_1)(1+W_1)}}\cdot \frac{1}{\vert c_1(\widehat{y}_0)\vert}\text{\ if $Q_1\geq 1$}
     \end{array}\ 
	\right.
	$$
	and
	$$
	\quad\quad\quad \quad\quad\quad \min\limits_{t\in  \mathbb{R}}K_\infty(t,y_0)=\frac{1}{\vert c_1(\widehat{y}_0)\vert}.
	$$
	Compare to (\ref{Kyz2})-(\ref{Kyz22}).
	
		We can conclude the following.
	\begin{itemize}
		\item It is expected that $\mathrm{OSF}(y_0)$ is not small;
		\item $\mathrm{OT}(t,y_0)$ assumes large values, as $t$ varies, if and only if $\frac{1-V_1}{1-W_1}$ is small.
		\item $\mathrm{OT}(t,y_0)$ does not assume small values, as $t$ varies.
		\item $K_\infty(t,y_0)$ assumes large values, as $t$ varies, if $\frac{1-V_1}{1-W_1}$ is small. 
		\item It is expected that $K_\infty(t,y_0)$ assumes large values, as $t$ varies, only if $\frac{1-V_1}{1-W_1}$ is small.
		\end{itemize}
\end{itemize}

\subsection{Numerical examples} \label{test}
In this subsection, we consider three examples of a real ODE (\ref{ode}) of dimension 3, where the matrix $A$ has a pair of imaginary conjugate simple eigenvalues as rightmost eigenvalues. We are interested in illustrating the behavior of the asymptotic pointwise condition number $K_\infty(t,y_0)$, product of the oscillating scale factor $\mathrm{OSF}(y_0)$ and the oscillating term $\mathrm{OT}(t,y_0)$, when the Euclidean norm is used as vector norm.

In all three examples, we consider as initial value $y_0$ the second and first singular vectors of the matrix $R_1$ of Proposition \ref{euclidean}, therefore considering the specific cases a) and b), respectively, of Subsection \ref{69}. 

\bigskip

\underline{First example}

\bigskip

The first example has both $V_1$ and $W_1$ not close to $1$.

Thus, for $y_0$ belonging to the subspace spanned by the first and second right singular vectors of the matrix $R_1$, $\mathrm{OSF}(y_0)$ is not large (see (\ref{OSF2})). Moreover, for all $y_0$, $\mathrm{OT}(t,y_0)$ does not assume neither large values nor small values, as $t$ varies (see Subsection \ref{V1notclose}). As a consequence, for $y_0$ belonging to the subspace spanned by the first and second right singular vectors of the matrix $R_1$, $K_\infty(t,y_0)$ does not assume large values, as $t$ varies.

Consider
\begin{equation*}
	A=\left[ 
	\begin{array}{rrr}
		1 & 0 & -2 \\ 
		1 &  -1 & 0\\ 
		1 & 0 & -1
	\end{array}
	\right]=V\mathrm{diag}(\mathrm{i},-\mathrm{i},-1)V^{-1},  \label{AV1}
\end{equation*}
with
\begin{equation*}
	V=\left[ 
	\begin{array}{rrr}
		1+\mathrm{i} & 1-\mathrm{i} & 0 \\ 
		1 & 1 & 1 \\ 
		1 & 1 & 0
	\end{array}
	\right]. 
\end{equation*}
We have $V_1=0.7071$ and $W_1=0.7454$. The matrix $R_1$ of Proposition \ref{euclidean} is  
\begin{equation*}
	R_1=\left[
	\begin{array}{c}
		\mathrm{Re}(\widehat{w}^{(1)}) \\ 
		\mathrm{Im}(\widehat{w}^{(1)})
	\end{array}
	\right]=\left[ 
	\begin{array}{rrr}
		0 & 0     & 0.5774 \\ 
		-0.5744 & 0 & 0.5774 
	\end{array}
	\right]
\end{equation*}
and has $(0.5257,0,-0.8507)$ and $(0.8507,0,0.5257)$ as first and second, respectively, right singular vectors.

\begin{figure}[tbp]
	\centering
	\begin{subfigure}[b]{1\textwidth}
		\centering
		\includegraphics[width=0.8\textwidth]{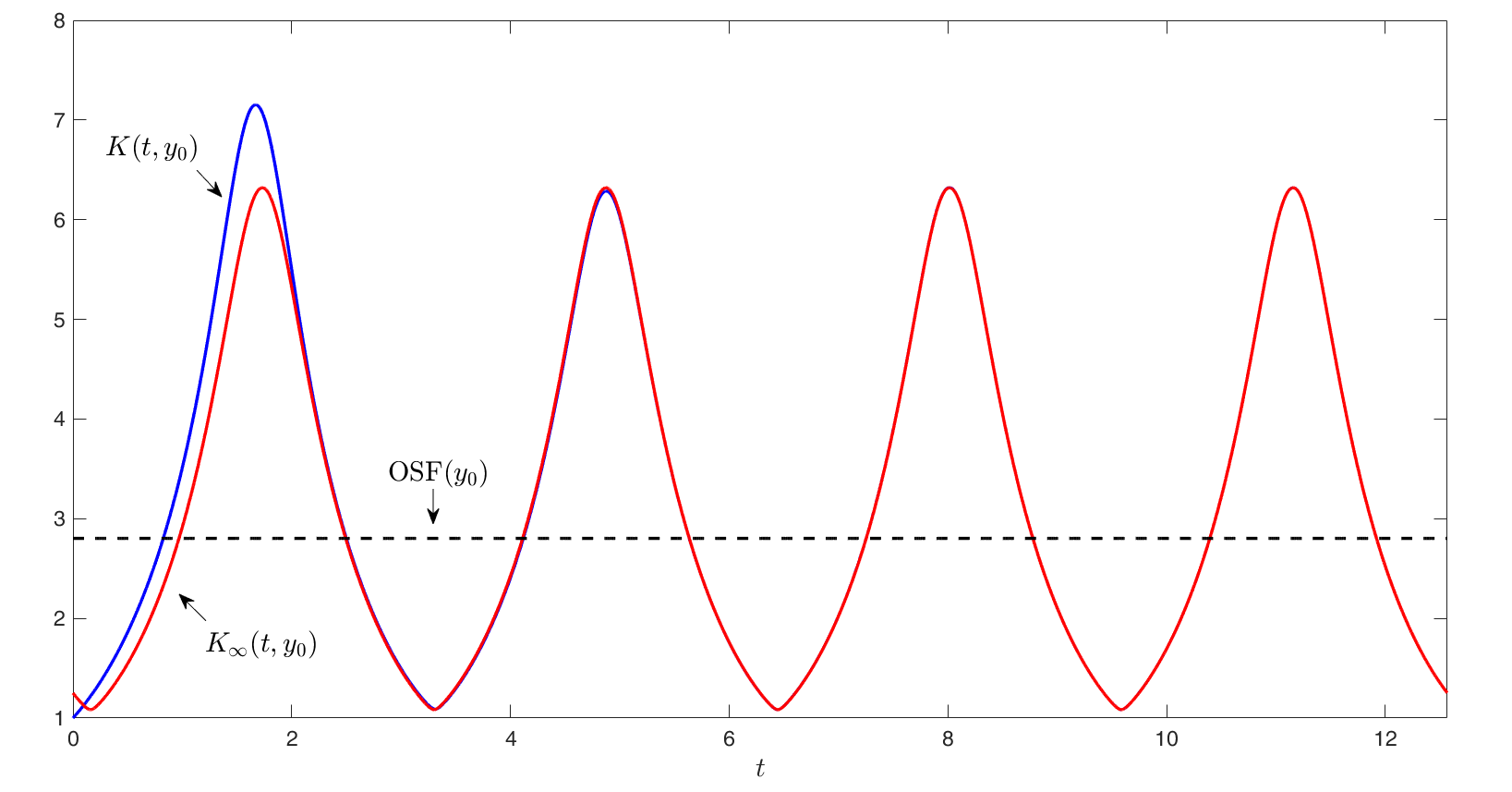}
		\caption{$y_0=(0.8507,0,0.5257)$ second right singular vector of $R_1$.} 
	\end{subfigure}
	\hfill
	\begin{subfigure}[b]{1\textwidth}
		\centering
		\includegraphics[width=0.8\textwidth]{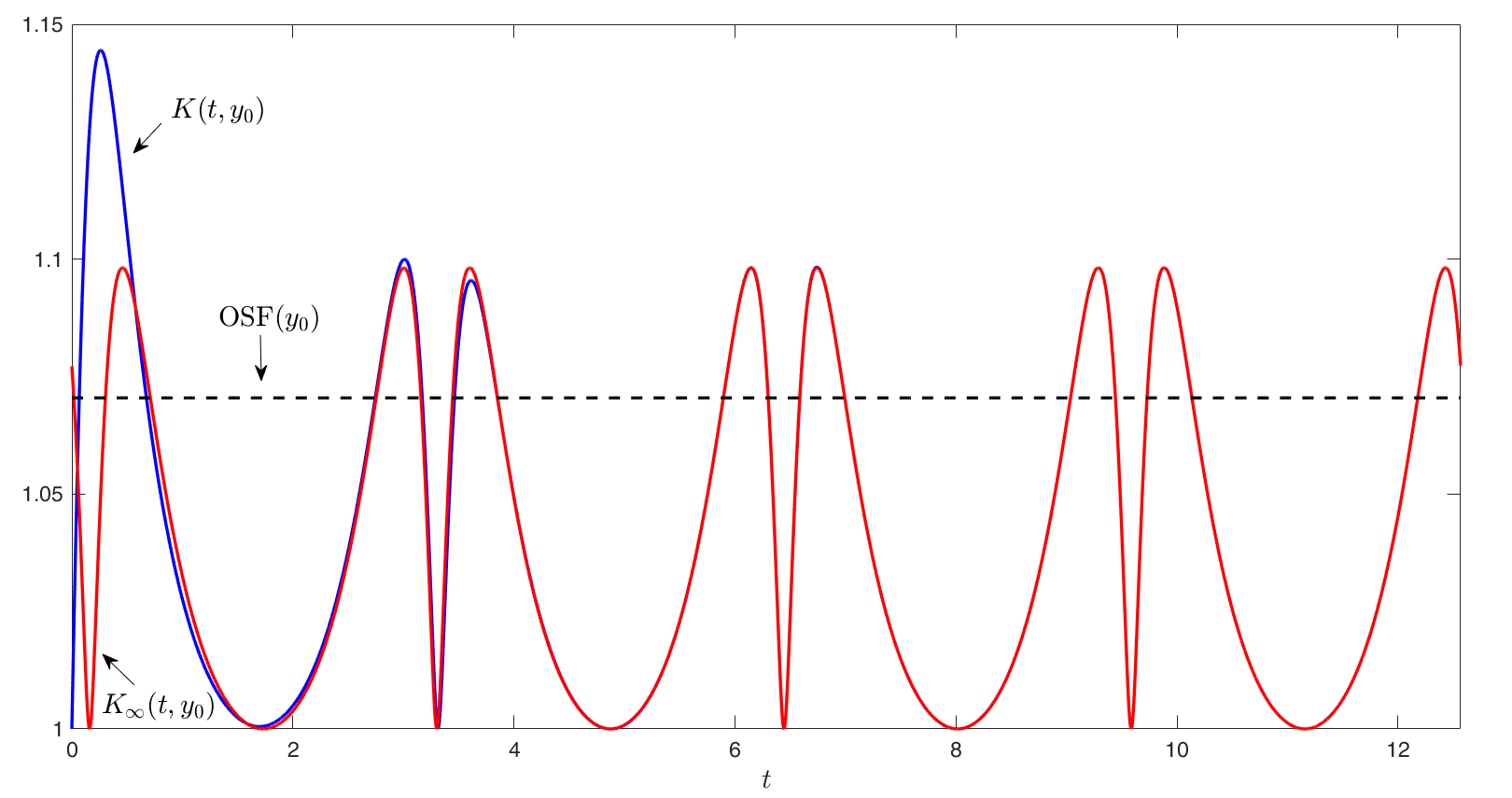}
		\caption{$y_0=(0.5257,0,-0.8507)$ first right singular vecto of $R_1$.}
	\end{subfigure}  
	\caption{Condition numbers $K(t,y_0)$ and $K_\infty(t,y_0)$, $t\in[0,4\pi]$, and oscillation scale factor $\mathrm{OSF}(y_0)$, for the first example.}
	\label{FigureV0}
\end{figure}

In Figure \ref{FigureV0}, for $y_0$ equal to the second (upper part) and  the first (lower part) right singular vectors, we see the condition numbers $K(t, y_0)$ and $K_\infty(t, y_0)$ along with the oscillation scale factor $\mathrm{OSF}(y_0)$ for $t \in [0, 4\pi]$, covering four periods of the function $K_\infty(t, y_0)$.

Table \ref{table1} reports the key quantities for Figure \ref{FigureV0}, computed using Subsection \ref{69}.

\begin{table}
	\begin{tabular}{|l|l|l|}
		\hline
		Key quantity & Upper part & Lower part\\
		\hline
		$\mathrm{OSF}(y_0)$ & $\sqrt{\frac{2}{1-W_1}}=2.8025$ & $\sqrt{\frac{2}{1+W_1}}=1.0705$ \\
		\hline 
		$\max\limits_{t\in  \mathbb{R}}\mathrm{OT}(t,y_0)$ & \makecell[l]{$a^{\mathrm{maxmax}}_{V_1W_1}$\\$= \sqrt{\frac{(1+V_1)(1+W_1)}{2(1-V_1)}}=2.2553$}   & 
		\makecell[l]{$a^{\mathrm{\min\max}}_{V_1W_1}$\\$=\sqrt{\frac{(1-V_1^2)(1+W_1)}{2(1-Q_1V_1)}}= 1.0259$\\$\text{(since $Q_1=0.8279\leq 1$)}$}\\
		\hline
		$\min\limits_{t\in  \mathbb{R}}\mathrm{OT}(t,y_0)$ & \makecell[l]{$a^{\mathrm{minmin}}_{V_1W_1}$\\$=\sqrt{\frac{(1-V_1)(1+W_1)}{2(1+V_1)}}= 0.3869$\\$  \text{(since $V_1\leq W_1$)}$}   & \makecell[l]{$a^{\mathrm{\max\min}}_{V_1W_1}$\\$=\sqrt{\frac{1+W_1}{2}}=0.9342$}\\
		\hline
		$\max\limits_{t\in  \mathbb{R}}K_\infty(t,y_0)$ & \makecell[l]{$\mathrm{OSF}(y_0)\cdot \max\limits_{t\in  \mathbb{R}}\mathrm{OT}(t,y_0)$\\$=6.3205$}& \makecell[l]{$\mathrm{OSF}(y_0)\cdot \max\limits_{t\in  \mathbb{R}}\mathrm{OT}(t,y_0)$\\$=1.0982$}\\
		\hline
		$\min\limits_{t\in  \mathbb{R}}K_\infty(t,y_0)$ & \makecell[l]{$\mathrm{OSF}(y_0)\cdot \min\limits_{t\in  \mathbb{R}}\mathrm{OT}(t,y_0)$\\$=1.0844$}& \makecell[l]{$\mathrm{OSF}(y_0)\cdot \min\limits_{t\in  \mathbb{R}}\mathrm{OT}(t,y_0)$\\$=1$} \\
		\hline
	\end{tabular}
	\caption{Key quantities for Figure \ref{FigureV0}}
	\label{table1}
\end{table}

\bigskip

\underline{Second example}

\bigskip

The second example has both $V_1$ and $W_1$ close to $1$, with the ratio $\frac{1-V_1}{1-W_1}$ not small. This is the critical, but not supercritical, case (see Subsection \ref{cnscc}).

Thus, for $y_0$ equal to the second right singular vector of the matrix $R_1$, $\mathrm{OSF}(y_0)$ is large and $\mathrm{OT}(t,y_0)$ assumes both large and small values, as $t$ varies (see  point a) of Subsection \ref{69}). As a consequence, for such $y_0$, $K_\infty(t,y_0)$ assumes large values, as $t$ varies.

Moreover, for $y_0$ equal to the first right singular vector of the matrix $R_1$, $\mathrm{OSF}(y_0)$ is not large and $\mathrm{OT}(t,y_0)$ assumes neither  large values nor small values, as $t$ varies (see  point b) of Subsection \ref{69}). As a consequence, for such $y_0$, $K_\infty(t,y_0)$ does not assume large values, as $t$ varies.

Consider
\begin{equation*}
	A=\left[ 
	\begin{array}{rrr}
		-1 & 20 & -20 \\ 
		0 &  19 & -20\\ 
		0 & 18.1 & -19
	\end{array}
	\right]=V\mathrm{diag}(\mathrm{i},-\mathrm{i},-1)V^{-1},  \label{AV2}
\end{equation*}
with
\begin{equation*}
	V=\left[ 
	\begin{array}{rrr}
		1+\mathrm{i} & 1-\mathrm{i} & 1 \\ 
		1+\mathrm{i} & 1-\mathrm{i} & 0 \\ 
		1+0.9\mathrm{i} & 1-0.9\mathrm{i} & 0
	\end{array}
	\right]. 
\end{equation*}
We have $V_1=0.9988$, $W_1=0.9986$ and $\frac{1-V_1}{1-W_1}=0.8600$. The matrix $R_1$ is
\begin{equation*}
	R_1=\left[
	\begin{array}{c}
		\mathrm{Re}(\widehat{w}^{(1)}) \\ 
		\mathrm{Im}(\widehat{w}^{(1)})
	\end{array}
	\right]=\left[ 
	\begin{array}{rrr}
		0 & -0.4611     & 0.5123 \\ 
		0 & -0.5123 & 0.5123 
	\end{array}
	\right].
\end{equation*}
and has $(0,0.6892,-0.7245)$ and $(0,-0.7245,-0.6892)$ as first and second, respectively, right singular vectors.

\begin{figure}[tbp]
	\centering
	\begin{subfigure}[b]{1\textwidth}
		\centering
		\includegraphics[width=0.8\textwidth]{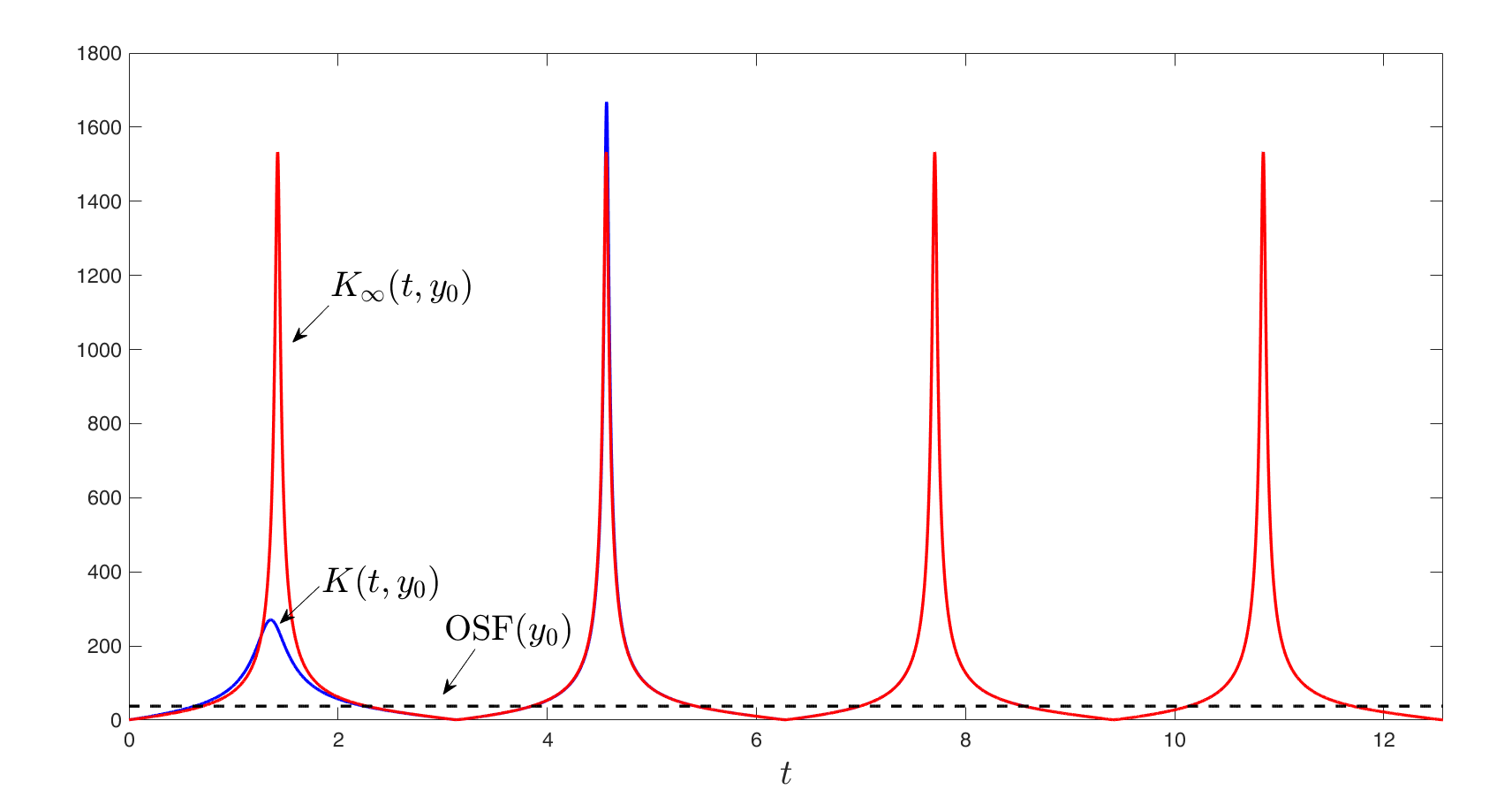}
		\caption{$y_0=(0,-0.7245,-0.6892)$ second right singular vector of $R_1$.} 
	\end{subfigure}
	\hfill
	\begin{subfigure}[b]{1\textwidth}
		\centering
		\includegraphics[width=0.8\textwidth]{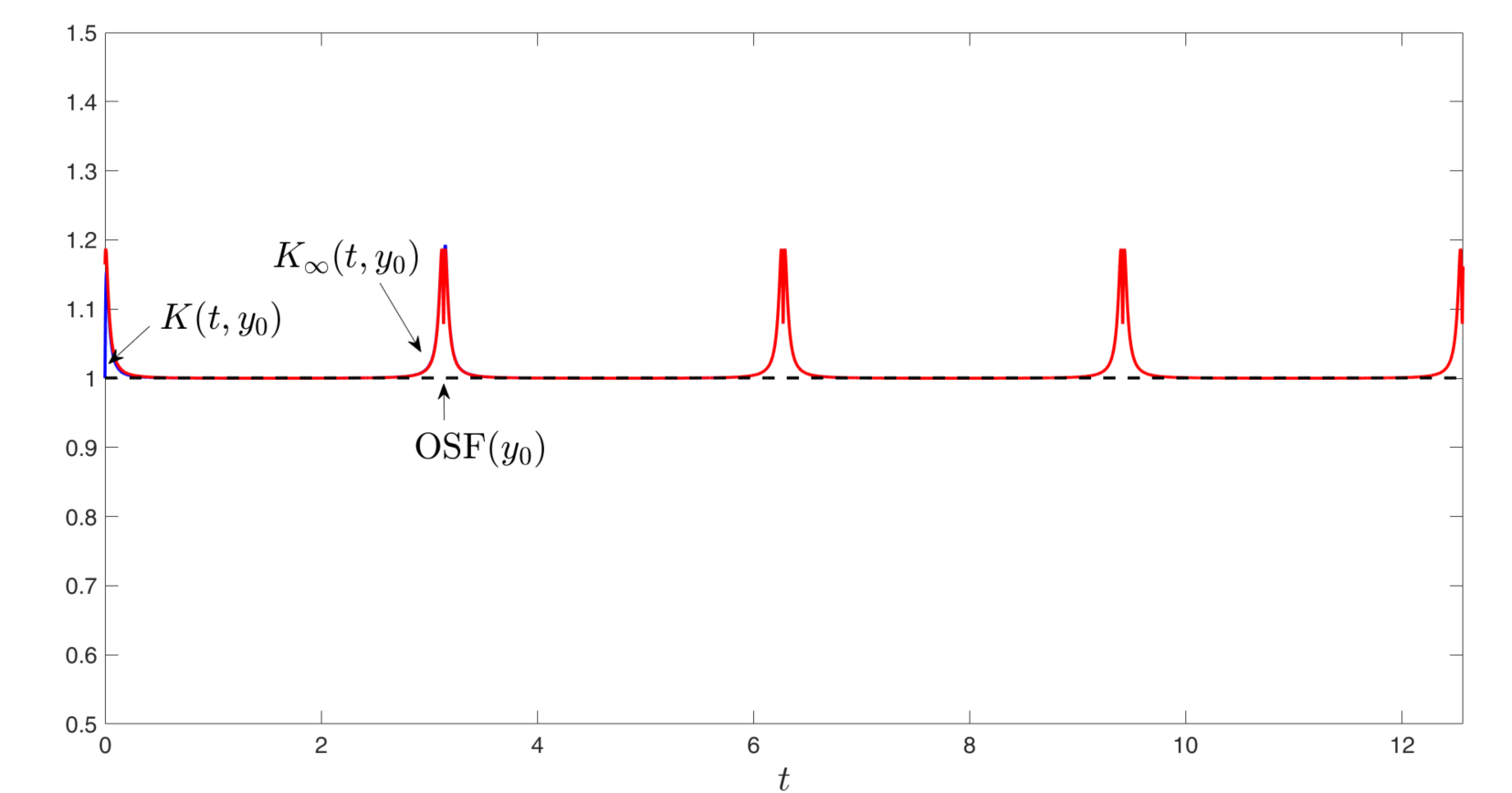}
		\caption{$y_0=(0,0.6892,-0.7245)$ first right singular vector of $R_1$.}
	\end{subfigure}  
	\caption{Condition numbers $K(t,y_0)$ and $K_\infty(t,y_0)$, $t\in[0,4\pi]$, and oscillation scale factor $\mathrm{OSF}(y_0)$, for the second example.}
	\label{FigureV1}
\end{figure}

Figure \ref{FigureV1}, similarly to Figure \ref{FigureV0} for the first example, shows $K(t, y_0)$, $K_\infty(t, y_0)$ and $\mathrm{OSF}(y_0)$, $t \in [0, 4\pi]$, for $y_0$ equal to the second (upper part) and first (lower part) right singular vectors of $R_1$.

Table \ref{table2} reports the key quantities for Figure \ref{FigureV1}.

\begin{table}
	\begin{tabular}{|l|l|l|}
		\hline
		Key quantity & Upper part & Lower part\\
		\hline
		$\mathrm{OSF}(y_0)$ & $\sqrt{\frac{2}{1-W_1}}=38.1$ & $\sqrt{\frac{2}{1+W_1}}=1.0003$ \\
		\hline 
		$\max\limits_{t\in  \mathbb{R}}\mathrm{OT}(t,y_0)$ & \makecell[l]{$a^{\mathrm{maxmax}}_{V_1W_1}$\\$= \sqrt{\frac{(1+V_1)(1+W_1)}{2(1-V_1)}}=41.0$}   & 
		\makecell[l]{$a^{\mathrm{\min\max}}_{V_1W_1}$\\$=\sqrt{\frac{(1-V_1^2)(1+W_1)}{2(1-Q_1V_1)}}= 1.1869$\\$\text{(since $Q_1=0.9995\leq 1$)}$}\\
		\hline
		$\min\limits_{t\in  \mathbb{R}}\mathrm{OT}(t,y_0)$ & \makecell[l]{$a^{\mathrm{minmin}}_{V_1W_1}$\\$=\sqrt{\frac{1-W_1}{2}}=0.0263$\\$  \text{(since $V_1\geq W_1$)}$}   & \makecell[l]{$a^{\mathrm{\max\min}}_{V_1W_1}$\\$=\sqrt{\frac{1+W_1}{2}}=0.99965$}\\
		\hline
		$\max\limits_{t\in  \mathbb{R}}K_\infty(t,y_0)$ & \makecell[l]{$\mathrm{OSF}(y_0)\cdot \max\limits_{t\in  \mathbb{R}}\mathrm{OT}(t,y_0)$\\$=1563$}& \makecell[l]{$\mathrm{OSF}(y_0)\cdot \max\limits_{t\in  \mathbb{R}}\mathrm{OT}(t,y_0)$\\=1.1873$$}\\
		\hline
		$\min\limits_{t\in  \mathbb{R}}K_\infty(t,y_0)$ & \makecell[l]{$\mathrm{OSF}(y_0)\cdot \min\limits_{t\in  \mathbb{R}}\mathrm{OT}(t,y_0)$\\$=1$}& \makecell[l]{$\mathrm{OSF}(y_0)\cdot \min\limits_{t\in  \mathbb{R}}\mathrm{OT}(t,y_0)$\\$=1$} \\
		\hline
	\end{tabular}
	\caption{Key quantities for Figure \ref{FigureV1}}
	\label{table2}
\end{table}

\bigskip

\bigskip

\underline{Third example}

\bigskip

The third example has $V_1$ close to $1$ and $W_1=0$. This is the supercritical case (see Subsection \ref{scc}).

We have $\mathrm{OSF}(y_0)\geq\sqrt{2}$, for all $y_0$, and $\mathrm{OSF}(y_0)=\sqrt{2}$
for $y_0$ belonging to the subspace spanned by the first and second right singular vectors of the matrix $R_1$ (see (\ref{OSF2})). Moreover, Remark \ref{V1W1zero} applies to this example and so, as $y_0$ varies, the functions $\mathrm{OT}(t,y_0)$ of $t$ are all time shifts of the same function oscillating monotonically between the extreme values
\begin{eqnarray*}
	\max\limits_{t\in\mathbb{R}}\mathrm{OT}\left(t,y_0\right)=\sqrt{\frac{1+V_1}{2(1-V_1)}}\text{\ \ and\ \ }\min\limits_{t\in\mathbb{R}}\mathrm{OT}\left(t,y_0\right)=\sqrt{\frac{1}{2}}.
\end{eqnarray*}
Therefore, for all $y_0$, $\mathrm{OT}(t,y_0)$ assumes large values, as $t$ varies, and consequently $K_\infty(t,y_0)$ assumes large values, as $t$ varies.

Consider
\begin{equation*}
	A=\left[ 
	\begin{array}{rrr}
		0 & -1 & 0 \\ 
		1 &  0 & 0\\ 
		200 & -200 & -1
	\end{array}
	\right]=V\mathrm{diag}(\mathrm{i},-\mathrm{i},-1)V^{-1},  \label{AV3}
\end{equation*}
with
\begin{equation*}
	V=\left[ 
	\begin{array}{rrr}
		0.5 & 0.5 & 0 \\ 
		-0.5\mathrm{i} & 0.5\mathrm{i} & 0 \\ 
		100 & 100 & 1
	\end{array}
	\right]. 
\end{equation*}
We have $V_1=0.99995$ and $W_1=0$. The matrix $R_1$ is  
\begin{equation*}
	R_1=\left[
	\begin{array}{c}
		\mathrm{Re}(\widehat{w}^{(1)}) \\ 
		\mathrm{Im}(\widehat{w}^{(1)})
	\end{array}
	\right]=\left[ 
	\begin{array}{rrr}
		0.7071 & 0     & 0 \\ 
		0 & 0.7071 & 0 
	\end{array}
	\right]
\end{equation*}
and has $(1,0,0)$ and $(0,1,0)$ as first and second, respectively, right singular vectors.

\begin{figure}[tbp]
	\centering
	\begin{subfigure}[b]{1\textwidth}
		\centering
		\includegraphics[width=0.85\textwidth]{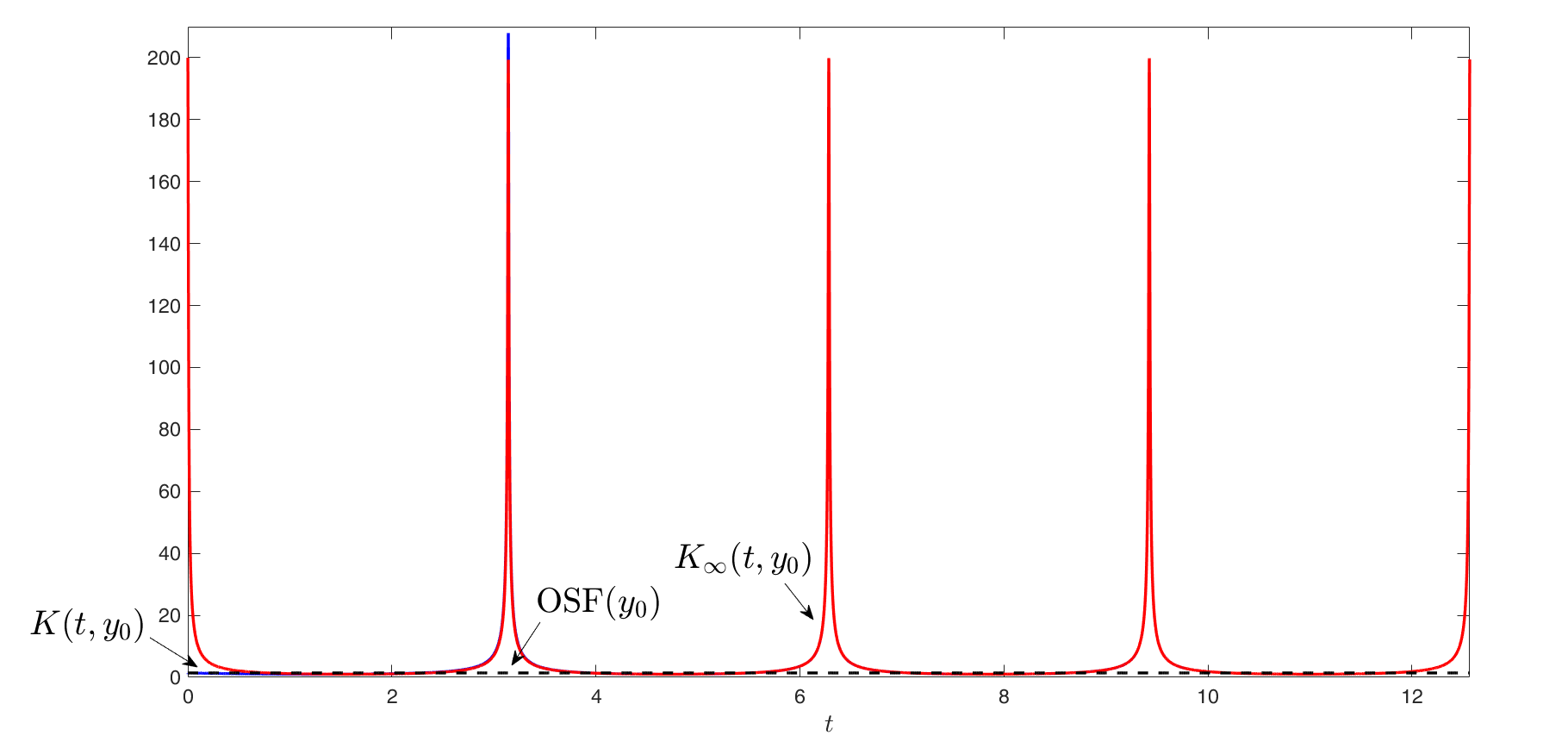}
		\caption{$y_0=(0,1,0)$ second right singular vector of $R_1$.} 
	\end{subfigure}
	\hfill
	\begin{subfigure}[b]{1\textwidth}
		\centering
		\includegraphics[width=0.8\textwidth]{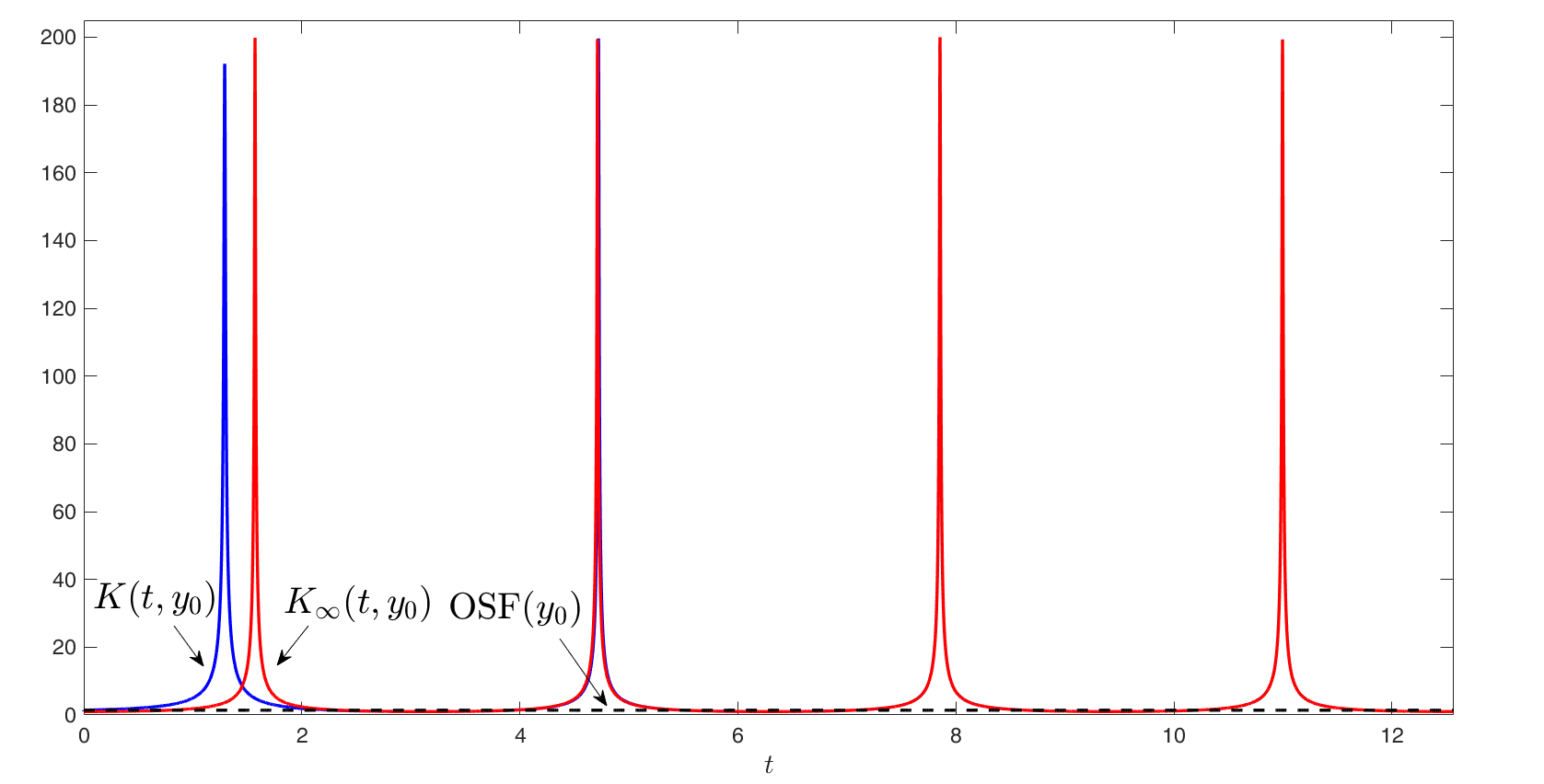}
		\caption{$y_0=(1,0,0)$ first right singular vectorof $R_1$.}
	\end{subfigure}  
	\caption{Condition numbers $K(t,y_0)$ and $K_\infty(t,y_0)$, $t\in[0,4\pi]$, and oscillation scale factor $\mathrm{OSF}(y_0)$, for the third example.}
	\label{FigureV2}
\end{figure}

Figure \ref{FigureV2}, similarly to Figure \ref{FigureV0} for the first example and Figure \ref{FigureV1} for the second example, shows $K(t, y_0)$, $K_\infty(t, y_0)$ and $\mathrm{OSF}(y_0)$, $t \in [0, 4\pi]$, for $y_0$ equal to the second (upper part) and first (lower part) right singular vectors of $R_1$.

Table \ref{table3} reports the key quantities for Figure \ref{FigureV2}.

\begin{table}
	\begin{tabular}{|l|l|l|}
		\hline
		Key quantity & Upper part & Lower part\\
		\hline
		$\mathrm{OSF}(y_0)$ & $\sqrt{2}$ & $\sqrt{2}$ \\
		\hline 
		$\max\limits_{t\in  \mathbb{R}}\mathrm{OT}(t,y_0)$ & $\sqrt{\frac{1+V_1}{2(1-V_1)}}=141.4231$   & $\sqrt{\frac{1+V_1}{2(1-V_1)}}=141.4231$\\
		\hline
		$\min\limits_{t\in  \mathbb{R}}\mathrm{OT}(t,y_0)$ & $\sqrt{\frac{1}{2}}$   & $\sqrt{\frac{1}{2}}$\\
		\hline
		$\max\limits_{t\in  \mathbb{R}}K_\infty(t,y_0)$ & \makecell[l]{$\mathrm{OSF}(y_0)\cdot \max\limits_{t\in  \mathbb{R}}\mathrm{OT}(t,y_0)$\\$=200.0025$}& \makecell[l]{$\mathrm{OSF}(y_0)\cdot \max\limits_{t\in  \mathbb{R}}\mathrm{OT}(t,y_0)$\\$=200.0025$}\\
		\hline
		$\min\limits_{t\in  \mathbb{R}}K_\infty(t,y_0)$ & \makecell[l]{$\mathrm{OSF}(y_0)\cdot \min\limits_{t\in  \mathbb{R}}\mathrm{OT}(t,y_0)$\\$=1$}& \makecell[l]{$\mathrm{OSF}(y_0)\cdot \min\limits_{t\in  \mathbb{R}}\mathrm{OT}(t,y_0)$\\$=1$} \\
		\hline
	\end{tabular}
	\caption{Key quantities for Figure \ref{FigureV2}}
	\label{table3}
\end{table}

\section{Dominance of the asymptotic condition numbers at a finite time}\label{Section6}

In this final section before the conclusion, we apply the results of Sections 2 and 3 in the context of Theorem \ref{two}, which establishes the dominance of the asymptotic condition numbers at finite time. This also informs how quickly the asymptotic behavior of the condition numbers sets in. 

Theorem \ref{two} states that
the dominance of the asymptotic condition numbers at a finite time is governed by how small $
\epsilon(t,\widehat{z}_0)$, $\epsilon(t,\widehat{y}_0)$ and $
\epsilon(t)$ are. The following
theorem addresses this question, where we consider a general vector $u\in\mathbb{R}^n$ such that $\widehat{w}^{(1)}u\neq 0$, as is the case  for $\widehat{z}_0$ and $\widehat{y}_0$.  

\begin{theorem}\label{onset}
		Suppose that, for any $j\in\{1,\ldots,q\}$, $\Lambda_j$ is simple real or simple complex. For $u\in\mathbb{R}^n$ such that $
	\widehat{w}^{(1)}u\neq 0$, we have 
	\begin{equation*}
		\epsilon(t,u)=\sum\limits_{j=2}^{q}\mathrm{e}^{\left( r_{j}-r_{1
			}\right) t}\frac{f_j}{f_1}\cdot\frac{\left\vert \widehat{w}^{(j)}u\right\vert
		}{\left\vert \widehat{w}^{(1)}u\right\vert}G_j(t,u),
	\end{equation*}
	where
	$$
	G_j(t,u):=\frac{g_j\left(t,u\right)}{
		g_1\left(t,u\right) }.
	$$
	Moreover, we have 
	\begin{eqnarray*}
		&&\epsilon(t)=\sum\limits_{j=2}^{q}\mathrm{e}^{\left( r_{j}-r_{1
			}\right) t}\frac{f_j}{f_1}G_j(t),
	\end{eqnarray*}
	where
	$$
	G_j(t):=\frac{g_j\left(t\right)}{g_1\left(t\right) },\ j\in\{1,\ldots,q\}.
	$$
	Here, $f_j$ and $g_j$ are the quantities introduced in Subsections \ref{unitvectors} and \ref{thenumbers}, respectively.
\end{theorem}

\begin{proof}
	We have (see Theorem \ref{two})
	\begin{equation*}
		\epsilon(t,u)=\sum\limits_{ j=2}^{q}\mathrm{e}^{\left( r_{j}-r_{1}\right)
			t}\frac{\left\Vert
			Q_{j}(t)u\right\Vert }{\left\Vert Q_{1}(t)u\right\Vert}
	\end{equation*}
	and
	\begin{equation*}
		\epsilon(t)=\sum\limits_{j=2}^{q}\mathrm{e}^{\left( r_{j}-r_{1}\right)
			t}\frac{\left\Vert
			Q_{j}(t)\right\Vert }{\left\Vert Q_{1}(t)\right\Vert }.
	\end{equation*}
	Now, use (\ref{normQJ0ureal}) and (\ref{normQJ0}).
\end{proof}

Regarding the ratios $G_j(t,u)$ and $G_j(t)$, $j\in\{2,\ldots,q\}$, in Theorem \ref{onset} we can make the following observations collected in the next theorem. 

\begin{theorem}\label{onset2}
	Suppose that, for any $j\in\{1,\ldots,q\}$, $\Lambda_j$ is simple real or simple complex. Let $j\in\{2,\ldots,q\}$. 
\begin{itemize}
	\item If both $\Lambda_j$ and $\Lambda_1$ are simple real, then 
	\begin{equation*}
		G_j(t,u)=1\text{\ \ and\ \ }G_j(t)=1.
	\end{equation*}
	\item If $\Lambda_j$ is simple complex and $\Lambda_1$ is simple real, then 
	\begin{equation*}
		G_j(t,u)=2\left\Vert \widehat{\Theta}
		_j(t,u)\right\Vert\text{\ \ and\ \ }G_j(t)=2\left\Vert \widehat{\Theta}_j(t)\right\Vert.
	\end{equation*}
	In particular, when the vector norm is a $p$-norm, we have
	\begin{equation}
		G_j(t,u)\leq 2\text{\ \ and\ \ }
		G_j(t)\leq 2.  \label{leq2}
	\end{equation}
	\item If $\Lambda_j$ is simple  real and $
	\Lambda_1$ is simple  complex, then 
	\begin{equation*}
		G_j(t,u)=\frac{1}{2\left\Vert 
			\widehat{\Theta}_1(t,u)\right\Vert}\text{\ \ and\ \ }G_j(t)=\frac{1}{2\left\Vert \widehat{\Theta}_1(t)\right\Vert}.
	\end{equation*}
	In particular, when the vector norm is the Euclidean norm, we have 
	\begin{equation}
		G_j(t,u)\leq \sqrt{\frac{1}{2(1-V_1)}}\text{\ \ and\ \ }G_j(t)\leq\sqrt{\frac{1}{a_{V_1W_1}}},  \label{jk}
	\end{equation}
	where
	\begin{equation*}
		a_{V_1W_1}=\left\{
		\begin{array}{l}
			(1-V_1)(1+W_1)\text{\ if\ }V_1\leq W_1\\
			\\
			(1+V_1)(1-W_1))\text{\ if\ }V_1\geq W_1
		\end{array}
		\right.
	\end{equation*}
	is defined in (\ref{aVjWj}).
	\item If both $\Lambda_j$ and $\Lambda_1$ are simple complex, then 
	\begin{equation*}
		G_j(t,u)=\frac{\left\Vert \widehat{
				\Theta}_j(t,u)\right\Vert}{\left\Vert \widehat{\Theta}_1(t,u)\right\Vert}
		\text{\ \ and\ \ }G_j(t)=\frac{
			\left\Vert \widehat{\Theta}_j(t)\right\Vert}{\left\Vert \widehat{\Theta}
			_1(t)\right\Vert}.
	\end{equation*}
	In particular, when the vector norm is the Euclidean norm, we have
	\begin{equation}
		G_j(t,u)\leq\sqrt{\frac{1+V_j}{1-V_1}}\text{\ \ and\ \ }G_j(t)\leq\sqrt{\frac{(1+V_j)(1+W_j)}{a_{V_1W_1}}}. \label{jk2}
	\end{equation}
\end{itemize}
\end{theorem}

\begin{proof}
	Recall Lemma \ref{Lemmareal} and Lemma \ref{Lemmacomplex}. Morover:
	\begin{itemize}
		\item (\ref{leq2}) follows by point 2 and 3 in Remark \ref{remTheta};
		\item (\ref{jk}) follows by (\ref{mintheta}) and (\ref{minTheta});
		\item (\ref{jk2}) follows by (\ref{maxtheta}), (\ref{mintheta}), (\ref{maxTheta}) and (\ref{minTheta}).
	\end{itemize}
\end{proof}

The previous Theorem \ref{onset} addresses several factors influencing the decay to zero of $\epsilon(t,u)$ and $\epsilon(t)$ as $t \rightarrow +\infty$. These factors are, for $j \in \{2, \ldots, q\}$:
\begin{itemize}
	\item the decreasing exponentials $\mathrm{e}^{(r_j-r_1)t}$: whether the decrease is fast or slow depends on how well separated the rightmost eigenvalues are from the others;
	\item the ratios $\frac{f_j}{f_1}$: non-normal matrices may have large values $f_j$, $j\in\{1,2,\ldots,n\}$, but what matters here is that large ratios $\frac{{f}_j}{f_1}$, rather than large values of $f_j$, can slow the decay to zero of $\epsilon(t,u)$ and $\epsilon(t)$ as $t \rightarrow +\infty$.
	\item the ratios $\frac{\vert \widehat{w}^{(j)}u\vert}{\vert \widehat{w}^{(1)}u\vert}$: these ratios are large, and thus can slow the decay to zero of $\epsilon(t,u)$, only if $u$ has a small normalized projection $\widehat{w}^{(1)}\widehat{u}$, $\widehat{u}=\frac{u}{\Vert u\Vert}$, onto the eigenspace associated with the rightmost eigenvalues;
	\item the ratios $G_j(t,u)$ and $G_j(t)$: for the Euclidean norm as the vector norm, the ratios $G_j(t,u)$ assume large values, and thus the decay to zero of $\epsilon(t,u)$ can be slow, only if $V_1$ is close to $1$; likewise, the ratios $G_j(t)$ assume large values, and thus the decay to zero of $\epsilon(t)$ can be slow, only if $V_1$ is close to $1$ and $W_1$ is close to $1$. 
\end{itemize}
A discussion about these factors with examples for the asymptotic condition number $K_\infty(t,y_0)$ can be found in \cite{Maset2022}.

\section{Conclusion}

In the present paper, we have analyzed in depth the asymptotic directional pointwise  condition number $K_\infty(t,y_0,\widehat{z}_0)$ (corresponding to a specific initial value $y_0$ and a specific direction of perturbation $\widehat{z}_0$ of $y_0$) and the asymptotic pointwise  condition number $K_\infty(t,y_0)$ (corresponding to a specific initial value $y_0$ and the worst perturbation of $y_0$) for a real ODE (\ref{ode}) in a generic case, namely when  the set of the rightmost eigenvalues of the matrix $A$ consists of a real eigenvalue, or a complex conjugate pair of eigenvalues, with a single Jordan mini-block of maximum order.

Here is a summary of the results obtained in the present paper

\begin{itemize}
	\item Starting from Theorems \ref{one} and \ref{two}, which were included in the results of the paper \cite{M0}, we have studied the norms $\Vert Q_j(t)u\Vert$, $u\in\mathbb{R}^n$, and $\Vert Q_j(t)\Vert$ to analyze $K_\infty(t,y_0,\widehat{z}_0)$ and $K_\infty(t,y_0)$.
	
In particular, when the rightmost eigenvalues form a complex conjugate pair and the Euclidean norm is used as the vector norm, we have derived  in Propositions \ref{oscillation} and \ref{oscillation2} expressions for $\Vert \widehat{\Theta}_j(t,u)\Vert_2$ and $\Vert \widehat{\Theta}_j(t)\Vert_2$ appearing in the formulas (\ref{normQJ0ureal}) for $\Vert Q_j(t)u\Vert_2$ and (\ref{normQJ0}) for $\Vert Q_j(t)\Vert_2$, respectively, via Lemma \ref{Lemmacomplex}.

\item The asymptotic condition numbers $K_\infty(t,y_0,\widehat{z}_0)$ and $K_\infty(t,y_0)$ are determined by the first right $\widehat{v}^{\left( 1\right)} $ and last left $\widehat{w}^{\left( 1\right) }$ generalized eigenvectors  in the longest Jordan chain corresponding to the real rightmost eigenvalue, or to the complex conjugate pair of rightmost eigenvalues, of the real matrix $A$. In case of a complex conjugate pair, the imaginary part of the eigenvalues is also involved. See Theorem \ref{Ar1} and the subsequent paragraph.

\item For a complex conjugate pair of rightmost eigenvalues, we have decomposed $K_\infty(t,y_0,\widehat{z}_0)$ and $K_\infty(t,y_0)$ as the product of  time-constant \emph{oscillation scale factors} $\mathrm{OSF}(y_0,\widehat{z}_0)$ and $\mathrm{OSF}(y_0)$, respectively, and time-periodic \emph{oscillating terms} $\mathrm{OT}(t,y_0,\widehat{z}_0)$ and $\mathrm{OT}(t,y_0)$, respectively. See Subsection \ref{sectionOSFOT}.

\item Oscillation scale factors and oscillating terms have been analyzed in depth for the Euclidean norm as vector norm. See Section \ref{Section5}. 

In particular, it is determined how the extreme values (maximum and minimum values) of the oscillating terms, as $t$ varies, depend on $y_0$ and $\widehat{z}_0$. See Subsections \ref{How1} and \ref{How2}. The conclusions in these subsections require all the extensive material presented in the appendices.

The ranges over which these extreme values vary, by varying $y_0$ and $\widehat{z}_0$, depend on the moduli $V_1$ and $W_1$ of the complex numbers $\left(\widehat{v}^{\left( 1\right) }\right) ^{T}\widehat{v}^{\left( 1\right)}$ and $\widehat{w}^{(1)}\left(\widehat{w}^{(1)}\right)^T$, respectively.

These ranges are confined to nonlarge and nonsmall values if $V_1$ is not close to $1$. See Subsections \ref{V1nc1} and \ref{V1notclose}.

On the other hand, critical cases appear when $V_1$ is close to $1$ or when the ratio $\frac{1-V_1}{1-W_1}$ is small. See Subsections \ref{critcase}, \ref{critcase2}, \ref{scc} and \ref{cnscc}.

\item The dominance of the asymptotic condition numbers $K_\infty(t,y_0,\widehat{z}_0)$ and $K_\infty(t,y_0)$ at a finite time is analyzed in Section \ref{Section6}. See Theorems \ref{onset} and \ref{onset2}, in combination with Theorem \ref{two}.

The last paragraph of Section \ref{Section6} explains which factors influence how quickly the asymptotic behavior sets in. 
\end{itemize}

Practical issues related to the asymptotic condition numbers $K_\infty(t,y_0,\widehat{z}_0)$ and $K_\infty(t,y_0)$ are addressed in \cite{Maset2022}, where one can find, in addition to several real-world examples, discussions of the onset of the asymptotic behavior, the effect of the non-normality of the matrix $A$, and the non-asymptotic phase.

\bigskip

\noindent {\bf Acknowledgements:} the research was supported by the INdAM
Research group GNCS (Gruppo Nazionale di Calcolo Scientifico).

\newpage

\appendix
\clearpage
\pagenumbering{arabic}
\setcounter{page}{1}

\section{The functions $f_{V_jW_j}$, $f_{V_jW_j}^{\max}$ and $f_{V_jW_j}^{\min}$} \label{AA}

In this appendix, we analyze the functions $f_{V_jW_j}$, $f_{V_jW_j}^{\max}$ and $f_{V_jW_j}^{\min}$ introduced in Subsection \ref{fVWsect} of the main paper. These functions play a crucial role when the Euclidean norm is used as the vector norm. To simplify the notation, we omit the index $j$ in $V_j$ and $W_j$.

\subsection{The function $f_{VW}$} \label{AAfVW}

Let $V,W\in[0,1)$ and let $f_{VW}:\mathbb{R}^2\rightarrow \mathbb{R}$ be the function given by
\begin{equation}
	f_{VW}\left(\alpha,x\right):=\frac{1+V\cos\left( x+\alpha\right)}{1-W\cos\alpha}
	,\ \left(\alpha,x\right)\in\mathbb{R}^2. \label{fVW}
\end{equation}
The first and second partial derivatives of $f_{VW}$ are given by  
\begin{eqnarray}
	&& \frac{\partial f_{VW}}{\partial\alpha}\left(\alpha,x\right)=\frac{-V\sin(x+\alpha)-W\sin \alpha+VW\sin x}{\left(1-W\cos \alpha\right)^2}\notag\\
	&& \frac{\partial f_{VW}}{\partial x}\left(\alpha,x\right)=\frac{-V\sin(x+\alpha)}{1-W\cos \alpha}\notag\\
	&&\frac{\partial^2 f_{VW}}{\partial \alpha^2}\left(\alpha,x\right)=\frac{-V\cos(x+\alpha)-W\cos \alpha}{\left(1-W\cos \alpha\right)^2},\text{\ \ \ if\ \ } \frac{\partial f_{VW}}{\partial\alpha}\left(\alpha,x\right)=0,\notag\\ 
	&&\frac{\partial^2 f_{VW}}{\partial \alpha\partial x}\left(\alpha,x\right)= \frac{\partial^2 f_{VW}}{\partial x\partial\alpha}\left(\alpha,x\right)=\frac{-V\cos(x+\alpha)+VW\cos x}{\left(1-W\cos \alpha\right)^2}\notag\\
	&&\frac{\partial^2 f_{VW}}{\partial x^2}\left(\alpha,x\right)=\frac{-V\cos(x+\alpha)}{1-W\cos \alpha}\notag\\
	&&(\alpha,x)\in\mathbb{R}^2.\notag \\
	\label{df}
\end{eqnarray}

For $x\in\mathbb{R}$, let
\begin{equation*}
	U_{VW}(x):=V\mathrm{e}^{\mathrm{i}x}+W. \label{Uvalue}
\end{equation*}
If $U_{VW}(x)=0$, i.e., $V=W=0$ or $V=W\neq 0$ and $x$ is an odd multiple of $\pi$, then 
$$
f_{VW}\left(\alpha,x\right)=1,\ \alpha\in\mathbb{R}.
$$
The next theorem states the form of the $2\pi$-periodic function $f_{VW}(\cdot, x)$, for $x \in \mathbb{R}$ such that $U_{VW}(x) \neq 0$.

\begin{theorem}	\label{LemmafvW} Let $V,W\in[0,1)$ and let $x\in\mathbb{R}$ such that $U_{VW}(x)\neq 0$. The function $f_{VW}(\cdot,x)$ oscillates monotonically between its extreme values
$$
\max\limits_{\alpha\in\mathbb{R}}f_{VW}(\alpha,x)=f_{VW}(\alpha^{\max}_{VW}(x),x)
$$
and
$$
\min\limits_{\alpha\in\mathbb{R}}f_{VW}(\alpha,x)=f_{VW}(\alpha^{\min}_{VW}(x),x),
$$
where $\alpha^{\max}_{VW}(x)$ and $\alpha^{\min}_{VW}(x)$ are given by
	\begin{equation*}
		\alpha^{\max}_{VW}(x)=\arcsin\frac{VW\sin x}{\vert U_{VW}(x)\vert}-\theta_{VW}(x)
	\end{equation*}
	and
	\begin{equation*}
		\alpha^{\min}_{VW}(x)=\pi-\arcsin\frac{VW\sin x}{\vert U_{VW}(x)\vert}-\theta_{VW}(x)
	\end{equation*}
	with $\theta_{VW}(x)$ the angle in $(-\pi,\pi]$ of the polar form of $U_{VW}(x)$.
\end{theorem}

\begin{proof}
	In this proof, to simplify the notation, we write $U(x)$ for $U_{VW}(x)$ and $\theta(x)$  for $\theta_{VW}(x)$.
	
	We determine the stationary points of $f_{VW}(\cdot,x)$.
	By (\ref{df}), we see that $\alpha\in\mathbb{R}$ is a stationary point of $f_{V_jW_j}(\cdot,x)$ if and only if
	\begin{equation*}
		V\sin(x+\alpha)+W\sin \alpha=VW\sin x, \label{sinsin}
	\end{equation*}
	which can equivalently be rewritten as
	$$
	\mathrm{Im}\left(U(x)\mathrm{e}^{\mathrm{i}\alpha}\right)=VW\sin x, 
	$$
	since
	\begin{equation}
		\mathrm{Im}\left(U(x)\mathrm{e}^{\mathrm{i}\alpha}\right)=\mathrm{Im}\left(V\mathrm{e}^{\mathrm{i}(x+\alpha)}+W\mathrm{e}^{\mathrm{i}\alpha}\right)=V\sin(x+\alpha)+W\sin \alpha. \label{star}
	\end{equation}
	By writing $U(x)=\vert U(x)\vert \mathrm{e}^{\mathrm{i}\theta(x)}$ in polar form, we obtain
	\begin{equation}
	\sin(\theta(x)+\alpha)=\frac{VW\sin x}{\vert U(x)\vert}.\label{equationfVW}
	\end{equation}
	In Appendix \ref{AD}, it is proved that
	\begin{equation*}
		\left\vert\frac{VW\sin x}{\vert U(x)\vert}\right\vert< 1.  \label{AppA}
	\end{equation*} 
	Therefore, the equation (\ref{equationfVW}) has the solutions, and then the function $f_{VW}(\cdot,x)$ has the stationary points, 
	\begin{equation}
		\alpha=\arcsin\frac{VW\sin x}{\vert U(x)\vert}-\theta(x)+2\pi k,\ k\in\mathbb{Z}, \label{maxx}
	\end{equation}
	and
	\begin{equation}
		\alpha=\pi-\arcsin\frac{VW\sin x}{\vert U(x)\vert}-\theta(x)+2\pi k,\ k\in\mathbb{Z}. \label{minn}
	\end{equation}
Since in a period of the function $f_{VW}(\cdot,x)$ there are two stationary points, they are the maximum and minimum points in such a period. Therefore the set of the maximum or minimum points of $f_{VW}(\cdot,x)$ coincide with the set of the stationary points.

By deriving (\ref{star}) with respect to $\alpha$, we obtain
\begin{eqnarray*}
V\cos(x+\alpha)+W\cos \alpha&=&\frac{d}{d\alpha}\mathrm{Im}\left(U(x)\mathrm{e}^{\mathrm{i}\alpha}\right)=\mathrm{Im}\left(\frac{d}{d\alpha}U(x)\mathrm{e}^{\mathrm{i}\alpha}\right)\\
&=&\mathrm{Im}\left(\mathrm{i}U(x)\mathrm{e}^{\mathrm{i}\alpha}\right)=\mathrm{Im}\left(\mathrm{i}\vert U(x)\vert \mathrm{e}^{\mathrm{i}(\theta(x)+\alpha)}\right)\\
&=&\vert U(x)\vert \cos (\theta(x)+\alpha).
\end{eqnarray*}
Therefore, when $\alpha$ is a stationary point of $f_{VW}(\cdot,x)$, the second derivative\\ $\frac{\partial^2 f_{VW}}{\partial\alpha^2}\left(\alpha,x\right)$ in (\ref{df}) reads 
	\begin{equation}
	\frac{\partial^2 f_{VW}}{\partial\alpha^2}\left(\alpha,x\right)=\frac{-\vert U(x)\vert\cos(\theta(x)+\alpha)}{\left(1-W\cos \alpha\right)^2}. \label{secondder}
	\end{equation}
	By looking at the sign of (\ref{secondder}), we see that (\ref{maxx}) are maximum points and (\ref{minn}) are minimum points. 
\end{proof}

The $2\pi$-periodic functions $\alpha_{VW}^{\max}$ and $\alpha_{VW}^{\min}$, which are defined in the case $V\neq 0$ or $W\neq 0$, have domain $\mathbb{R}$ if $V\neq W$, and $\mathbb{R}$ without the odd multiples of $\pi$ if $V=W$. They are studied in Appendix \ref{AE}. For $V<W$, they are real analytic on $\mathbb{R}$, and, for $V\geq W$, their restrictions to $(-\pi,\pi)$ are real analytic and such restrictions can be extended to $C^1$  functions on $[-\pi,\pi]$. For $V>W$, if the polar angle function $\theta_{VW}$ is regarded as taking values in the torus $\mathbb{R}/(2\pi\mathbb{Z})$, rather than in $(-\pi,\pi]$, then $\alpha_{VW}^{\max}$ and $\alpha_{VW}^{\min}$ are real analytic on $\mathbb{R}$ and $2\pi$-periodic only modulo $2\pi$.  

\subsection{The functions $f_{VW}^{\max}$ and $f_{VW}^{\min}$} \label{fVWmaxmin}

Let $V,W\in[0,1)$. The $2\pi$-periodic functions $f_{VW}^{\max},f_{VW}^{\min}:\mathbb{R}\rightarrow \mathbb{R}$ are given by
\begin{equation}
	f_{VW}^{\max}(x)=\max\limits_{\alpha\in\mathbb{R}} f_{VW}\left(\alpha,x\right)\text{\ \ and\ \ }f_{VW}^{\min}(x)=\min\limits_{\alpha\in\mathbb{R}} f_{VW}\left(\alpha,x\right),\ x\in\mathbb{R}.  \label{fmax}
\end{equation}
They are Lipschitz continuous functions: see Appendix \ref{ACMAXMIN}.

If $V=0$ or $W=0$, then
\begin{equation}
	f_{VW}^{\max}(x)=\frac{1+V}{1-W}\text{\ \ and\ \ } f_{VW}^{\min}(x)=\frac{1-V}{1+W},\ x\in\mathbb{R},\ \label{fVWmaxzero}
\end{equation}
are constant functions.

Thus, assume $V\neq 0$ and $W\neq 0$. Under this assumption, we can write
\begin{equation}
f_{VW}^{\max}(x)=f_{VW}(\alpha_{VW}^{\max}(x),x)\text{\ \ and\ \ }f_{VW}^{\min}(x)=f_{VW}(\alpha_{VW}^{\min}(x),x) \label{sss}
\end{equation}
for $x$ in the domain of $\alpha_{VW}^{\max}$ and $\alpha_{VW}^{\min}$, i.e., for $x$ in $\mathbb{R}$ if $V\neq W$ and for $x$ in $\mathbb{R}$ without the odd multiples of $\pi$ if $V=W$.

By recalling the last paragraph of the previous Subsection \ref{AAfVW}, we see that $f_{VW}^{\max}$ and $f_{VW}^{\min}$ are real analytic on $\mathbb{R}$ if $V\neq W$, and real analytic on $\mathbb{R}$ without the odd multiples of $\pi$ if $V=W$. 

The first and second derivatives of the functions $f_{VW}^{\max}$ and $f_{VW}^{\min}$ are (see (\ref{sss}))
\begin{eqnarray}
&&\left(f^{\max}_{VW}\right)^\prime(x)=\underset{=0}{\underbrace{\frac{\partial f_{VW}}{\partial \alpha}\left(\alpha^{\max}_{VW}(x),x\right)}}\left(\alpha^{\max}_{VW}\right)^\prime(x)+\frac{\partial f_{VW}}{\partial x}\left(\alpha^{\max}_{VW}(x),x\right)\notag\\
&&\quad\quad\quad\quad\quad  =\frac{\partial f_{VW}}{\partial x}\left(\alpha^{\max}_{VW}(x),x\right)\notag\\
&&\left(f^{\max}_{VW}\right)^{\prime\prime}(x)=\frac{\partial^2 f_{VW}}{\partial \alpha\partial x}\left(\alpha^{\max}_{VW}(x),x\right)\left(\alpha^{\max}_{VW}\right)^\prime(x)+\frac{\partial^2 f_{VW}}{\partial x^2}\left(\alpha^{\max}_{VW}(x),x\right)\notag\\
&&\left(f^{\min}_{VW}\right)^\prime(x)=\underset{=0}{\underbrace{\frac{\partial f_{VW}}{\partial \alpha}\left(\alpha^{\min}_{VW}(x),x\right)}}\left(\alpha^{\min}_{VW}\right)^\prime(x)+\frac{\partial f_{VW}}{\partial x}\left(\alpha^{\min}_{VW}(x),x\right)\notag\\
&&\quad\quad\quad\quad\quad  =\frac{\partial f_{VW}}{\partial x}\left(\alpha^{\min}_{VW}(x),x\right)\notag\\
&&\left(f^{\min}_{VW}\right)^{\prime\prime}(x)=\frac{\partial^2 f_{VW}}{\partial \alpha\partial x}\left(\alpha^{\min}_{VW}(x),x\right)\left(\alpha^{\min}_{VW}\right)^\prime(x)+\frac{\partial^2 f_{VW}}{\partial x^2}\left(\alpha^{\min}_{VW}(x),x\right)\notag \\
&&x\in\mathbb{R}\text{\ if $V\neq W$},\ x\in\mathbb{R}\setminus\{k\pi:k\in\mathbb{Z}\text{\ is odd}\}\text{\ if $V=W$}. \notag\\
\label{dfmaxmin}
\end{eqnarray}

In the case $V=W$, $f^{\max}_{VW}$ and $f^{\min}_{VW}$ are not differentiable at odd multiples of $\pi$. In fact 
\begin{eqnarray*}
&&\lim\limits_{x\rightarrow -\pi^{+}}\left(f^{\max}_{VV}\right)^{\prime}(x)=\lim\limits_{x\rightarrow -\pi^{+}}\frac{\partial f}{\partial x}\left(\alpha^{\max}_{VV}(x),x\right)=\frac{V}{\sqrt{1-V^2}}\\
	&&\lim\limits_{x\rightarrow \pi^{-}}\left(f^{\max}_{VV}\right)^{\prime}(x)=\lim\limits_{x\rightarrow \pi^{-}}\frac{\partial f}{\partial x}\left(\alpha^{\max}_{VV}(x),x\right)=\frac{-V}{\sqrt{1-V^2}}\\
	&&\lim\limits_{x\rightarrow -\pi^{+}}\left(f^{\min}_{VV}\right)^{\prime}(x)=\lim\limits_{x\rightarrow -\pi^{+}}\frac{\partial f}{\partial x}\left(\alpha^{\min}_{VV}(x),x\right)=\frac{-V}{\sqrt{1-V^2}}\\
	&&\lim\limits_{x\rightarrow \pi^{-}}\left(f^{\min}_{VV}\right)^{\prime}(x)=\lim\limits_{x\rightarrow \pi^{-}}\frac{\partial f}{\partial x}\left(\alpha^{\min}_{VV}(x),x\right)=\frac{V}{\sqrt{1-V^2}}.
\end{eqnarray*}
(see (\ref{df}) and,  in Appendix \ref{AE}, table (\ref{tablealpha2})). In the case $V=W$, $f^{\max}_{VW}$ and $f^{\min}_{VW}$ are only Lipschitz continuous functions on $\mathbb{R}$, but not $C^1$.

The forms of the functions $f_{VW}^{\max}$ and $f_{VW}^{\min}$ are stated in the next Theorem.

\begin{theorem} \label{maxminf}
	Let $V,W\in[0,1)$. 	The $2\pi$-periodic function $f_{VW}^{\max}$ is even and it oscillates monotonically between the extreme values
	\begin{equation*}
\max\limits_{x\in\mathbb{R}} f_{VW}^{\max}\left(x\right)
		=\frac{1+V}{1-W}\text{\ \ and\ \ }
		\min\limits_{x\in\mathbb{R}} f_{VW}^{\max}\left(x\right)
		=\left\{
		\begin{array}{l}
			\frac{1-V}{1-W}\text{\ if\ }V\leq W\\
			\\
			\frac{1+V}{1+W}\text{\ if\ }V\geq W.
		\end{array}
		\right.
	\end{equation*}
	The maximum value is achieved if $x$ is an even multiple of $\pi$ and the minimum value is achieved if $x$ is an odd multiple of $\pi$.
	
	Moreover, the $2\pi$-periodic function $f_{VW}^{\min}$ is even and it oscillates monotonically between the extreme values
	\begin{equation*}
		\max\limits_{x\in\mathbb{R}} f_{VW}^{\min}\left(x\right)=\left\{
		\begin{array}{l}
			\frac{1+V}{1+W}\text{\ if\ }V\leq W\\
			\\
			\frac{1-V}{1-W}\text{\ if\ }V\geq W
		\end{array}
		\right.
		\text{\ \ and\ \ }\min\limits_{x\in\mathbb{R}} f_{VW}^{\min}\left(x\right)
		=\frac{1-V}{1+W}.
	\end{equation*}
	The maximum value is achieved if $x$ is an odd multiple of $\pi$ and the minimum value is achieved if $x$ is an even multiple of $\pi$.
\end{theorem}

\begin{proof}
	We prove that $f^{\max}$ is even. For $x\in\mathbb{R}$, we have
	\begin{eqnarray*}
	f_{VW}^{\max}(-x)&=&\max\limits_{\alpha\in\mathbb{R}}\frac{1+V\cos\left(-x+\alpha\right)}{1-W\cos\alpha}=\max\limits_{\beta\in\mathbb{R}}\frac{1+V\cos\left(-x-\beta\right)}{1-W\cos(-\beta)}\\
	&=&\max\limits_{\beta\in\mathbb{R}}\frac{1+V\cos\left(x+\beta\right)}{1-W\cos\beta}=f_{VW}^{\max}(x).
\end{eqnarray*}
Similarly, we can prove that $f^{\min}_{VW}$ is even.

For the part of the proof regarding the oscillatory behavior, we suppose $V\neq 0$, $W\neq 0$ and $V\neq W$. In the case $V=0$ or $W=0$ the functions $f_{VW}^{\max}$ and $f_{VW}^{\min}$ are constants and given in (\ref{fVWmaxzero}). In the case $V=W\neq 0$, a continuity argument works: see Appendix \ref{ACMAXMIN}.

Since $V\neq W$, $f_{VW}^{\max}$ and $f_{VW}^{\min}$ are real analytic on $\mathbb{R}$ and their first derivatives are
\begin{equation*}
	\left(f_{VW}^{\max}\right)^\prime\left(x\right)=\frac{\partial f_{VW}}{\partial x}(\alpha^{\max}_{VW}(x),x)\text{\ \ and\ \ }\left(f_{VW}^{\min}\right)^\prime\left(x\right)=\frac{\partial f_{VW}}{\partial x}(\alpha^{\min}_{VW}(x),x),\ x\in\mathbb{R},
\end{equation*}
(see (\ref{dfmaxmin})).
Therefore, if $x$ is a stationary point of $f_{VW}^{\max}$ or $f^{\min}_{VW}$, then (recall (\ref{df}))
\begin{equation}
	\left\{
	\begin{array}{l}
		-V\sin(x+\alpha)-W\sin \alpha+VW\sin x=0\\
		\sin(x+\alpha)=0,
	\end{array}
	\right.  \label{firstsystem}
\end{equation}
where $\alpha=\alpha^{\max}_{VW}(x)$ or $\alpha=\alpha^{\min}_{VW}(x)$. The first equation is the equation for the stationary points $\alpha$ of $f_{VW}(\cdot,x)$: $\alpha=\alpha^{\max}_{VW}(x)$ or $\alpha=\alpha^{\min}_{VW}(x)$ is a stationary point of $f_{VW}(\cdot,x)$. The second equation is equivalent to $\left(f_{VW}^{\max}\right)^\prime(x)=0$ or $\left(f_{VW}^{\min}\right)^\prime(x)=0$. It is not difficult to show that $(\alpha,x)\in\mathbb{R}^2$ is a solution of (\ref{firstsystem}) if and only if $\alpha$ and $x$ are multiples of $\pi$. 

Therefore, between two consecutive multiples of $\pi$, the derivatives of $f_{VW}^{\max}$ and $f_{VW}^{\min}$ are never zero and then $f_{VW}^{\max}$ and $f_{VW}^{\min}$ are monotonic between two multiples of $\pi$. Since in a period of $f_{VW}^{\max}$ and $f_{VW}^{\min}$ there are only two multiples of $\pi$, it turns out that, at such two multiples of $\pi$, $f_{VW}^{\max}$ and $f_{VW}^{\min}$ take the extreme values. 

It remains to determine the extreme values of $f_{VW}^{\max}$ and $f_{VW}^{\min}$. If $x$ is an even multiple of $\pi$, we have (see table (\ref{tablesalpha}) in Appendix \ref{AE})
$$\alpha_{VW}^{\max}(x)=0\text{\ \ and\ \ }\alpha_{VW}^{\min}(x)=\pi
$$
and then
$$
f^{\max}_{VW}(x)=f_{VW}(\alpha_{VW}^{\max}(x),x)=\frac{1+V}{1-W}\text{\ \ and\ \ }f^{\min}_{VW}(x)=f_{VW}(\alpha_{VW}^{\min}(x),x)=\frac{1-V}{1+W}.
$$
If $x$ is an odd multiple of $\pi$, we have (see table (\ref{tablesalpha}) in Appendix \ref{AE})
$$\alpha_{VW}^{\max}(x)=\left\{
\begin{array}{l}
	0\text{\ if\ }V<W\\
	\\
	-\pi \text{\ if\ }V>W
\end{array}
\right.\text{\ \ and\ \ }\alpha_{VW}^{\min}(x)=\left\{
\begin{array}{l}
	\pi\text{\ if\ }V<W\\
	\\
	0 \text{\ if\ }V>W
\end{array}
\right.
$$
and then
$$
f^{\max}_{VW}(x)=f_{VW}(\alpha_{VW}^{\max}(x),x)=\left\{
\begin{array}{l}
	\frac{1-V}{1-W}\text{\ if\ }V<W\\
	\\
	\frac{1+V}{1+W} \text{\ if\ }V>W
\end{array}
\right.
$$
and
$$
f^{\min}_{VW}(x)=f_{VW}(\alpha_{VW}^{\min}(x),x)=\left\{
\begin{array}{l}
	\frac{1+V}{1+W}\text{\ if\ }V<W\\
	\\
	\frac{1-V}{1-W} \text{\ if\ }V>W.
\end{array}
\right.
$$
By comparing the values of $f^{\max}_{VW}$ at even and odd multiples of $\pi$, we see that the maximum value is achieved at even multiples and the minimum at odd multiples. By comparing the values of $f^{\min}_{VW}$ at even and odd multiples of $\pi$, we see that the maximum value is achieved at odd multiples and the minimum at even multiples.

\end{proof}

In Appendix \ref{AB}, the function $f_{VW}^{\max}$ is involved in the definition of the function $H_{VW}$. Regarding the function $f_{VW}^{\max}$, we have (see (\ref{dfmaxmin}), (\ref{df}) and, in Appendix \ref{AE}, (\ref{alphamaxprime}))
\begin{eqnarray}
	&&f_{VW}^{\max}(x)=f_{VW}(\alpha_{VW}^{\max}(x),x)=\frac{1+V\cos(x+\alpha_{VW}^{\max}(x))}{1-W\cos\alpha_{VW}^{\max}(x)}\notag \\
	&&\left(f_{VW}^{\max}\right)^\prime(x)=\frac{\partial f_{VW}}{\partial x}(\alpha_{VW}^{\max}(x),x)=\frac{-V\sin(x+\alpha^{\max}_{VW}(x))}{1-W\cos \alpha^{\max}_{VW}(x)}\notag \\
	&&\left(f^{\max}_{VW}\right)^{\prime\prime}(x)=\frac{\partial^2 f_{VW}}{\partial \alpha\partial x}\left(\alpha^{\max}_{VW}(x),x\right)\left(\alpha^{\max}_{VW}\right)^\prime(x)+\frac{\partial^2 f_{VW}}{\partial x^2}\left(\alpha^{\max}_{VW}(x),x\right)\notag\\
	&&\quad\quad\quad\quad\quad  =\frac{\left(-V\cos\left(x+\alpha^{\max}_{VW}(x)\right)+VW\cos x\right)^2}{\left(1-W\cos \alpha^{\max}_{VW}(x)\right)^2\left(V\cos\left(x+\alpha^{\max}_{VW}(x)\right)+W\cos\alpha^{\max}_{VW}(x)\right)} \notag \\
	&&\quad\quad\quad\quad\quad\quad\ \  -\frac{V\cos\left(x+\alpha^{\max}_{VW}(x)\right)}{1-W\cos \alpha^{\max}_{VW}(x)}\notag\\
	&&x\in\mathbb{R}\text{\ if $V\neq W$},\ x\in\mathbb{R}\setminus\{k\pi:k\in\mathbb{Z}\text{\ is odd}\}\text{\ if $V=W$}. \notag\\
	\label{dfmax}
\end{eqnarray}
Next table shows, for $V\neq W$, the values of $f_{VW}^{\max}$, $\left(f_{VW}^{\max}\right)^\prime$ and $\left(f^{\max}_{VW}\right)^{\prime\prime}$ at $0$ and $\pi$ (recall table (\ref{tablesalpha}) of Appendix \ref{AE}).
\begin{equation}
	\begin{tabular}{|c|c|c|c|}
		\hline
		$x$ & $f_{VW}^{\max}(x)$ & $\left(f_{VW}^{\max}\right)^\prime(x)$ & $\left(f_{VW}^{\max}\right)^{\prime\prime}(x)$ \\
		\hline 
		$0$ & $\frac{1+V}{1-W}$ & $0$ & $-\frac{VW(1+V)}{(V+W)(1-W)}$\\
		\hline
		$\pi$ & $\left\{\begin{array}{l}\frac{1-V}{1-W}\text{\ if $V<W$}\\ \frac{1+V}{1+W}\text{\ if $V>W$}\end{array}\right.$  & $0$ &  $\left\{\begin{array}{l}\frac{VW(1-V)}{(W-V)(1-W)}\text{\ if $V<W$}\\\frac{VW(1+V)}{(V-W)(1+W)}\text{\ if $V>W$}\end{array}\right.$\\
		\hline
	\end{tabular}	\label{tablefmax}
\end{equation}

\section{The functions $H_{V_1W_1}$, $H_{V_1W_1}^{\max}$ and $H_{V_1W_1}^{\min}$} \label{AB}

In this appendix, we study  the functions $H_{V_1W_1}$, $H_{V_1W_1}^{\max}$ and $H_{V_1W_1}^{\min}$, introduced in Subsection \ref{OTy0} of the main paper. These functions are fundamental in the analysis of the oscillating term $\mathrm{OT}(y_0)$, when the Euclidean norm is used as the vector norm. To simplify the notation, we omit the index $1$ in $V_1$ and $W_1$.

\subsection{ The function $H_{VW}$}
Let $V,W\in[0,1)$ and let $H_{VW}:\mathbb{R}^2\rightarrow \mathbb{R}$ be the function given by
\begin{equation}
	H_{VW}(x,\beta)=\frac{f^{\max}_{VW}(x)}{1+V\cos(x+\beta)},\ (x,\beta)\in\mathbb{R}^2,  \label{Halphax}
\end{equation}
where the $2\pi$-periodic function $f^{\max}_{V W}$ is defined and studied in Appendix \ref{AA}.

For $V=0$ or $W=0$, the function $f^{\max}_{V W}$ is constant (see (\ref{fVWmaxzero}) of Appendix \ref{AA}). For $V\neq 0$ and $W\neq 0$, $f^{\max}_{VW}$ is real analytic on $\mathbb{R}$ if $V\neq W$, and real analytic on $\mathbb{R}$ without the odd multiples of $\pi$ if $V=W$. It is given in (\ref{sss}) of Appendix \ref{AA} and their derivatives are given in (\ref{dfmax}) of Appendix \ref{AA}.

Consequently, for any $\beta\in\mathbb{R}$, the $2\pi$-periodic function $H_{VW}(\cdot,\beta)$ is real analytic on $\mathbb{R}$ if $V\neq W$, and real analytic on $\mathbb{R}$ without the odd multiples of $\pi$ if $V=W$.

The first and second partial derivatives of $H_{VW}$ are given by
\begin{eqnarray}
	&&\frac{\partial H_{VW}}{\partial x}(x,\beta)=\frac{\left(f^{\max}_{VW}\right)^\prime(x)(1+V\cos(x+\beta))+f_{VW}^{\max}(x)V\sin(x+\beta)}{(1+V\cos(x+\beta))^2}\notag \\
	&&\frac{\partial H_{VW}}{\partial \beta}(x,\beta)=\frac{f^{\max}_{VW}(x)V\sin(x+\beta)}{(1+V\cos(x+\beta))^2}\notag\\
	&&\frac{\partial^2 H_{VW}}{\partial x^2}(x,\beta)\notag \\
	&&=\frac{\left(f^{\max}_{VW}\right)^{\prime\prime}(x)(1+V\cos(x+\beta))+f_{VW}^{\max}(x)V\cos(x+\beta)}{(1+V\cos(x+\beta))^2}\text{\ if\ \ } \frac{\partial H_{VW}}{\partial x}\left(x,\beta\right)=0,\notag\\
	&&\frac{\partial^2 H_{VW}}{\partial x\partial \beta}\left(x,\beta\right)= \frac{\partial^2 H_{VW}}{\partial\beta\partial x}\left(x,\beta\right)\notag \\
	&&=\frac{-\left(f^{\max}_{VW}\right)^\prime(x)V\sin(x+\beta)+f_{VW}^{\max}(x)V\cos(x+\beta)}{(1+V\cos(x+\beta))^2}\text{\ if\ \ } \frac{\partial H_{VW}}{\partial x}\left(x,\beta\right)=0,\notag\\
	&&\frac{\partial^2 H_{VW}}{\partial \beta^2}(x,\beta)=\frac{f^{\max}_{VW}(x)V\left(\cos(x+\beta)+V(1+\sin^2(x+\beta))\right)}{(1+V\cos(x+\beta))^3}\notag \\
&&(x,\beta)\in\mathbb{R}\times \mathbb{R}\text{\ if $V\neq W$},\ (x,\beta)\in\left(\mathbb{R}\setminus\{k\pi:k\in\mathbb{Z}\text{\ is odd}\}\right)\times\mathbb{R}\text{\ if $V=W$}. \notag\\
	\label{dH}
\end{eqnarray}
 where $f^{\max}_{VW}(x)$, $\left(f^{\max}_{VW}\right)^\prime(x)$ and $\left(f^{\max}_{VW}\right)^{\prime\prime}(x)$ are given in (\ref{dfmax}).
 
\subsection{ The functions $H_{VW}^{\max}$ and $H_{VW}^{\min}$}\label{HmaxHminsect}

Let $V,W\in[0,1)$. The $2\pi$-periodic functions $H_{VW}^{\max},H_{VW}^{\min}:\mathbb{R}\rightarrow \mathbb{R}$ are given by
\begin{equation}
	H_{VW}^{\max}(\beta):=\max\limits_{x\in\mathbb{R}}H_{VW}(x,\beta)\text{\ \ and\ \ }H_{VW}^{\min}(\beta):=\min\limits_{x\in\mathbb{R}}H_{VW}(x,\beta),\ \beta\in\mathbb{R}.
	\label{HmaxHmindef}
\end{equation}
The functions $H_{VW}^{\max}$ and $H_{VW}^{\min}$ are Lipschitz continuous: see Appendix \ref{ACMAXMIN}.

The forms of the functions $H_{VW}^{\max}$ and $H_{VW}^{\min}$ are stated in the next Theorem.

\begin{theorem} \label{ThAB2}
	Let $V,W\in[0,1)$. 	The $2\pi$-periodic function $H_{VW}^{\max}$ is even and it oscillates monotonically between the extreme values
	\begin{equation*}
		\max\limits_{\beta\in\mathbb{R}} H_{VW}^{\max}\left(\beta\right)
		=\frac{1+V}{(1-V)(1-W)}
	\end{equation*}
	and
	\begin{equation*}
		\min\limits_{\beta\in\mathbb{R}} H_{VW}^{\max}\left(\beta\right)
		=\left\{
		\begin{array}{l}
			\frac{1-V^2}{(1-QV)(1-W)}\text{\ if\ }Q\leq 1\\
			\\
			\frac{1+V}{(1-V)(1+W)}\text{\ if\ }Q\geq 1,
		\end{array}
		\right.
	\end{equation*}
		where
		\begin{equation}
		Q:=\frac{V(1+W)}{2W}. \label{QQ0}
		\end{equation}
	In the case $V\neq 0$ and $W=0$, we set $Q=+\infty$. In the case $V=0$ and $W=0$, $Q$ can be chosen arbitrarily. The maximum value is achieved if $\beta$ is an odd multiple of $\pi$ and the minimum value is achieved if $\beta$ is an even multiple of $\pi$.

	Moreover, the $2\pi$-periodic function $H_{VW}^{\min}$ is even and it oscillates monotonically between the extreme values
	\begin{equation*}
		\max\limits_{\beta\in\mathbb{R}} H_{VW}^{\min}\left(\beta\right)=\frac{1}{1-W}
		\text{\ \ and\ \ }\min\limits_{\beta\in\mathbb{R}} H_{VW}^{\min}\left(\beta\right)
		=\left\{
		\begin{array}{l}
			\frac{1-V}{(1+V)(1-W)}\text{\ if\ }V\leq W\\
			\\
			\frac{1}{1+W}\text{\ if\ }V\geq W.
		\end{array}
		\right.
	\end{equation*}
	The maximum value is achieved if $\beta$ is an even multiple of $\pi$ and the minimum value is achieved if $\beta$ is an odd multiple of $\pi$.
\end{theorem}

The rest of this appendix constitutes the proof of Theorem \ref{ThAB2}.

Since $f^{\max}_{VW}$ is even, we have, for $\beta\in\mathbb{R}$, 
\begin{eqnarray*}
	H_{VW}^{\max}(-\beta)&=&\max\limits_{x\in\mathbb{R}}\frac{f^{\max}_{VW}(x)}{1+V\cos(x-\beta)}=\max\limits_{y\in\mathbb{R}}\frac{f^{\max}_{VW}(-y)}{1+V\cos(-y-\beta)}\\
	&=&\max\limits_{y\in\mathbb{R}}\frac{f^{\max}_{VW}(y)}{1+V\cos(y+\beta)}=H_{VW}^{\max}(\beta).
\end{eqnarray*}
This  proves that $H^{\max}_{VW}$ is even. Similarly, we can prove that $H^{\min}_{VW}$ is even. By the $2\pi$ periodicity and the fact that $H_{VW}^{\max}$ and $H_{VW}^{\min}$ are even, we have 
\begin{equation}
	H^{\max}_{VW}(\pi+\beta)=H^{\max}_{VW}(\pi-\beta)\text{\ and\ }H^{\min}_{VW}(\pi+\beta)=H^{\min}_{VW}(\pi-\beta),\ \beta\in\mathbb{R}. \label{propHmaxmin}
\end{equation}
This shows that the analysis of the $2\pi$-periodic functions $H_{VW}^{\max}$ and $H_{VW}^{\min}$ can be confined to the interval $[0,\pi]$.

In the case $V=0$ or $W=0$, the function $f_{VW}^{\max}$ is constant: see (\ref{fVWmaxzero}) of Appendix \ref{AA}. Therefore, we have
\begin{equation*}
	H_{VW}(x,\beta)=\frac{1+V}{(1+V\cos(x+\beta))(1-W)},\ (x,\beta)\in\mathbb{R}^2,
\end{equation*}
and then
\begin{equation*}
	H_{VW}^{\max}(\beta)=\frac{1+V}{(1-V)(1-W)}\text{\ \ and\ \ }H_{VW}^{\min}(\beta)=\frac{1}{1-W},\ \beta\in\mathbb{R},
\end{equation*}
are constant functions. Theorem \ref{ThAB2} is immediate in this case.

From now on, we assume $V\neq 0$ and $W\neq 0$. Moreover, to avoid the annoying issue that, for $V=W$, the function $f^{\max}_{VW}$ is not differentiable at odd multiples of $\pi$ and, consequently, the same holds for the functions $H_{VW}(\cdot,\beta)$, $\beta\in\mathbb{R}$, we also assume $V\neq W$.

Since $V\neq W$, the functions $H_{VW}(\cdot,\beta)$, $\beta\in\mathbb{R}$, are real analytic on $\mathbb{R}$. As noted above (see (\ref{propHmaxmin})), we may restrict our attention to $\beta\in[0,\pi]$.

For a given $\beta\in[0,\pi]$, the values $H^{\max}_{VW}(\beta)$ and $H^{\min}_{VW}(\beta)$ are values $H_{VW}(x,\beta)$, where $x$ is stationary point of the function $H_{VW}(\cdot,\beta)$. The stationary points of $H_{VW}(\cdot,\beta)$ are the points $x\in\mathbb{R}$ such that
\begin{equation}
	\frac{\partial H_{VW}}{\partial x}(x,\beta)=0. \label{dHdx=0}
\end{equation}
By using the expression of this partial derivative given in (\ref{dH}) and the expressions of $f^{\max}_{VW}(x)$ and $\left(f^{\max}_{VW}\right)^\prime(x)$ given in (\ref{dfmax}) of Appendix \ref{AA}, we see that (\ref{dHdx=0}) is equivalent to 
\begin{equation}
	-\sin(x+\alpha^{\max}_{VW}(x))+\sin(x+\beta)-V\sin(\alpha^{\max}_{VW}(x)-\beta)=0. \label{eq}
\end{equation}

In order to simplify the proof of Proposition \ref{monbr0} below, the next assumption is introduced.
\begin{assumption}\label{Acondneqzero}	(NDR) For any stationary point $x_0\in\mathbb{R}$ of  $H_{VW}(\cdot,0)$, we have 
\begin{equation}
	\frac{\partial^2 H_{VW}}{\partial x^2}(x_0,0)\neq 0, \label{condneqzero}
\end{equation}
i.e., any root of (\ref{dHdx=0}) with $\beta=0$ is not a double root.
\end{assumption}

Therefore, in the rest of this appendix, we assume $V\neq 0$, $W\neq 0$, $V\neq W$ and NDR. If $V=0$ or $W=0$, Theorem \ref{ThAB2} is immediate (as we have seen above). If $V\neq 0$, $W\neq 0$ and $V=W$ or NDR does not hold, we use a continuity argument: see Appendix \ref{ACMAXMIN}.

\subsection{The branches $x_i$ and $h_i$}\label{Brxihi}
For any fixed $\beta\in [0,\pi]$, the equation (\ref{dHdx=0}) has finitely many solutions $x\in (-\pi,\pi]$. This follows by the fact that the function $\frac{\partial H_{VW}}{\partial x}(\cdot,\beta)$ is real analytic on $\mathbb{R}$ (it is the derivative of $H(\cdot,\beta)$ real analytic on $\mathbb{R}$) and it is not identically zero (consider, in (\ref{dH}), $\frac{\partial H_{VW}}{\partial x}(0,\beta)$ for $\beta\neq 0$ and $\beta\neq \pi$, and $\frac{\partial^2 H_{VW}}{\partial x^2}(0,\beta)$ for $\beta=0$ or $\beta=\pi$). Viewing these solutions as functions of $\beta$, we obtain a family of branches, as explained below.

For a given $\beta_0\in[0,\pi]$, let $x_0$ be a solution of (\ref{dHdx=0}) with $\beta=\beta_0$. If
$$
\frac{\partial^2 H_{VW}}{\partial x^2}(x_0,\beta_0)\neq 0,
$$
then, by the implicit function theorem, there exist a function $x(\beta)$ of $\beta$, $\beta$ in a neighborhhod of $\beta_0$, such that $x(\beta_0)=x_0$ and $x(\beta)$ is a solution of   
(\ref{dHdx=0}). Such function is the solution of the Cauchy problem
\begin{equation*}
\left\{
\begin{array}{l}
\frac{d}{d\beta}\frac{\partial H_{VW}}{\partial x}(x(\beta),\beta)=\frac{\partial^2 H_{VW}}{\partial x^2}(x(\beta),\beta)x^\prime(\beta)+\frac{\partial^2 H_{VW}}{\partial \beta\partial x}(x(\beta),\beta)=0\\
x(\beta_0)=x_0
\end{array}
\right.
\end{equation*}
The solution $x$, defined on an open interval $D$ of $[0,\pi]$, admits a continuous extension to the endpoints of $D$. The existence of the limit of $x(\beta)$, as $\beta$ approaches an endpoint of $D_i$ follows from the monotonicity of $x$ in a neighborhood of that endpoint. This monotonicity is a consequence of the fact that
\begin{equation}
	\frac{\partial^2 H_{VW}}{\partial \beta \partial x}(x(\beta),\beta)=0 \label{d2Hdx2}
\end{equation}
has only finitely many solutions in $[0,\pi]$: see Appendix \ref{AF}; and then there are only finitely many zeros of $x^\prime$. Therefore, for $\overline{\beta}$ endpoint of $D$, we have, with $\beta$ interior point of $D$, 
$$
\frac{\partial H_{VW}}{\partial x}(x\left(\overline{\beta}\right),\overline{\beta})=\lim\limits_{\beta\rightarrow \overline{\beta}}\frac{\partial H_{VW}}{\partial x}(x\left(\beta\right),\beta)=0,
$$
since $\frac{\partial H_{VW}}{\partial x}$ is continuous. This shows that $x\left(\overline{\beta}\right)$ is a solution of (\ref{dHdx=0}). The endpoints of the interval $D$ are $0$, $\pi$ or points $\beta\in(0,\pi)$ such that
$$
\frac{\partial^2 H_{VW}}{\partial x^2}(x(\beta),\beta)=0.
$$

Therefore, the stationary points of $H_{VW}(\cdot,\beta)$, i.e., the solutions of (\ref{dHdx=0}), form, as a function of $\beta\in[0,\pi]$, a family of branches
$$
x_i:D_i\rightarrow  \mathbb{R}/(2\pi\mathbb{Z}),\ \beta\mapsto x_i(\beta),\ i\in I,
$$
where $I$ is a finite set of indices and $D_i$ is a closed interval included in $[0,\pi]$. The  set $I$ is a finite, since there are finitely many $(x,\beta)\in(-\pi,\pi]\times [0,\pi]$ such that
$$
\frac{\partial H_{VW}}{\partial x}(x,\beta)=0\text{\ \ and\ \ }\frac{\partial^2 H_{VW}}{\partial x^2}(x,\beta)=0
$$
(and then there are finitely many endpoints of intervals $D_i$): see  Appendix \ref{A2H}. On the interior of $D_i$, the branch $x_i$ satisfies the differential equation
\begin{equation}
	\frac{\partial^2 H_{VW}}{\partial x^2}(x_i(\beta),\beta)x_i^\prime(\beta)+\frac{\partial^2 H_{VW}}{\partial \beta\partial x}(x_i(\beta),\beta)=0,
	\label{odexi}
\end{equation}
where
$$
\frac{\partial^2 H_{VW}}{\partial x^2}(x_i(\beta),\beta)\neq 0
$$
holds.

The values of $H_{VW}(\cdot,\beta)$ at the stationary points also form a family of branches
$$
h_i:D_i\rightarrow \mathbb{R},\ \beta\mapsto H_{VW}(x_i(\beta),\beta),\ i\in I.
$$

\begin{figure}
	\centering
	\subfloat[Branches $x_i$: $V=0.45$, $W=0.5$.]{\includegraphics[width=0.45\textwidth]{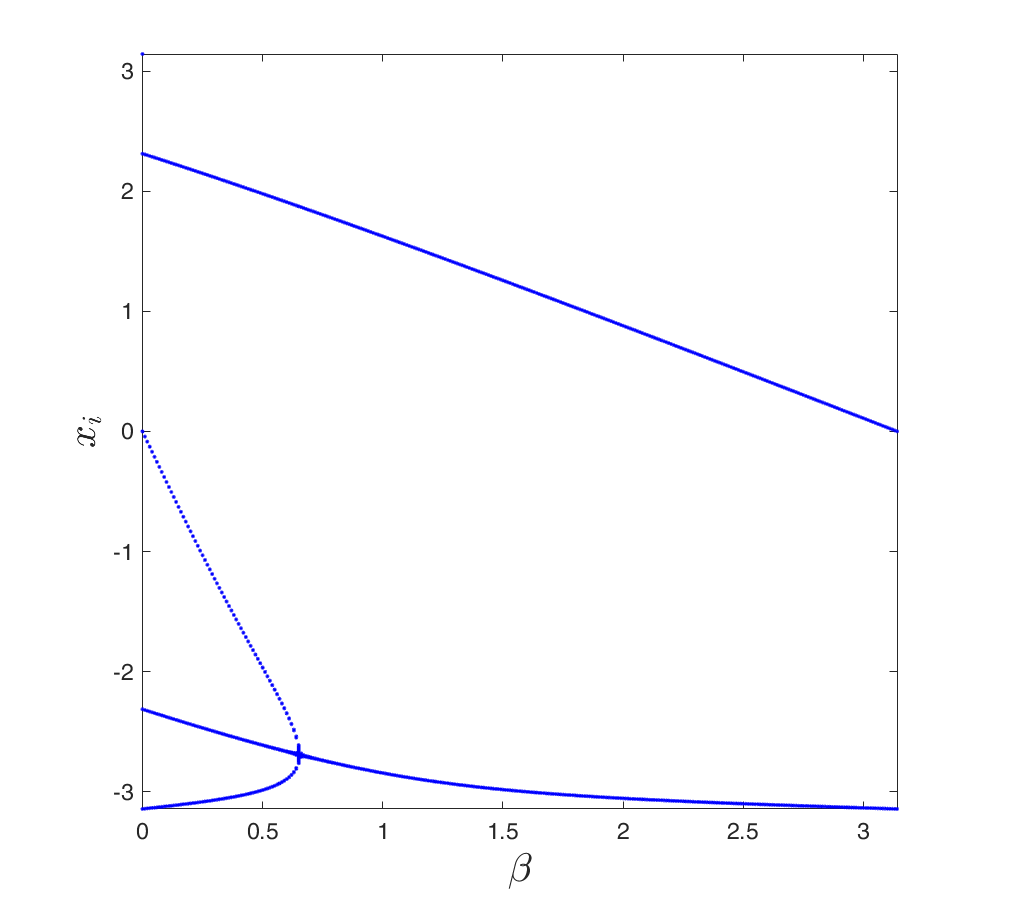}}
	\hfill
	\subfloat[Branches $h_i$: $V=0.45$, $W=0.5$.]{\includegraphics[width=0.45\textwidth]{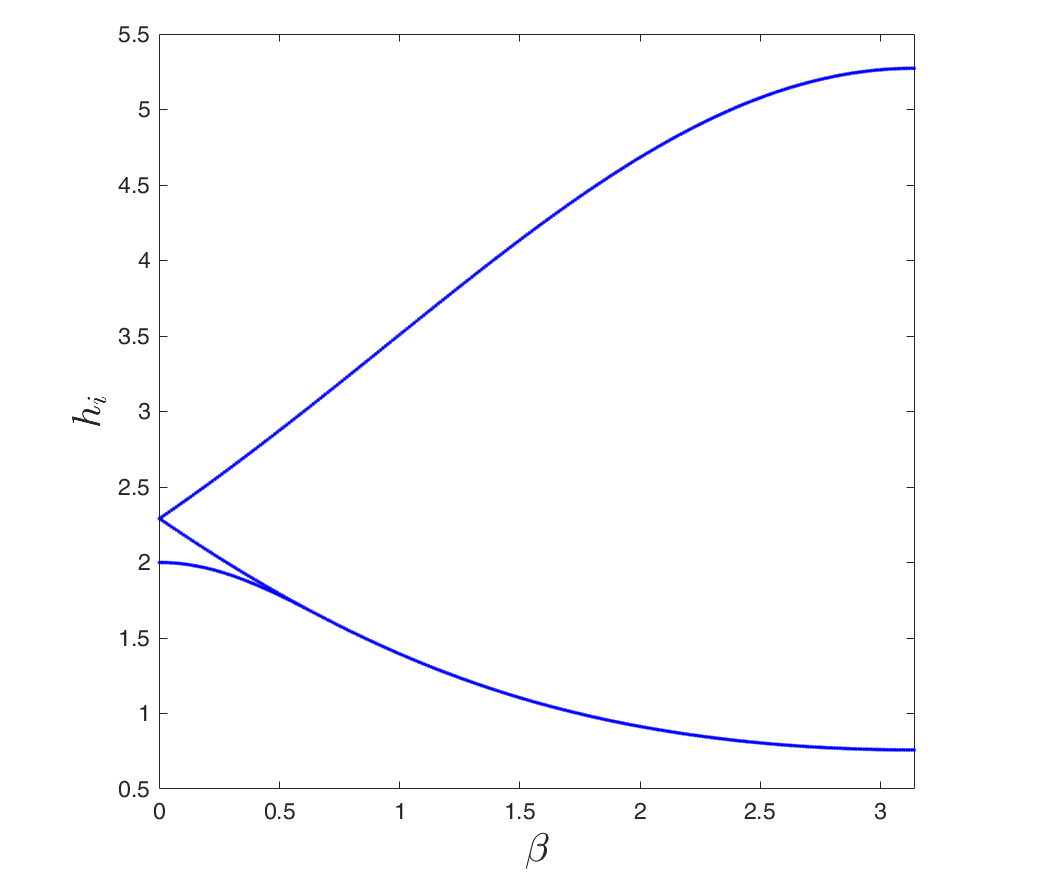}}
	
	\vspace{0.5cm}

	\subfloat[Branches $x_i$: $V=0.55$, $W=0.5$.]{\includegraphics[width=0.45\textwidth]{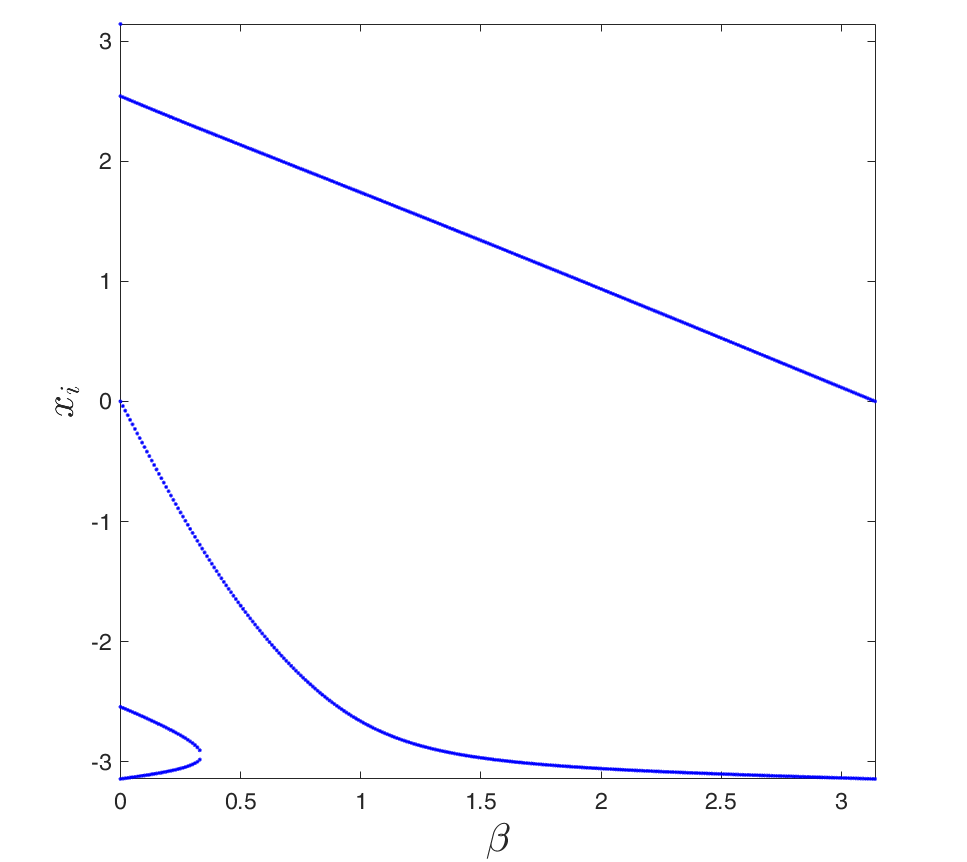}}
	\hfill
	\subfloat[Branches $h_i$: $V=0.55$, $W=0.5$.]{\includegraphics[width=0.45\textwidth]{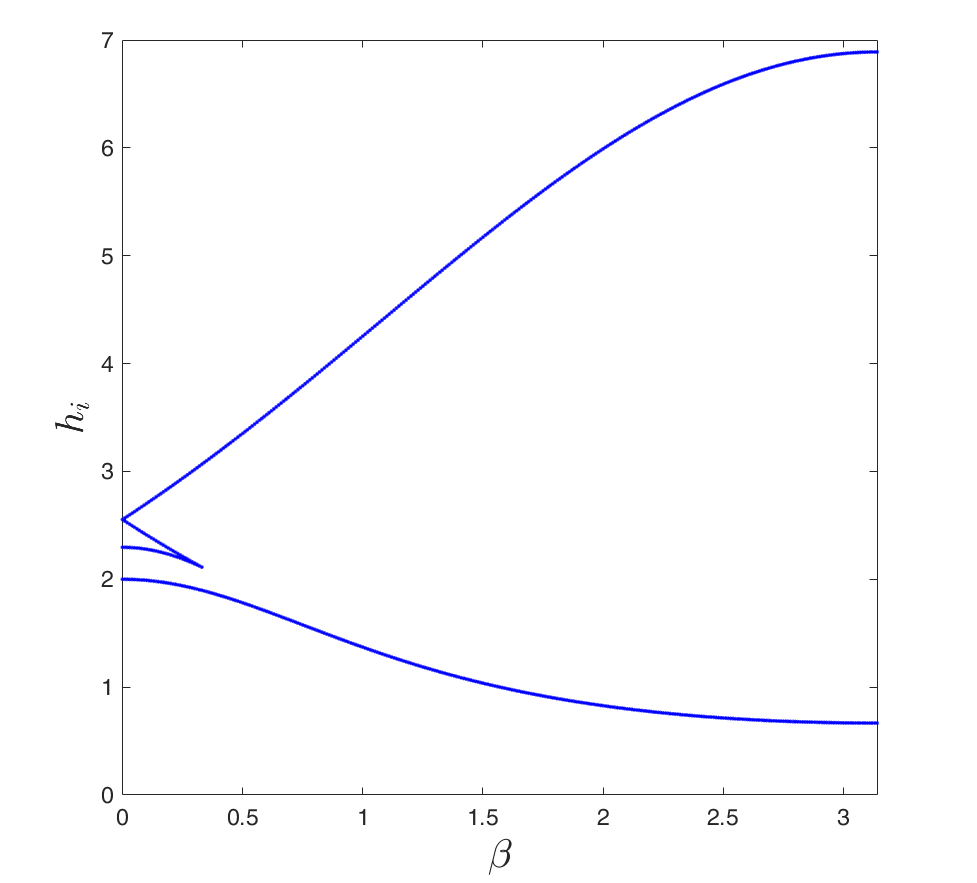}}
	
	\vspace{0.5cm} 
	
	\subfloat[Branches $x_i$: $V=0.7$, $W=0.5$.]{\includegraphics[width=0.45\textwidth]{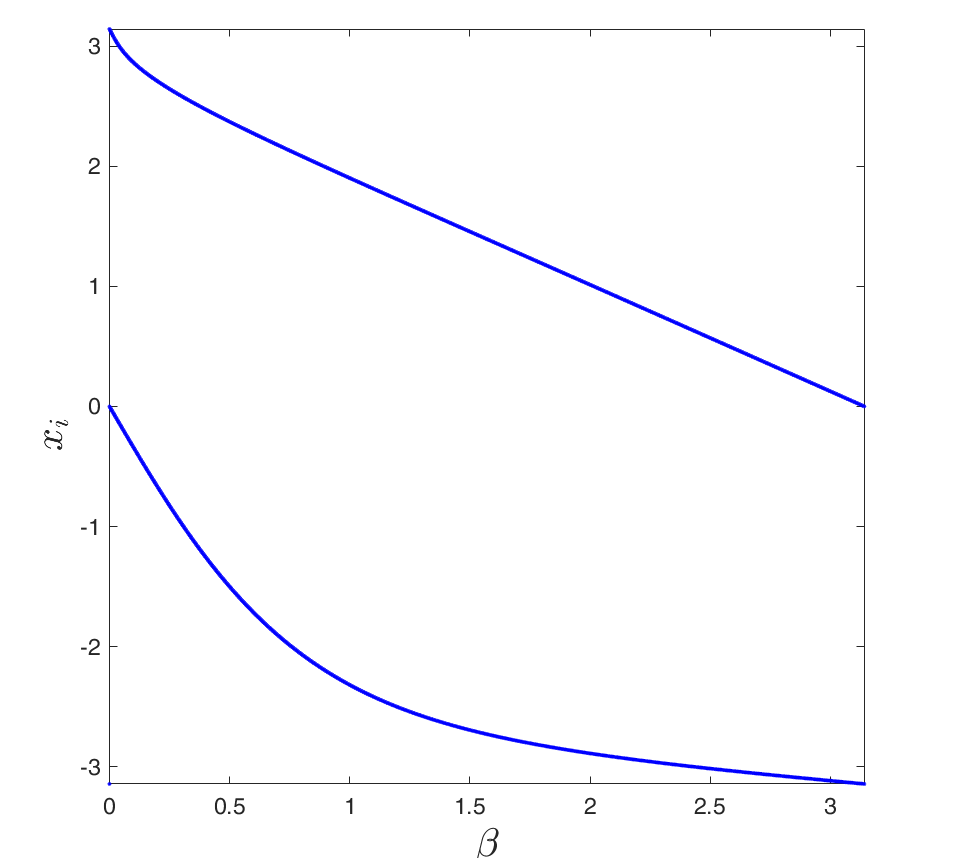}}
	\hfill
	\subfloat[Branches $h_i$: $V=0.7$, $W=0.5$.]{\includegraphics[width=0.45\textwidth]{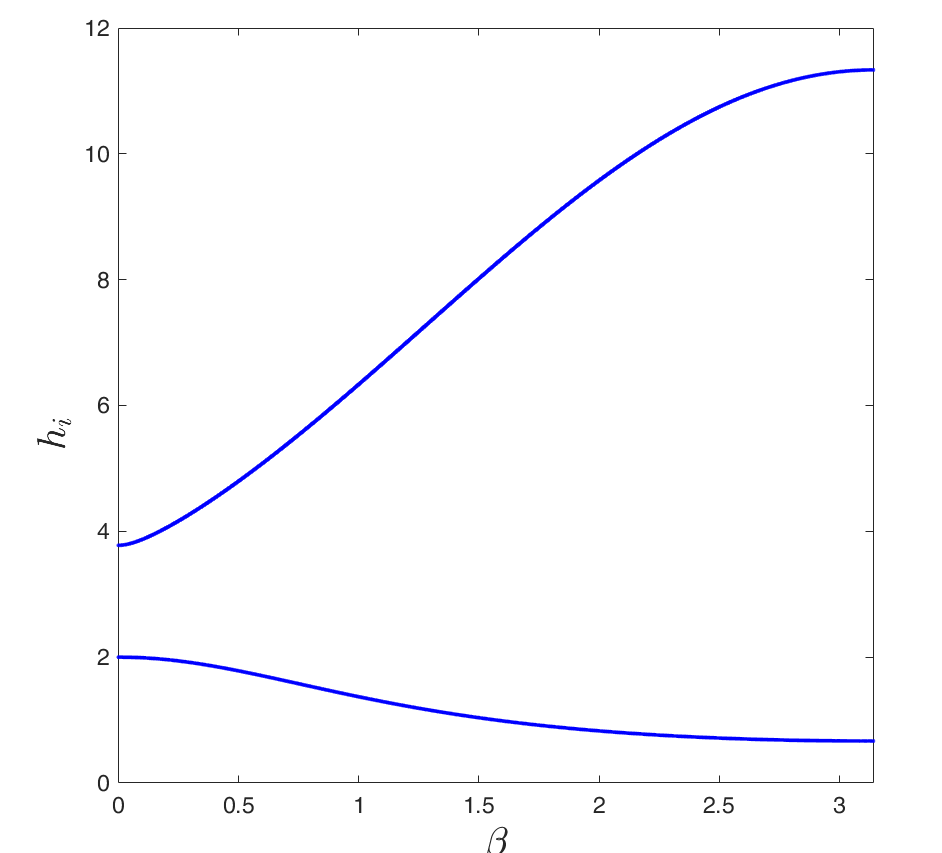}}
	
	\caption{Left: branches $x_i$ for $(\beta,x)\in[0,\pi]\times[-\pi,\pi]$. Right: branches $h_i$ for $\beta\in[0,\pi]$.}
	
	\label{branches}
\end{figure}

In Figure \ref{branches}, we see examples of such branches $x_i$ and $h_i$.

\begin{remark} \label{twosets}
	Observe that the equation (\ref{eq}), equivalent to (\ref{dHdx=0}), can be factorized as
	$$
	-2\sin\frac{\alpha^{\max}_{VW}(x)-\beta}{2}\left(\cos\left(x+\frac{\alpha^{\max}_{VW}(x)+\beta}{2}\right)+V\cos\frac{\alpha^{\max}_{VW}(x)-\beta}{2}\right)=0
	$$
	and then its solutions are the union of the solutions of
	\begin{equation}
		\alpha^{\max}_{VW}(x)-\beta\text{\ is a multiple of $2\pi$}\label{2pi}
	\end{equation}
	and
	\begin{equation}
	\cos\left(x+\frac{\alpha^{\max}_{VW}(x)+\beta}{2}\right)+V\cos\frac{\alpha^{\max}_{VW}(x)-\beta}{2}=0. \label{2piother}
	\end{equation}
	The branches $h_i$ corresponding to solutions of (\ref{2pi}) are given by
	\begin{equation}
	h_i(\beta)=\frac{1}{1-W\cos \beta} \label{hiWcos}
	\end{equation}
	as can be immediately verified by using (\ref{Halphax}) and (\ref{dfmax}) in Appendix \ref{AA}. In Figure \ref{branchespartial}, we see examples of such branches $x_i$ and $h_i$ corresponding to solutions of (\ref{2pi}). Compare with Figure \ref{branches}.
\end{remark}
\begin{figure}
	\centering
	\subfloat[Branches $x_i$: $V=0.45$, $W=0.5$.]{\includegraphics[width=0.45\textwidth]{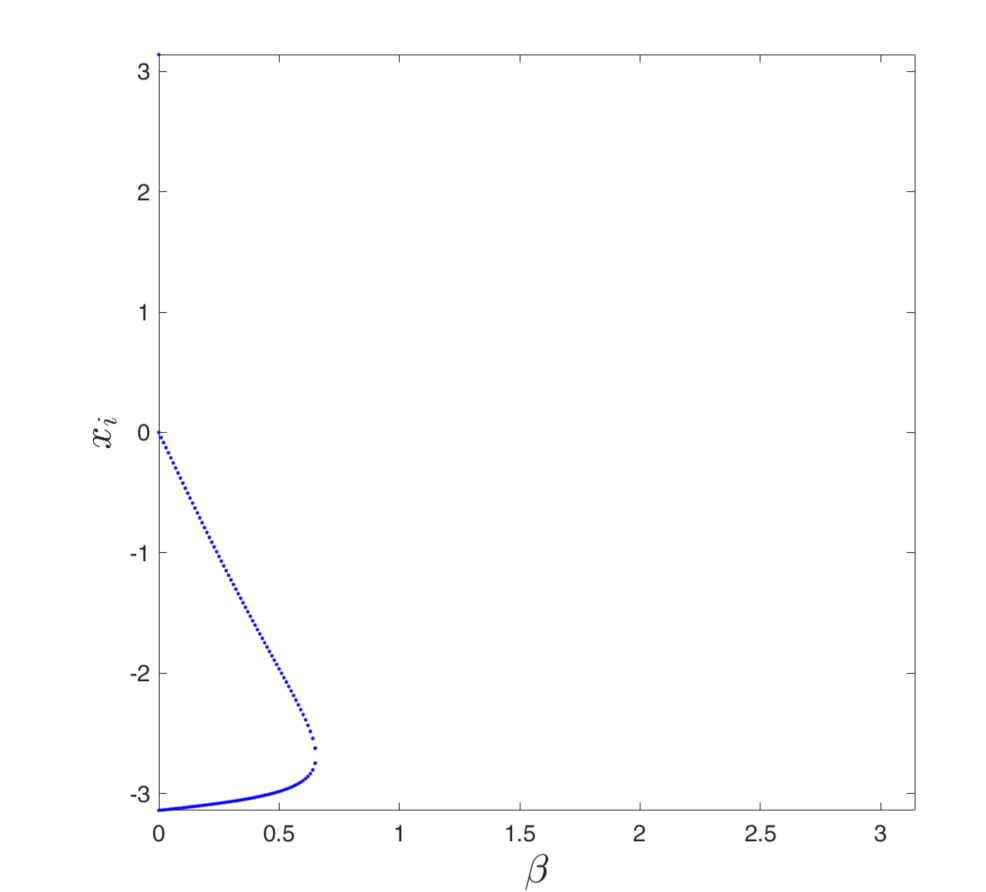}}
	\hfill
	\subfloat[Branches $h_i$: $V=0.45$, $W=0.5$.]{\includegraphics[width=0.47\textwidth]{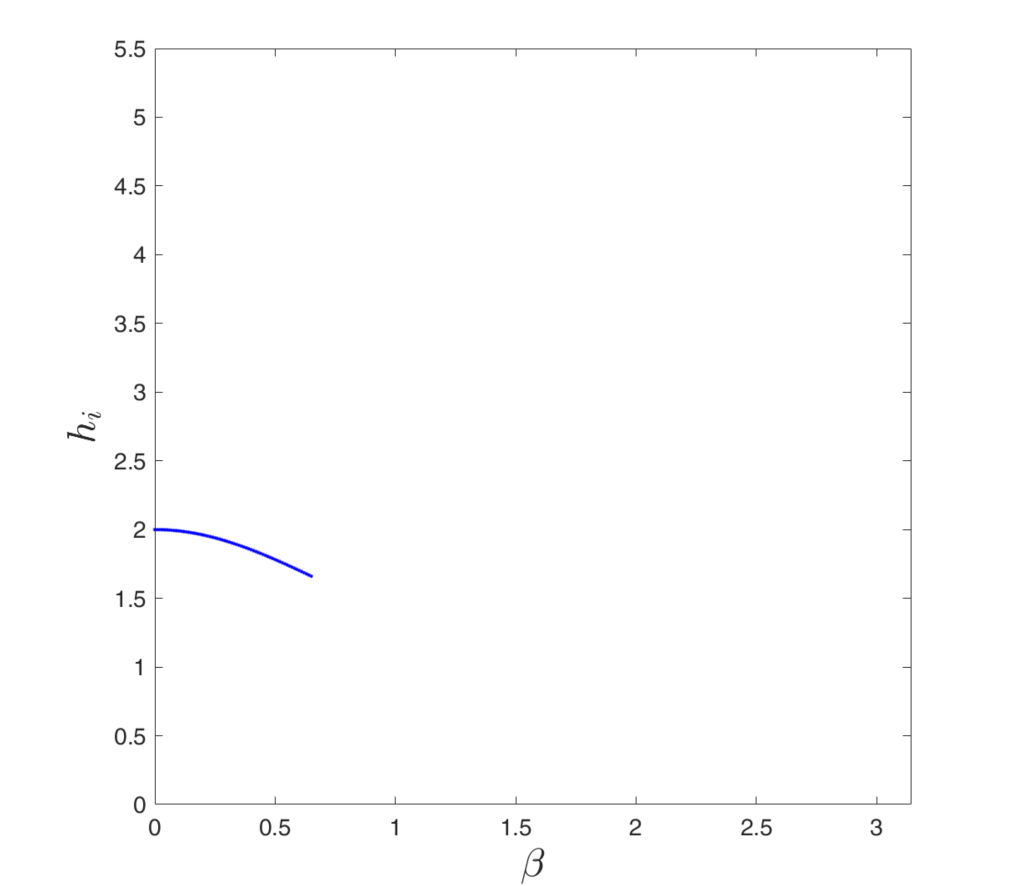}}
	
	\vspace{0.5cm} 
	
	\subfloat[Branches $x_i$: $V=0.55$, $W=0.5$.]{\includegraphics[width=0.45\textwidth]{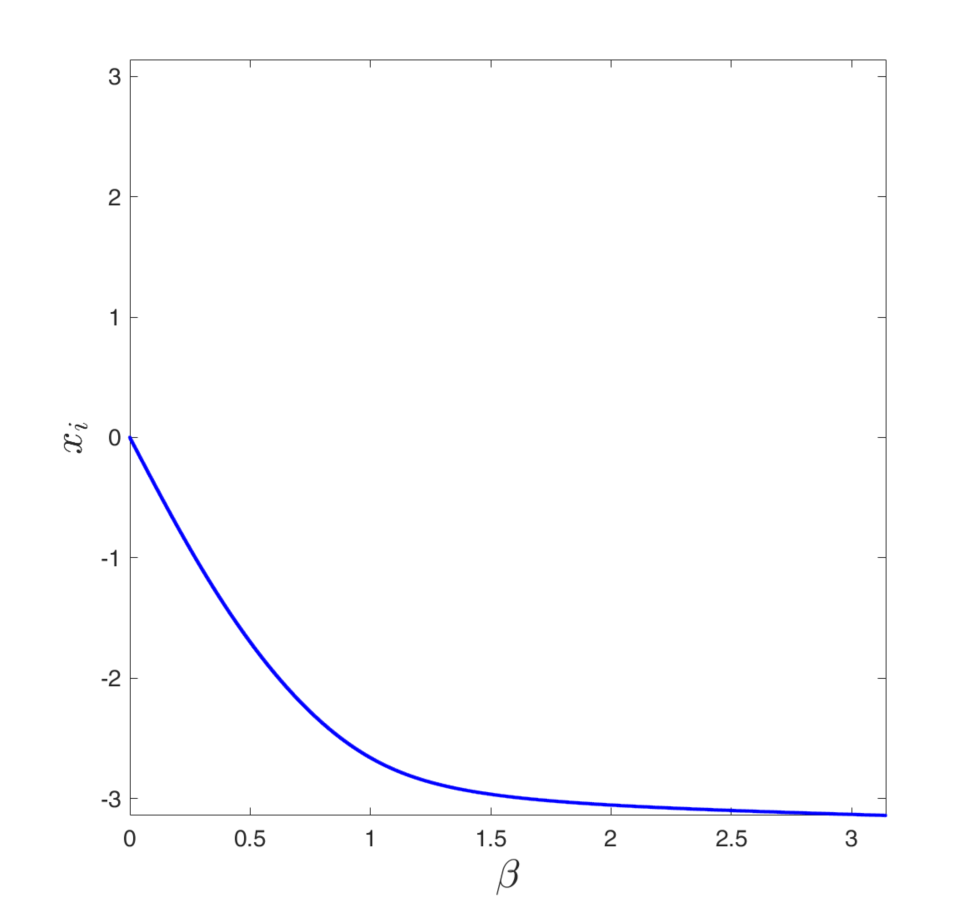}}
	\hfill
	\subfloat[Branches $h_i$: $V=0.55$, $W=0.5$.]{\includegraphics[width=0.48\textwidth]{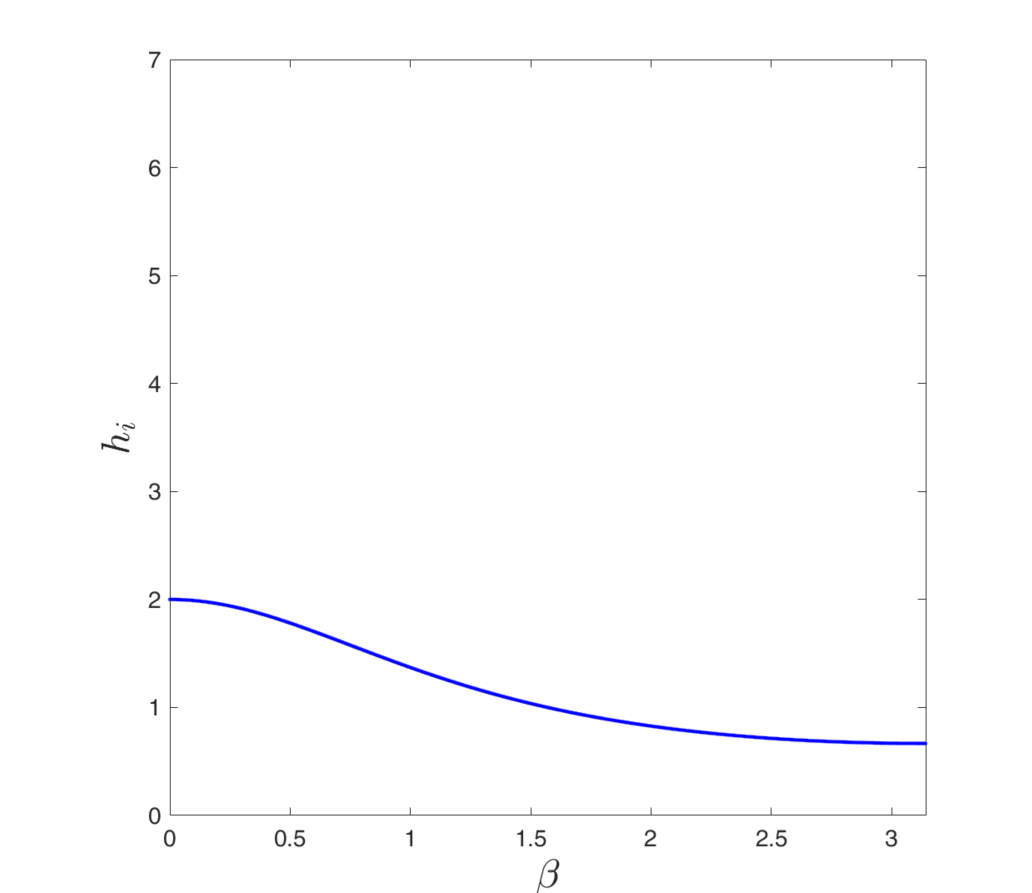}}
	
	\vspace{0.5cm} 
	
	\subfloat[Branches $x_i$: $V=0.7$, $W=0.5$.]{\includegraphics[width=0.48\textwidth]{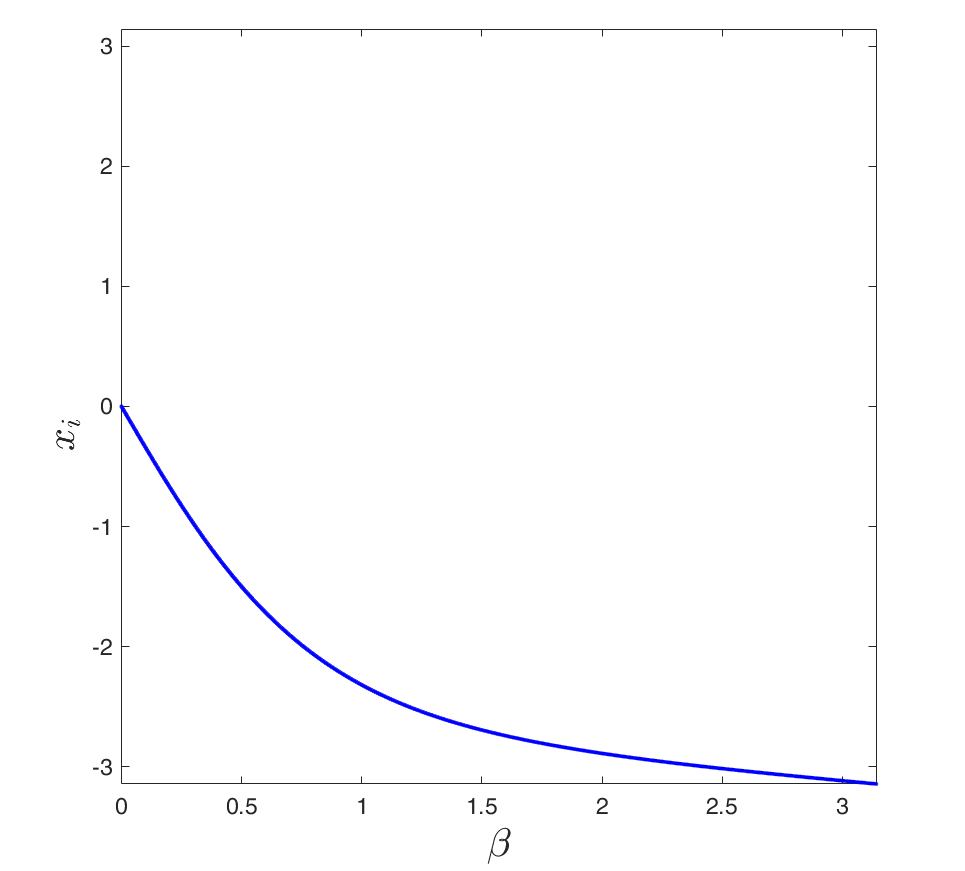}}
	\hfill
	\subfloat[Branches $h_i$: $V=0.7$, $W=0.5$.]{\includegraphics[width=0.48\textwidth]{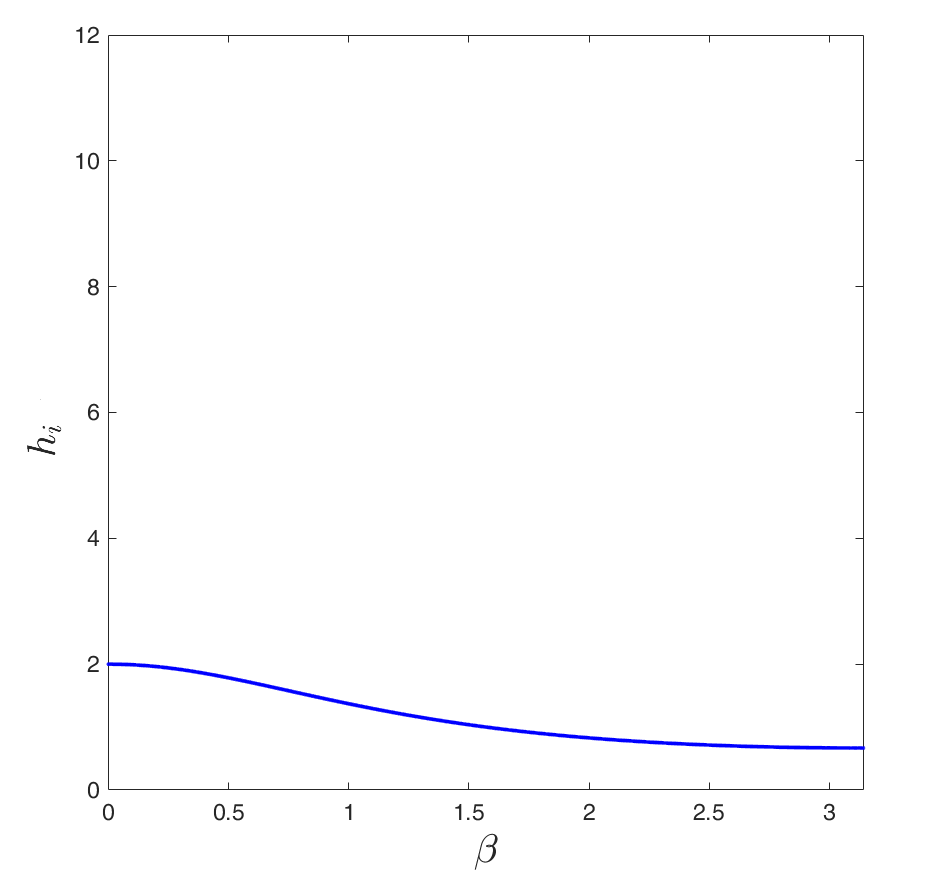}}

	\caption{Left: branches $x_i$ for $(\beta,x)\in[0,\pi]\times[-\pi,\pi]$ corresponding to solutions of (\ref{2pi}). Right: branches $h_i$ for $\beta\in[0,\pi]$ corresponding to solutions of (\ref{2pi}).}
	
	\label{branchespartial}
\end{figure}

The next three propositions regard the branches $h_i$.

\begin{proposition} \label{diff}
	Any branch $h_i$, $i\in I$, is continuously differentiable with derivative
	\begin{equation}
		h_i^\prime(\beta)=\frac{\partial H_{VW}}{\partial \beta}(x_i(\beta),\beta)=\frac{f^{\max}_{VW}(x_i(\beta))V\sin(x_i(\beta)+\beta)}{(1+V\cos(x_i(\beta)+\beta))^2},\ \beta\in D_i.	\label{dbeta}
	\end{equation}
\end{proposition}
\begin{proof}
	For $\beta$ interior point of $D_i$, $x_i$ is differentiable in $\beta$ and then 
\begin{eqnarray*}
	h_i^\prime(\beta)&=&\frac{d}{d\beta}H_{VW}(x_i(\beta),\beta)=\underset{=0}{\underbrace{\frac{\partial H_{VW}}{\partial x}(x_i(\beta),\beta)}}x_i^{\prime}(\beta)+\frac{\partial H_{VW}}{\partial \beta}(x_i(\beta),\beta)\\
	&=&\frac{\partial H_{VW}}{\partial \beta}(x_i(\beta),\beta),
\end{eqnarray*}
where the partial derivative $\frac{\partial H_{VW}}{\partial \beta}$ is given in (\ref{dH}). For $\overline{\beta}\in D_i$ endpoint of $D_i$, we have, with $\beta$ interior point of $D_i$, 
$$
\lim\limits_{\beta\rightarrow \overline{\beta}}h_i^\prime(\beta)=\lim\limits_{\beta\rightarrow \overline{\beta}}\frac{\partial H_{VW}}{\partial \beta}(x_i(\beta),\beta)=\frac{\partial H_{VW}}{\partial \beta}(x_i\left(\overline{\beta}\right),\overline{\beta}),
$$
 by the continuity of $\frac{\partial H_{VW}}{\partial \beta}$ and $x_i$. We conclude that $h_i$ is differentiable at the endpoint $\overline{\beta}$ and
$$
h_i^\prime\left(\overline{\beta}\right)=\frac{\partial H_{VW}}{\partial \beta}(x_i\left(\overline{\beta}\right),\overline{\beta}).
$$

\end{proof}

\begin{proposition} \label{diffneq0}
	For any branch $h_i$, $i\in I$, we have $h_i^\prime(\beta)=0$ only if $\beta=0$ or $\beta=\pi$.
\end{proposition}

\begin{proof}
If $\beta\in D_i$ is such that $h^\prime_i(\beta)=0$, then
	\begin{equation}
		\left\{
		\begin{array}{l}
		-V\sin(x+\alpha)-W\sin \alpha+VW\sin x=0\\
		-\sin(x+\alpha)+\sin(x+\beta)-V\sin(\alpha-\beta)=0\\
		\sin(x+\beta)=0,
		\end{array}
		\right.
		\label{systemthree}
	\end{equation}
	where $x=x_i(\beta)$ and $\alpha=\alpha^{\max}_{VW}(x_i(\beta))$. Regarding the equations in (\ref{systemthree}), observe that:
	\begin{itemize}
		\item the first equation is the equation for the stationary points $\alpha$ of $f_{VW}(\cdot,x_i(\beta))$ (recall (\ref{df}) in Appendix \ref{AA}): $\alpha=\alpha^{\max}_{VW}(x_i(\beta))$ is a stationary point of $f_{VW}(\cdot,x_i(\beta))$;
		\item the second equation is the equation for the stationary points of $H_{VW}(\cdot,\beta)$ (recall (\ref{eq})): $x=x_i(\beta)$ is a stationary point of $H_{VW}(\cdot,\beta)$;
		\item the third equation is equivalent to the equation $h_i^{\prime}(\beta)=0$: recall (\ref{dbeta}) in Proposition  \ref{diff} and observe that $f^{\max}_{VW}(x_i(\beta))>0$ (see, in Appendix \ref{AA}, Theorem \ref{maxminf}, or the definition (\ref{fVW}) of $f_{VW}$).
	\end{itemize} 
	
	Then, $h^\prime_i(\beta)=0$ is equivalent to
	
	It is not difficult to show that $(\alpha,x,\beta)\in\mathbb{R}^3$ is a solution of (\ref{systemthree}) if and only if $\alpha$, $x$ and $\beta$ are multiples of $\pi$. The theorem follows.
\end{proof}

\begin{proposition}\label{miormc}
	Any branch $h_i$, $i\in I$, is either  monotonically increasing or monotonically decreasing.
\end{proposition}
\begin{proof}
The proposition is a consequence of the previous Proposition \ref{diffneq0}. Since the interior points of $D_i$ cannot be $0$ or $\pi$, there are no interior points of $D_i$ where the continuous derivative $h_i^\prime$ is zero. 
\end{proof}
\subsection{Monotonicity of $H_{VW}^{\max}$ and $H_{VW}^{\min}$}
In this subsection, we show that the functions $H_{VW}^{\max}$ and $H_{VW}^{\min}$ are monotonically increasing and decreasing, respectively, in $[0,\pi]$.

Since, for $\beta\in[0,\pi]$, we have $H_{VW}^{\max}(\beta)=H_{VW}\left(x,\beta\right)$ for some stationary point $x$ of $H_{VW}(\cdot,\beta)$, and the same is true for $H_{VW}^{\min}(\beta)$, it follows that $H^{\max}_{VW}(\beta)$ and $H_{VW}^{\min}(\beta)$ are values $h_i(\beta)=H_{VW}(x_i(\beta),\beta)$ for some $i\in I$.

\begin{figure}[tbp]
	\includegraphics[width=0.8\textwidth]{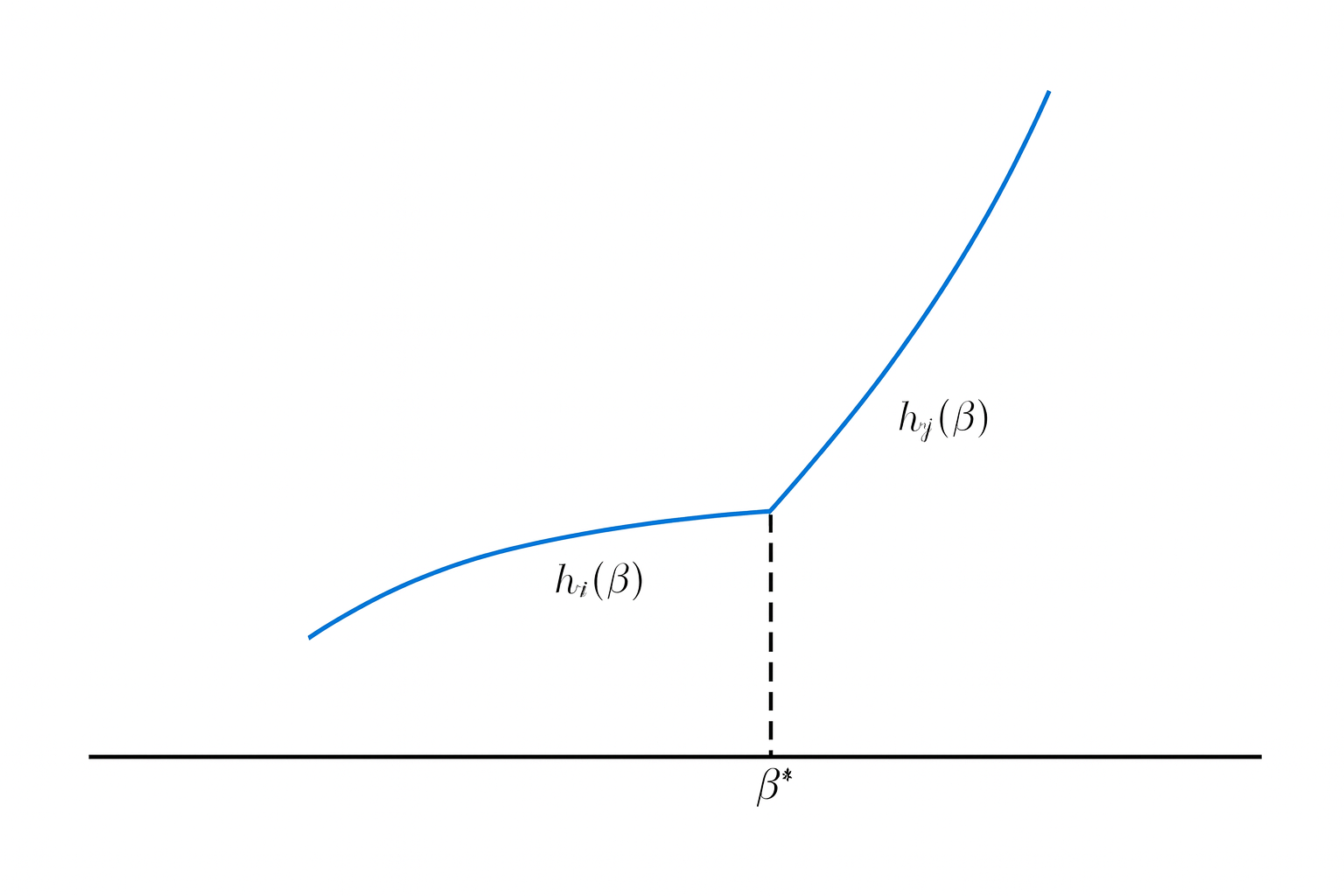}
	\caption{Junction in $\beta^\ast$ of two consecutive branches $h_i$ and $h_j$.}
	\label{Hmin} 
\end{figure}

Therefore, we obtain
\begin{equation} H^{\max}_{VW}(\beta)=\max\limits_{\substack{i\in I\text{ such that}\\ \beta\in D_i}}h_i(\beta)\text{\ \ and\ \ }H^{\min}_{VW}(\beta)=\min\limits_{\substack{i\in I\text{\ such that}\\ \beta\in D_i}}h_i(\beta),\ \beta\in[0,\pi]. \label{Hmaxmin}
\end{equation}
By (\ref{Hmaxmin}), we see that the continuous functions $H_{VW}^{\max}$ and $H_{VW}^{\min}$, considered in $[0,\pi]$, are the union of branches $h_i$ attached one after the other as in Figure \ref{Hmin}, where we see the junction of two consecutive branches $h_i$ and $h_j$.

Next proposition concerns the first branch of $H_{VW}^{\max}$ and $H_{VW}^{\min}$. Here, we use the Assumption NDR.

\begin{proposition}\label{monbr0}
	The first branch of $H^{\max}_{VW}$ is monotonically increasing and the first branch of $H^{\min}_{VW}$ is monotonically decreasing.
\end{proposition}
\begin{proof}
	Consider $x_0\in (-\pi,\pi]$ such that $x_0$ is a stationary point of $H_{VW}(\cdot,0)$. A branch $x_i(\beta)$, $\beta$ in a right neighborhood of zero, such that $x_i(0)=x_{0}$ exists if
	\begin{equation*}
	\frac{\partial^2 H_{VW}}{\partial x^2}(x_0,0)\neq 0 
	\end{equation*}
	(recall the previous subsection).  By (\ref{dbeta}), the corresponding branch $h_i$ has the derivative
	\begin{equation*}
		h_{i}^\prime(0)=\frac{f^{\max}_{VW}(x_0)V\sin x_0}{(1+V\cos x_0)^2},
	\end{equation*}
	where $f^{\max}_{VW}( x_0)>0$. 
	
	Therefore, the branch $h_i$ is monotonically increasing if $\sin  x_0>0$ and it is monotonically decreasing if $\sin  x_0<0$.
	
	If $\sin x_0=0$, i.e., $x_0=0$ or $x_0=\pi$, to determine whether the branch $h_i$ is monotonically increasing or decreasing, we examine the sign of $\sin\!\left(x_{i}(\beta)+\beta\right)$ in (\ref{dbeta}), for $\beta$ in a right neighborhood of $0$: $h_i$ is monotonically increasing if the sign is positive and it is monotonically decreasing if the sign is negative. This sign is determined by
	$$
	x_{i}(\beta)+\beta|_{\beta=0}=x_0,
	$$
	which is $0$ or $\pi$, and the sign of
	\begin{equation}
		\frac{d}{d\beta}\left(x_{i}(\beta)+\beta\right)|_{\beta=0}=x_{i}^\prime(0)+1 \label{derivative}
	\end{equation}
		with the following table:
	\begin{equation}
		\begin{tabular}{|c|c|c|}
			\hline
			\makecell[l]{$x_{i}(\beta)+\beta|_{\beta=0}$\\ $=x_0$ }& \makecell[l]{$\frac{d}{d\beta}\left(x_{i}(\beta)+\beta\right)|_{\beta=0}$\\ $=x_{i}^\prime(0)+1$}& \makecell[l]{$\sin(x_{i}(\beta)+\beta)$ for $\beta>0$ in a \\ right neighborhood of zero}\\
			\hline
			$0$ & positive & positive \\
			\hline
			$0$ & negative &  negative \\
			\hline
			$\pi$ & positive & negative \\
			\hline
			$\pi$ & negative &   positive \\
			\hline
		\end{tabular}
		\label{tablesigns}
	\end{equation}
	Observe that in (\ref{derivative}), we have (see (\ref{odexi}) and pass to the limit as $\beta\rightarrow 0$)
	\begin{equation}
		x_{i}^\prime(0)+1=\frac{\frac{\partial^2 H_{VW}}{\partial x^2}(x_0,0)-\frac{\partial^2 H_{VW}}{\partial \beta\partial x}(x_0,0)}{\frac{\partial^2 H_{VW}}{\partial x^2}(x_0,0)}, \label{xiprime+1}
	\end{equation}
	where
	\begin{equation}
		\frac{\partial^2 H_{VW}}{\partial x^2}(x_0,0)-\frac{\partial^2 H_{VW}}{\partial \beta\partial x}(x_0,0)=\frac{\left(f^{\max}_{VW}\right)^{\prime\prime}(x_0)}{1+V\cos x_0} \label{partialminus}
	\end{equation}
	(see (\ref{dH}) and recall that $x_0=0$ or $x_0=\pi$).
	
	The proof considers four cases.

	    \underline{Case 1}
		
		Consider the function $H_{VW}^{\max}$ in the case $Q<1$, where $Q$ is defined in (\ref{QQ0}).
		
		For $x_0\in(-\pi,\pi]$, we have $H_{VW}^{\max}(0)=H_{VW}(x_0,0)$ if and only if $\cos x_0=-Q$ (see Theorem \ref{lastThAC} of Appendix \ref{AC}), i.e.,
		$$
		x_0=x_{01}:=\arccos(-Q)\in\left(\frac{\pi}{2},\pi\right)\text{\ \ or\ \ }x_0=x_{02}:=-\arccos(-Q)\in\left(-\pi,-\frac{\pi}{2}\right).
		$$
Since (\ref{condneqzero}) holds for $x_0=x_{01}$ or $x_0=x_{02}$,  there exists a branch $x_{i_1}(\beta)$, $\beta$ in a right neighborhood of $0$, such that $x_{i_1}(0)=x_{01}$ and a branch $x_{i_2}(\beta)$, $\beta$ in a right neighborhood of $0$, such that $x_{i_2}(0)=x_{02}$.

Let $i^\ast$ be the index of the first branch of $H^{\max}_{VW}$. We have $i^\ast=i_1$ or $i^\ast=i_2$. For $x_0=x_{01}$, we have $\sin x_{0}>0$ and then the branch $h_{i_1}$ is monotonically increasing. For $x_0=x_{02}$, we have $\sin x_{0}<0$ and then the branch $h_{i_2}$ is monotonically decreasing. This shows that $i^\ast=i_1$ and then the branch $h_{i^{\ast}}$ is monotonically increasing.

For an example of this situation, see the first two rows of Figure \ref{branches}, where $Q<1$ holds.

    \underline{Case 2}
	
	Consider the function $H^{\max}_{VW}$ in the case $Q\geq 1$.
	
	For $x_0\in(-\pi,\pi]$, we have $H_{VW}^{\max}(0)=H_{VW}(x_0,0)$ if and only if $x_0=\pi$ (see Theorem \ref{lastThAC} of Appendix \ref{AC}). Since (\ref{condneqzero}) holds for $x_0=\pi$,  there exists a branch $x_{i}(\beta)$, $\beta$ in a right neighborhood of $0$, such that $x_{i}(0)=\pi$.
	
	Let $i^\ast$ be the index of the first branch of $H^{\max}_{VW}$. We have $i^\ast=i$. Since $x_0=\pi$ is a maximum of $H_{VW}(\cdot,0)$, we have
	\begin{equation}
		\frac{\partial^2 H_{VW}}{\partial x^2}(x_0,0)<0. \label{max<0}
	\end{equation}
	Moreover, we have (see table (\ref{tablefmax}) of Appendix \ref{AA})
	\begin{equation}
	\left(f^{\max}_{VW}\right)^{\prime\prime}(x_0)>0. \label{fmax>0}
	\end{equation}
	Then, by (\ref{max<0}), (\ref{fmax>0}) and (\ref{xiprime+1})-(\ref{partialminus}), we obtain
	$$
	x_{i^\ast}^\prime(0)+1<0
	$$
	and, then, since $x_0=\pi$, we have $\sin(x_{i^\ast}(\beta)+\beta)$ positive for $\beta$ in a right neighborhood of zero (see table (\ref{tablesigns})). Therefore, the branch $h_{i^\ast}$ is monotonically increasing.
	
	For an example of this situation, see the third row of Figure \ref{branches}, where $Q\geq 1$ holds.

		\underline{Case 3}
	
	Consider the function $H_{VW}^{\min}$ in the case $V<W$.
	
	For $x_0\in(-\pi,\pi]$, we have $H_{VW}^{\min}(0)=H_{VW}(x_0,0)$ if and only if $x_0=0$ or $x_0=\pi$ (see Theorem \ref{lastThAC} of Appendix \ref{AC}).  Since (\ref{condneqzero}) holds for $x_0=0$ or $x_0=\pi$,  there exists a branch $x_{i_1}(\beta)$, $\beta$ in a right neighborhood of $0$, such that $x_{i_1}(0)=0$ and a branch $x_{i_2}(\beta)$, $\beta$ in a right neighborhood of $0$, such that $x_{i_2}(0)=\pi$.
	
	Let $i^\ast$ be the index of the first branch of $H^{\min}_{VW}$. We have $i^\ast=i_1$ or $i^\ast=i_2$. Since, for $x_0=0$ or $x_0=\pi$, $x_0$ is a minimum of $H_{VW}(\cdot,0)$, we have
	\begin{equation}
		\frac{\partial^2 H_{VW}}{\partial x^2}(x_0,0)>0. \label{>>0}
	\end{equation}
	For $x_0=0$, we have (see table (\ref{tablefmax}) of Appendix \ref{AA})
	\begin{equation}
		\left(f^{\max}_{VW}\right)^{\prime\prime}(x_0)<0, \label{2>>>0}
	\end{equation}
	and then, by (\ref{>>0}), (\ref{2>>>0}) and (\ref{xiprime+1})-(\ref{partialminus}), we obtain
	$$
	x_{i_1}^\prime(0)+1<0.
	$$
	Then, since $x_0=0$, we have $\sin(x_{i_1}(\beta)+\beta)$ negative for $\beta$ in a right neighborhood of zero (see table (\ref{tablesigns})). Therefore, the branch $h_{i_1}$ is monotonically decreasing.  
	For $x_0=\pi$, we have (see table (\ref{tablefmax}) of Appendix \ref{AA})
	\begin{equation}
		\left(f^{\max}_{VW}\right)^{\prime\prime}(x_0)>0, \label{>>>0}
	\end{equation}
	and then, by (\ref{>>0}), (\ref{>>>0}) and (\ref{xiprime+1})-(\ref{partialminus}), we obtain
	$$
	x_{i_2}^\prime(0)+1>0.
	$$
	Then, since $x_0=\pi$, we have $\sin(x_{i_2}(\beta)+\beta)$ negative for $\beta$ in a right neighborhood of zero (see table (\ref{tablesigns})). Therefore, the branch $h_{i_2}$ is monotonically decreasing. Since both branches $h_{i_1}$ and $h_{i_2}$ are monotonically decreasing, the branch $h_{i^\ast}$, which is one between $h_{i_1}$ and $h_{i_2}$,  is monotonically decreasing. Indeeed, as it shown below in Remark \ref{Rembranches}, the two branches $h_{i_1}$ and $h_{i_2}$ coincide, although the two branches $x_{i_1}$ and $x_{i_2}$ are different.
	
	For an example of this situation, see the first row of Figure \ref{branches}, where $V<W$ holds.
	
	\underline{Case 4}
	
	Consider the function $H^{\min}_{VW}$ in the case $V>W$.
	
	For $x_0\in(-\pi,\pi]$, we have $H_{VW}^{\min}(0)=H_{VW}(x_0,0)$ if and only if $x_0=0$ (see Theorem \ref{lastThAC} of Appendix \ref{AC}).  Since (\ref{condneqzero}) holds for $x_0=0$,  there exists a branch $x_{i}(\beta)$, $\beta$ in a right neighborhood of $0$, such that $x_{i}(0)=0$.
	
	Let $i^\ast$ be the index of the first branch of $H^{\min}_{VW}$. We have $i^\ast=i$. Since $x_0=0$ is a minimum of $H_{VW}(\cdot,0)$, we have 
	\begin{equation}
		\frac{\partial^2 H_{VW}}{\partial x^2}(x_0,0)>0, \label{min<0} 
	\end{equation}
Moreover, we have  (see table (\ref{tablefmax}) of Appendix \ref{AA})
	\begin{equation}
	\left(f^{\max}_{VW}\right)^{\prime\prime}(x_0)<0. \label{fmax<0}
	\end{equation}
	Then, by (\ref{min<0}), (\ref{fmax<0}) and (\ref{xiprime+1})-(\ref{partialminus}), we obtain
	$$
	x_{i^\ast}^\prime(0)+1<0
	$$
	and, then, since $x_0=0$, we have $\sin(x_{i^\ast}(\beta)+\beta)$ negative for $\beta$ in a right neighborhood of zero (see table (\ref{tablesigns})). Therefore, the branch $h_{i^\ast}$ is monotonically decreasing.
	
	For an example of this situation, see the last two rows of Figure \ref{branches}, where $V>W$ holds.

\end{proof}

\begin{remark}\label{Rembranches}

	In the case $V<W$, the branch $x_i(\beta)$, $\beta$ in a right neighborhood of zero, such that $x_i(0)=0$ and the branch $x_i(\beta)$, $\beta$ in a right neighborhood of zero,  such that $x_i(0)=\pi$ correspond to a unique branch $h_i$.

	In fact, recalling Remark \ref{twosets}, $(x_0,0)$ with $x_0$ multiple of $\pi$ is a solution of (\ref{2pi}), not of (\ref{2piother}) (we have $\alpha^{\max}_{VW}(x_0)=0$: see table (\ref{tablesalpha}) of Appendix \ref{AE}). Since the branches $x_i(\beta)$ are continuous functions of $\beta$, this remains valid for $\beta$ in a right neighborhood of zero, and, then, for both branches the corresponding branch $h_i$ is given by (\ref{hiWcos}).
	
	For an example of this situation, see the left side of the first row of Figures \ref{branches} and \ref{branchespartial}, where $V<W$ holds.

\end{remark}
\begin{theorem} \label{monbr}
 The functions $H_{VW}^{\max}$ and $H_{VW}^{\min}$, considered in $[0,\pi]$, are monotonically increasing and decreasing, respectively.
\end{theorem}

\begin{proof}
	We prove that the branches forming the functions $H_{VW}^{\max}$ in $[0,\pi]$ are all monotonically increasing and the branches forming the functions $H_{VW}^{\min}$ in $[0,\pi]$ are all monotonically decreasing.  The proof is given for $H_{VW}^{\max}$. The proof for $H_{VW}^{\min}$ is similar.
	
	We show that if a branch forming $H_{VW}^{\max}$ is monotonically increasing, then the successive branch is also monotonically increasing. Since  the previous Proposition \ref{monbr0} states that the first branch of $H_{VW}^{\max}$ is monotonically increasing, it follows that all the branches of $H_{VW}^{\max}$ are monotonically increasing. 
	
	Consider two consecutive branches $h_i$ and $h_j$ of $H_{VW}^{\max}$, with $h_i$ being monotonically increasing: see Figure \ref{Hmin}. Let $\beta^{\ast}\in(0,\pi)$ denote the point where these two branches meet. Since $D_i$ and $D_j$ are closed (recall Subsection \ref{Brxihi}), we have $\beta^\ast\in D_i\cap D_j$.
	
	By Proposition  \ref{diff}, we have
	$$
	h_i^\prime(\beta^\ast)=\frac{\partial H_{VW}}{\partial \beta}(x_i(\beta^\ast),\beta^\ast).
	$$
	By Proposition \ref{diffneq0} and the fact that $h_i$ is monotonically increasing, we obtain
	$$
	\frac{\partial H_{VW}}{\partial \beta}(x_i(\beta^\ast),\beta^\ast)=h_i^\prime(\beta^\ast)>0.
	$$
	This shows that, for $\beta>\beta^\ast$ in a right neighborhood of $\beta^\ast$, we have
	$$
	H_{VW}(x_i(\beta^\ast),\beta)>H_{VW}(x_i(\beta^\ast),\beta^\ast)
	$$
	and then
	$$
	H_{VW}^{\max}(\beta)\geq H_{VW}(x_i(\beta^\ast),\beta)>H_{VW}(x_i(\beta^\ast),\beta^\ast)=H_{VW}^{\max}(\beta^\ast).
	$$
	This shows that the branch $h_j$ is monotonically increasing.

\end{proof}

Since the $2\pi$-periodic functions $H^{\max}_{VW}$ and $H^{\min}_{VW}$ satisfy the property (\ref{propHmaxmin}), by the previous Theorem \ref{monbr} we conclude that $H^{\max}_{VW}$ and $H^{\min}_{VW}$ oscillate monotonically between their extreme values, with $H^{\max}_{VW}$ ($H^{\min}_{VW}$) monotonically increasing (decreasing) in $[0,\pi]$ and monotonically decreasing (increasing) in $[\pi,2\pi]$.

To complete the proof of Theorem \ref{ThAB2}, it remains to determine the maximum and minimum values of $H_{VW}^{\max}$ and $H_{VW}^{\min}$.

Since $H_{VW}^{\max}$ satisfies the property (\ref{propHmaxmin}) and it is monotonically increasing in $[0,\pi]$, we have (see Theorem \ref{lastThAC} of Appendix \ref{AC})
	\begin{equation*}
		\min\limits_{\beta\in\mathbb{R}}H_{VW}^{\max}(\beta)=\min\limits_{\beta\in[0,\pi]}H_{VW}^{\max}(\beta)=H_{VW}^{\max}(0)=\left\{
		\begin{array}{l}
			\frac{1-V^2}{(1-QV)(1-W)}\text{\ if\ }Q\leq 1\\
			\\
			\frac{1+V}{(1-V)(1+W)}\text{\ if\ }Q\geq 1
		\end{array}
		\right.,
	\end{equation*}
	where $Q$ is defined in (\ref{QQ0}), and
	\begin{eqnarray*}
	\max\limits_{\beta\in\mathbb{R}}H_{VW}^{\max}(\beta)&=&\max\limits_{\beta\in[0,\pi]}H_{VW}^{\max}(\beta)=H_{VW}^{\max}(\pi)\\
	&=&\max\limits_{x\in\mathbb{R}}H_{VW}(x,\pi)=\max\limits_{x\in\mathbb{R}}\frac{f^{\max}_{VW}(x)}{1-V\cos x}=\frac{1+V}{(1-V)(1-W)}
	\end{eqnarray*}
	since $f_{VW}^{\max}(x)$ achieves its maximum value
	$
	\frac{1+V}{1-W},
	$
	as $x$ varies, when $x$ is an even multiple of $\pi$ (see Theorem \ref{maxminf} of Appendix \ref{AA}), and $1-V\cos x$ achieves its minimum value $1-V$, as $x$ varies, when $x$ is an even multiple of $\pi$.
	
	Moreover, since $H_{VW}^{\min}$ satisfies the property (\ref{propHmaxmin}) and it is monotonically decreasing in $[0,\pi]$, we have
	\begin{equation*}
		\max\limits_{\beta\in\mathbb{R}}H_{VW}^{\min}(\beta)=\max\limits_{\beta\in[0,\pi]}H_{VW}^{\min}(\beta)=H_{VW}^{\min}(0)=\frac{1}{1-W}
	\end{equation*}
	(see Theorem \ref{lastThAC} of Appendix \ref{AC}) and
	\begin{eqnarray*}
	\min\limits_{\beta\in\mathbb{R}}H_{VW}^{\min}(\beta)&=&\min\limits_{\beta\in[0,\pi]}H_{VW}^{\min}(\beta)=H_{VW}^{\min}(\pi)\\
	&=&\min\limits_{x\in\mathbb{R}}H_{VW}(x,\pi)=\min\limits_{x\in\mathbb{R}}\frac{f^{\max}_{VW}(x)}{1-V\cos x}=\left\{
			\begin{array}{l}
					\frac{1-V}{(1+V)(1-W)}\text{\ if\ }V<W\\
					\\
					\frac{1}{1+W}\text{\ if\ }V>W,
				\end{array}
		\right.
	\end{eqnarray*}
	since $f_{VW}^{\max}(x)$ achieves its minimum value
	$$
	\left\{
	\begin{array}{l}
		\frac{1-V}{1-W}\text{\ if\ }V<W\\
		\\
		\frac{1+V}{1+W}\text{\ if\ }V>W.
	\end{array}
	\right.
	$$
	as $x$ varies, when $x$ is an odd multiple of $\pi$ (see Theorem \ref{maxminf} of Appendix \ref{AA}), and $1-V\cos x$ achieves its maximum value $1+V$, as $x$ varies, when $x$ is an odd multiple of $\pi$.

\section{}\label{AD}
In this appendix, we show that equation (\ref{equationfVW}) of Appendix \ref{AA} has solutions.
\begin{proposition}\label{arcsin}
	Let $V,W\in[0,1)$, let $x\in\mathbb{R}$ and let
	$$
	U_{VW}(x)=V\mathrm{e}^{\mathrm{i}x}+W.
	$$
	Suppose $U_{VW}(x)\neq 0$. We have
	\begin{equation}
		\left\vert\frac{VW\sin x}{\vert U_{VW}(x)\vert}\right\vert< 1.  \label{AppAA}
	\end{equation}
\end{proposition}
\begin{proof}
	The inequality (\ref{AppAA}) holds if $V=0$ or $W=0$. Therefore, now we suppose $V,W\neq 0$.
	
	Observe that
	$$
	U_{VW}(x)=V\mathrm{e}^{\mathrm{i} x}+W=V\cos x+W+\mathrm{i}V\sin x
	$$
	has squared modulus
	$$
	\vert U_{VW}(x) \vert^2=V^2+2VW\cos x+W^2.
	$$
	Hence, (\ref{AppAA}) is equivalent to
	$$
	V^2W^2\sin^2 x<V^2+2VW\cos x+W^2,
	$$
	i.e.,
	\begin{equation}
		V^2W^2\cos^2 x+2VW\cos x+V^2+W^2-V^2W^2>0.\label{ineqVW}
	\end{equation}

	The trinomial in the variable $y$
	$$
	V^2W^2y^2+2VWy+V^2+W^2-V^2W^2
	$$
	has discriminant
	$$
	\frac{\Delta}{4}=V^2W^2-V^2W^2(V^2+W^2-V^2W^2)=V^2W^2(1-V^2)(1-W^2)>0
	$$
	and then it has positive values for
	\begin{equation}
		y< y_{-}\text{\ \ or\ \ }y> y_{+},  \label{posval}
	\end{equation}
	where
	$$
	y_{-}=\frac{-1-\sqrt{(1-V^2)(1-W^2)}}{VW}
	$$
	and
	$$
	y_{+}=\frac{-1+\sqrt{(1-V^2)(1-W^2)}}{VW}.
	$$
	are the zeros of the trinomial. We have
	$$
	y_{-}<y_{+}\leq -1.
	$$
	In fact, the condition $y_{+}\leq -1$ is equivalent to
	$$
	\sqrt{(1-V^2)(1-W^2)}\leq 1-VW,
	$$
	i.e.,
	$$
	(1-V^2)(1-W^2)\leq (1-VW)^2
	$$
	that reduces to
	$$
	0\leq (V-W)^2.
	$$
	
	The condition (\ref{posval}), in particular $y> y_{+}$,  is satisfied for $y=\cos x$. In fact, we have $y=\cos x\geq -1\geq y_{+}$. Moreover, we have  $y=\cos x=-1=y_{+}$ if and only if $x$ is an odd multiple of $\pi$ and $V=W$, a case where we have  $U_{VW}(x)=0$, which is excluded.
	
	Since  (\ref{posval}) holds for $y=\cos x$, we have (\ref{ineqVW}) and then (\ref{AppAA}).
\end{proof}

\section{}\label{AE}

In this appendix, we study the functions $\alpha_{VW}^{\max}$ and $\alpha_{VW}^{\min}$ defined in Theorem \ref{LemmafvW} of Appendix \ref{AA}.

Let $V,W\in[0,1)$ with $V\neq 0$ or $W\neq 0$, let $x\in\mathbb{R}$ and let
\begin{equation}
	U_{VW}(x)=V\mathrm{e}^{\mathrm{i}x}+W=V\cos x +W+\mathrm{i}V\sin x. \label{Ux}
\end{equation}
Consider the $2\pi$-periodic functions
	\begin{equation}
	\alpha^{\max}_{VW}(x)=\arcsin\frac{VW\sin x}{\vert U_{VW}(x)\vert}-\theta_{VW}(x)
	\label{alphamax}
\end{equation}
and
\begin{equation}
	\alpha^{\min}_{VW}(x)=\pi-\arcsin\frac{VW\sin x}{\vert U_{VW}(x)\vert}-\theta_{VW}(x)
	\label{alphamin}
\end{equation}
where $\theta_{VW}(x)$ is the angle in $(-\pi,\pi]$ of the polar form of $U_{VW}(x)$. The domain of these functions is $\mathbb{R}$ if $V\neq W$, and $\mathbb{R}$ without the odd multiples of $\pi$ if $V=W$.

The function $\theta_{VW}(x)$ is discontinuos for $x$ such that $\theta_{VW}(x)$ crosses $-\pi$ or $\pi$. This happens only if (see (\ref{Ux}))
$$
\sin x=0\text{\ \ and\ \ }V\cos x +W<0,
$$
i.e., $x$ is an odd multiple of $\pi$ and $V>W$. We conclude that that $\theta_{VW}$ is continuous on $\mathbb{R}$ when $V<W$, and it is continuous on $(-\pi,\pi)$ when $V\geq W$. Consequently, the same is true for  $\alpha^{\max}_{VW}(x)$ and $\alpha^{\min}_{VW}(x)$.

Indeed, $\theta_{VW}$ is not only continuous, but it is real analytic on $\mathbb{R}$ when $V<W$, and real analytic on $(-\pi,\pi)$ when $V\geq W$. Consequently, the same is true for  $\alpha^{\max}_{VW}$ and $\alpha^{\min}_{VW}$.

The derivatives of the functions $\alpha^{\max}_{VW}$ and $\alpha^{\min}_{VW}$ are solutions of the differential equations 
obtained by deriving with respect to $x$ the equations
$$
\frac{\partial f}{\partial \alpha}\left(\alpha^{\max}_{VW}(x),x\right)=0\text{\ \ and\ \ }\frac{\partial f}{\partial \alpha}\left(\alpha^{\min}_{VW}(x),x\right)=0.
$$
with $x\in\mathbb{R}$ if $V<W$ and $x\in\mathbb{R}\setminus\{k\pi:k\in\mathbb{Z}\text{\ is odd}\}$ if $V\geq W$.
We have
$$
0=\frac{d}{dx}\frac{\partial f}{\partial \alpha}\left(\alpha^{\max}_{VW}(x),x\right)=\frac{\partial^2 f}{\partial \alpha^2}\left(\alpha^{\max}_{VW}(x),x\right)\left(\alpha^{\max}_{VW}\right)^\prime
(x)+\frac{\partial^2 f}{\partial x\partial  \alpha}\left(\alpha^{\max}_{VW}(x),x\right)
$$
where (see the end of proof of Theorem \ref{LemmafvW} in Appendix \ref{AA})
$$
\frac{\partial^2 f}{\partial \alpha^2}\left(\alpha^{\max}_{VW}(x),x\right)<0.
$$
Thus, we have
\begin{equation}
\left(\alpha^{\max}_{VW}\right)^\prime
(x)=-\frac{\frac{\partial^2 f}{\partial\ x\partial  \alpha}\left(\alpha^{\max}_{VW}(x),x\right)}{\frac{\partial^2 f}{\partial \alpha^2}\left(\alpha^{\max}_{VW}(x),x\right)}=\frac{-V\cos(x+\alpha^{\max}_{VW}(x))+VW\cos x}{V\cos(x+\alpha^{\max}_{VW}(x))+W\cos \alpha^{\max}_{VW}(x)}
\label{alphamaxprime}
\end{equation}
(see (\ref{df}) in Appendix \ref{AA}). Analogously, we have
\begin{equation}
\left(\alpha^{\min}_{VW}\right)^\prime
(x)=-\frac{\frac{\partial^2 f}{\partial\ x\partial  \alpha}\left(\alpha^{\min}_{VW}(x),x\right)}{\frac{\partial^2 f}{\partial \alpha^2}\left(\alpha^{\min}_{VW}(x),x\right)}=\frac{-V\cos(x+\alpha^{\min}_{VW}(x))+VW\cos x}{V\cos(x+\alpha^{\min}_{VW}(x))+W\cos \alpha^{\min}_{VW}(x)}.
\label{alphaminprime}
\end{equation}

When $V\geq W$, the restrictions of the functions $\alpha^{\max}_{VW}$ and $\alpha^{\min}_{VW}$  to $(-\pi,\pi)$ can be extended by continuity to $-\pi$ and $\pi$, since the limits of  $\alpha_{VW}^{\max}$ and $\alpha_{VW}^{\min}$ exist and are finite, and these extensions are continuously differentiable there, since the limits of the derivatives of $\alpha_{VW}^{\max}$ and $\alpha_{VW}^{\min}$ exist and are finite. Therefore, the restrictions of the functions $\alpha^{\max}_{VW}$ and $\alpha^{\min}_{VW}$  to $(-\pi,\pi)$ can be extended to $C^1$ functions on $[-\pi,\pi]$.

In fact, suppose $V>W$, so that $U_{VW}(x)\neq 0$ for $x=-\pi$ or $x=\pi$. We have
\begin{eqnarray*}
&&\lim\limits_{x\rightarrow -\pi^{+}}\vert U_{VW}(x)\vert=\lim\limits_{x\rightarrow \pi^{-}}\vert U_{VW}(x)\vert=\vert U_{VW}(-\pi)\vert=\vert U_{VW}(\pi)\vert=V-W\\
&&\lim\limits_{x\rightarrow -\pi^{+}}\theta_{VW}(x)=-\pi\text{\ \ and\ \ }\lim\limits_{x\rightarrow \pi^{-}}\theta_{VW}(x)=\theta_{VW}(\pi)=\pi
\end{eqnarray*}
and then by (\ref{alphamax}) and (\ref{alphamin}) we have
\begin{eqnarray}
	&&\lim\limits_{x\rightarrow -\pi^{+}}\alpha_{VW}^{\max}(x)=\pi\text{\ \ and\ \ }\lim\limits_{x\rightarrow \pi^{-}}\alpha_{VW}^{\max}(x)=\alpha_{VW}^{\max}(\pi)=-\pi\notag \\
	&&\lim\limits_{x\rightarrow -\pi^{+}}\alpha_{VW}^{\min}(x)=2\pi\text{\ \ and\ \ }\lim\limits_{x\rightarrow \pi^{-}}\alpha_{VW}^{\min}(x)=\alpha_{VW}^{\min}(\pi)=0. \notag\\
	\label{limalphamaxmin}
\end{eqnarray}
Regarding the derivatives, by (\ref{alphamaxprime}), (\ref{alphaminprime}) and (\ref{limalphamaxmin}) we have
\begin{eqnarray*}
	&&\lim\limits_{x\rightarrow -\pi^{+}}\left(\alpha_{VW}^{\max}\right)^\prime(x)=\lim\limits_{x\rightarrow \pi^{-}}\left(\alpha_{VW}^{\max}\right)^\prime(x)=-\frac{V(1+W)}{V-W}\\
	&&\lim\limits_{x\rightarrow -\pi^{+}}\left(\alpha_{VW}^{\min}\right)^\prime(x)=\lim\limits_{x\rightarrow \pi^{-}}\left(\alpha_{VW}^{\min}\right)^\prime(x)=-\frac{V(1-W)}{V-W}.
\end{eqnarray*}

It is unsurprising that the values $\left(\alpha_{VW}^{\max}\right)(x)$ and $\left(\alpha_{VW}^{\min}\right)(x)$ differ by $2\pi$ at $x=\pm \pi$, and the values of the derivatives $\left(\alpha_{VW}^{\max}\right)^\prime(x)$ and $\left(\alpha_{VW}^{\min}\right)^{\prime}(x)$ at $x=\pm \pi$ coincide. In fact, when the polar angle function $\theta_{VW}$ in (\ref{alphamax}) and (\ref{alphamin}) is regarded as taking values in the torus $\mathbb{R}/(2\pi\mathbb{Z})$, rather than in $(-\pi,\pi]$, it is real analytic on $\mathbb{R}$ and $2\pi$-periodic modulo $2\pi$. Consequently, also  $\alpha_{VW}^{\max}$ and $\alpha_{VW}^{\min}$ become real analytic on $\mathbb{R}$ and $2\pi$-periodic modulo $2\pi$.  

Next tables contain the values of  $\alpha_{VW}^{\max}$, $\alpha_{VW}^{\min}$ and the derivatives $\left(\alpha_{VW}^{\max}\right)^\prime$, $\left(\alpha_{VW}^{\min}\right)^\prime$ at $-\pi^{+}$, $0$ and $\pi$, for the cases $V<W$ and $V>W$.
\begin{eqnarray}
&&\text{Case $V<W$}: \begin{tabular}{|l|l|l|l|l|}
	\hline
	$x$ & $\alpha_{VW}^{\max}(x)$ & $\left(\alpha_{VW}^{\max}\right)^\prime(x)$ & $\alpha_{VW}^{\min}(x)$ & $\left(\alpha_{VW}^{\min}\right)^\prime(x)$ \\
	\hline
	$-\pi$ & $0$  & $\frac{V(1-W)}{W-V}$ & $\pi$ & $\frac{V(1+W)}{W-V}$\\
	\hline 
	$0$ & $0$ & $-\frac{V(1-W)}{V+W}$ & $\pi$ & $-\frac{V(1+W)}{V+W}$\\
	\hline
		$\pi$ & $0$  & $\frac{V(1-W)}{W-V}$ & $\pi$ & $\frac{V(1+W)}{W-V}$\\
	\hline
\end{tabular}\notag\\
\notag\\
&&\text{Case $V>W$}: \begin{tabular}{|l|l|l|l|l|}
	\hline
	$x$ & $\alpha_{VW}^{\max}(x)$ & $\left(\alpha_{VW}^{\max}\right)^\prime(x)$ & $\alpha_{VW}^{\min}(x)$ & $\left(\alpha_{VW}^{\min}\right)^\prime(x)$ \\
	\hline
	$-\pi^{+}$ & $\pi$  & $-\frac{V(1+W)}{V-W}$ & $2\pi$ & $-\frac{V(1-W)}{V-W}$\\
	\hline 
	$0$ & $0$ & $-\frac{V(1-W)}{V+W}$ & $\pi$ & $-\frac{V(1+W)}{V+W}$\\
	\hline
	$\pi^{-}$ & $-\pi$  & $-\frac{V(1+W)}{V-W}$ & $0$ & $-\frac{V(1-W)}{V-W}$\\
	\hline
\end{tabular}\notag\\
\label{tablesalpha}
\end{eqnarray}

Now, suppose $V=W$. For $x\in(-\pi,\pi)$, we have
\begin{eqnarray*}
\vert U_{VV}(x)\vert=V\sqrt{2(1+\cos x)}\text{\ \ and\ \ }\theta_{VV}(x)=\arcsin \frac{\sin x}{\sqrt{2(1+\cos x)}}
\end{eqnarray*}
(for the expression of $\theta_{VV}(x)$ observe that $\mathrm{Re}(U_{VV}(x))=V(1+\cos x)>0$). Therefore, since
$$
\frac{\sin x}{\sqrt{2(1+\cos x)}}=\sin\frac{x}{2},
$$
we obtain
\begin{equation*}
	\alpha^{\max}_{VV}(x)=\arcsin\left(V\sin\frac{x}{2}\right)-\frac{x}{2}\text{\ \ and\ \ }\alpha^{\min}_{VV}(x)=\pi-\arcsin\left(V\sin\frac{x}{2}\right)-\frac{x}{2}
\end{equation*}
and then
\begin{equation*}
	\lim\limits_{x\rightarrow -\pi^{+}}\alpha^{\max}_{VV}(x)=-\arcsin V+\frac{\pi}{2}\text{\ \ and\ \ }\lim\limits_{x\rightarrow \pi^{-}}\alpha^{\max}_{VV}(x)=\arcsin V-\frac{\pi}{2}
\end{equation*}
and
\begin{equation*}
	\lim\limits_{x\rightarrow -\pi^{+}}\alpha^{\min}_{VV}(x)=\arcsin V+\frac{3\pi}{2}\text{\ \ and\ \ }\lim\limits_{x\rightarrow \pi^{-}}\alpha^{\min}_{VV}(x)=-\arcsin V+\frac{\pi}{2}
\end{equation*}
Regarding the derivatives, we have, for $x\in(-\pi,\pi)$, 
\begin{equation*}
	\left(\alpha^{\max}_{VV}\right)^{\prime}(x)=\frac{\frac{V}{2}\cos\frac{x}{2}}{\sqrt{1-V^2\sin^2\frac{x}{2}}}-\frac{1}{2}\text{\ \ and\ \ }\left(\alpha^{\min}_{VV}\right)^\prime(x)=-\frac{\frac{V}{2}\cos\frac{x}{2}}{\sqrt{1-V^2\sin^2\frac{x}{2}}}-\frac{1}{2},
\end{equation*}
and then
\begin{equation*}
	\lim\limits_{x\rightarrow -\pi^{+}}\left(\alpha^{\max}_{VV}\right)^\prime(x)=-\frac{1}{2}\text{\ \ and\ \ }\lim\limits_{x\rightarrow  \pi^{-}}\left(\alpha^{\max}_{VV}\right)^\prime(x)=-\frac{1}{2}.
\end{equation*}
and
\begin{equation*}
	\lim\limits_{x\rightarrow -\pi^{+}}\left(\alpha^{\min}_{VV}\right)^\prime(x)=-\frac{1}{2}\text{\ \ and\ \ }\lim\limits_{x\rightarrow  \pi^{-}}\left(\alpha^{\min}_{VV}\right)^\prime(x)=-\frac{1}{2}.
\end{equation*}

Similarly to the case $V>W$, the values of the derivatives $\left(\alpha_{VV}^{\max}\right)^\prime(x)$ and $\left(\alpha_{VV}^{\min}\right)^{\prime}(x)$ at $x=\pm \pi$ are equal. Instead, the values $\left(\alpha_{VV}^{\max}\right)(x)$ and $\left(\alpha_{VV}^{\min}\right)(x)$ differ by $\pi-2\mathrm{arcsin}V$ at $x=\pm \pi$, not of $2\pi$ as in the case $V>W$.

Next table completes (\ref{tablesalpha}) with the case $V=W$.
\begin{eqnarray}
&&\text{Case $V=W$}:\notag\\
&&\notag \\
&&\begin{tabular}{|c|c|c|c|c|}
	\hline
	$x$ & $\alpha_{VW}^{\max}(x)$ & $\left(\alpha_{VW}^{\max}\right)^\prime(x)$ & $\alpha_{VW}^{\min}(x)$ & $\left(\alpha_{VW}^{\min}\right)^\prime(x)$ \\
	\hline
	$-\pi^{+}$ & $-\arcsin V+\frac{\pi}{2}$  & $-\frac{1}{2}$ & $\arcsin V+\frac{3\pi}{2}$ & $-\frac{1}{2}$\\
	\hline 
	$0$ & $0$ & $-\frac{1-V}{2}$ & $\pi$ & $-\frac{1+V}{2}$\\
	\hline
	$\pi^{-}$ & $\arcsin V-\frac{\pi}{2}$  & $-\frac{1}{2}$ & $-\arcsin V+\frac{\pi}{2}$ & $-\frac{1}{2}$\\
	\hline
\end{tabular}\ \ .\notag\\
\label{tablealpha2}
\end{eqnarray}

\section{} \label{ACMAXMIN}

In this appendix we present a result concerning the continuity of maximum and minimum values, which is used more times in Appendices \ref{AA} and \ref{AB}.

\begin{theorem}\label{Contmax}
	Let $g:A\times B\rightarrow \mathbb{R}$, where $A$ is an arbitrary non-empty set and $B$ is a metric space. Suppose that, for any $b\in B$, there exists
	$$
	\max\limits_{a\in A}g(a,b)
	$$
	and consider $g^{\max}:B\rightarrow \mathbb{R}$ given by
	$$
	g^{\max}(b)=\max\limits_{a\in A}g(a,b),\ b\in B.
	$$
	If the function $g$ is Lipschitz continuous of constant $L$ with respect to the argument in $B$, then the function $g^{\max}$ is Lipschitz continuous of constant $L$.
	
\end{theorem}

\begin{proof}
	Suppose that $g$ is Lipschitz continuous of constant $L$ with respect to the argument in $B$. Let $b_1,b_2\in B$ and let $a_1,a_2\in A$ be such that
	$$
	g^{\max}(b_1)=g(a_1,b_1)\text{\ \ and\ \ }g^{\max}(b_2)=g(a_2,b_2).
	$$
	We have
	\begin{eqnarray*}
		g^{\max}(b_1)-g^{\max}(b_2)&=&g(a_1,b_1)-g(a_2,b_2)\\
		&=&g(a_1,b_1)-g(a_1,b_2)+\underset{\leq 0}{\underbrace{g(a_1,b_2)-g(a_2,b_2)}}\\
		&\leq& g(a_1,b_1)-g(a_1,b_2)\leq \vert g(a_1,b_1)-g(a_1,b_2)\vert\\
		&\leq& L d_B(b_1,b_2),
	\end{eqnarray*}
	where $d_B$ is the distance in the metric space $B$. By exchanging the role of $b_1$ and $b_2$, we obtain
	$$
	g^{\max}(b_2)-g^{\max}(b_1) \leq Ld_B(b_2,b_1)=Ld_B(b_1,b_2).
	$$
	Thus
	$$
	\vert g^{\max}(b_1)-g^{\max}(b_2)\vert\leq Ld_B(b_1,b_2).
	$$
	
\end{proof}
By replacing $g$ with $-g$, we see that the theorem also holds when $\max$ is replaced with $\min$.

Theorem \ref{Contmax} is used more times. Below, we list the four applications of the theorem. 

\subsection{Lipschitz continuity of $f^{\max}_{VW}$ and $f^{\min}_{VW}$}
Theorem \ref{Contmax} is applied at the beginning of Subsection \ref{fVWmaxmin} of Appendix \ref{AA}, in both $\max$ and $\min$ versions,  to the function $g=f_{VW}:\mathbb{R}^2\rightarrow \mathbb{R}$ given in (\ref{fVW}). Since the partial derivative $\frac{\partial f_{VW}}{\partial x}$ given in (\ref{df}) is bounded in $\mathbb{R}^2$, the function $f_{VW}$ is Lipschitz continuous with respect to the variable $x$. We conclude that the functions $f^{\max}_{VW}$ and $f^{\min}_{VW}$ in  (\ref{fmax}) are Lipschitz continuous.

\subsection{Case $V=W$ in the proof of Theorem \ref{maxminf} of Appendix \ref{AA}}\label{caseV=W}

Theorem \ref{Contmax} is applied, in both $\max$ and $\min$ versions, to show that in the proof of Theorem \ref{maxminf} of Appendix \ref{AA}, we can include the case $V=W\neq 0$ by a continuity argument. Here is how this is done.

Let $\overline{V}\in (0,1)$. We consider the function $g_1:\mathbb{R}^2\times B\rightarrow \mathbb{R}$, where $B$ is an open ball of center $(\overline{V},\overline{V})$ whose closure is strictly included in $(0,1)^2$, given by
$$
g_1(\alpha,x,V,W)=f_{VW}(\alpha,x)=\frac{1+V\cos(x+\alpha)}{1-W\cos\alpha},\ (\alpha,x,V,W)\in \mathbb{R}^2\times B.
$$
Since the partial derivatives (for the first see (\ref{df}) in Appendix \ref{AA})
\begin{eqnarray*}
	&&\frac{\partial g_1}{\partial x}(\alpha,x,V,W)=\frac{\partial f_{VW}}{\partial x}\left(\alpha,x\right)=\frac{-V\sin(x+\alpha)}{1-W\cos \alpha}\\
	&&
	\frac{\partial g_1}{\partial V}(\alpha,\beta,V,W)
	=\frac{\cos(x+\alpha)}{1-W\cos \alpha}\\
	&&
	\frac{\partial g_1}{\partial W}(\alpha,\beta,V,W)
	=\frac{(1+V\cos(x+\alpha))\cos \alpha}{\left(1-W\cos \alpha\right)^2}\\
	&&(\alpha,x,V,W)\in \mathbb{R}^2\times B
\end{eqnarray*}
are bounded in $\mathbb{R}^2\times B$, $g_1$ is Lipschitz continuous with respect to $(x,V,W)$ and then, by Theorem \ref{Contmax} applied to the function $g_1$, we conclude that the functions $g_2,g_3:\mathbb{R}\times B\rightarrow \mathbb{R}$ given by
\begin{eqnarray*}
&&g_2(x,V,W)=\max\limits_{\alpha\in\mathbb{R}}g_1(\alpha,x,V,W)=f_{VW}^{\max}(x)\\
&&g_3(x,V,W)=\min\limits_{\alpha\in\mathbb{R}}g_1(\alpha,x,V,W)=f_{VW}^{\min}(x)\\
&&(x,V,W)\in \mathbb{R}\times B 
\end{eqnarray*}
are Lipschitz continuous.

Now, consider a sequence $\left\{\left(V^{(n)},W^{(n)}\right)\right\}$ in $B$ with $V^{(n)}\neq W^{(n)}$ for any $n$ and
$$
\left(V^{(n)},W^{(n)}\right)\rightarrow \left(\overline{V},\overline{V}\right),\ n\rightarrow \infty.
$$
Since $g_i$, $i\in\{2,3\}$, is a continuous function, we have, for any $x\in\mathbb{R}$. 
$$
g_i\left(x,\overline{V},\overline{V}\right)=\lim\limits_{n\rightarrow \infty}g_i\left(x,V^{(n)},W^{(n)}\right)
$$
Therefore, for $x_1,x_2\in\mathbb{R}$ such that $x_1\leq x_2$, if
$$
g_i\left(x_1,V^{(n)},W^{(n)}\right)\leq g_i\left(x_2,V^{(n)},W^{(n)}\right)\text{\ \ for any $n$}
$$
then
$$
g_i\left(x_1,\overline{V},\overline{V}\right)\leq g_i\left(x_2,\overline{V},\overline{V}\right),
$$
and if
$$
g_i\left(x_1,V^{(n)},W^{(n)}\right)\geq g_i\left(x_2,V^{(n)},W^{(n)}\right)\text{\ \ for any $n$}
$$
then
$$
g_i\left(x_1,\overline{V},\overline{V}\right)\geq g_i\left(x_2,\overline{V},\overline{V}\right).
$$

Since, for any $n$, $V^{(n)}\neq 0$, $W^{(n)}\neq 0$ and $V^{(n)}\neq W^{(n)}$, proof of Theorem \ref{maxminf}{maxminf} of Appendix \ref{AA} shows that $f^{\max}_{V^{(n)}W^{(n)}}$ is monotonically decreasing (increasing) and $f^{\min}_{V^{(n)}W^{(n)}}$ is monotonically increasing (decreasing) in the interval $[k\pi,(k+1)\pi]$ with $k$ even (odd). Therefore, by the previous observations about the functions $g_i$, $i\in\{1,2\}$, we conclude that the same is true for $f^{\max}_{\overline{V}\ \overline{V}}$ and $f^{\min}_{\overline{V}\ \overline{V}}$. Moreover, we have
$$
f^{\max}_{\overline{V}\ \overline{V}}(x)=\lim\limits_{n\rightarrow \infty}f^{\max}_{V^{(n)}W^{(n)}}(x)\text{\ \ and\ \ }f^{\min}_{\overline{V}\ \overline{V}}(x)=\lim\limits_{n\rightarrow \infty}f^{\min}_{V^{(n)}W^{(n)}}(x),\ x\in\{k\pi,(k+1)\pi\},
$$
and the values $f^{\max}_{V^{(n)}W^{(n)}}(x)$ and $f^{\min}_{V^{(n)}W^{(n)}}(x)$ are given by the expressions provided for the extreme values of $f_{V^{(n)}W^{(n)}}^{\max}$ and $f_{V^{(n)}W^{(n)}}^{\min}$ in the proof of Theorem \ref{ThAB2} of Appendix \ref{AA} for the case $V\neq W$. This show that such expressions remain also valid for the case $V=W$.

\subsection{Lipschitz continuity of $H^{\max}_{VW}$ and $H^{\min}_{VW}$}

Theorem \ref{Contmax} is applied at the beginning of Subsection \ref{HmaxHminsect} of Appendix \ref{AB}, in both $\max$ and $\min$ versions,  to the function $g=H_{VW}:\mathbb{R}^2\rightarrow \mathbb{R}$ given in (\ref{Halphax}). Since the partial derivative
$
\frac{\partial H}{\partial \beta}
$
given in (\ref{dH}) is bounded in $\mathbb{R}^2$, the function $H_{VW}$ is Lipschitz continuous with respect to the variable $\beta$. We conclude that the functions $H_{VW}^{\max}$ and $H_{VW}^{\min}$ in  (\ref{HmaxHmindef}) are Lipschitz continuous.

\subsection{Case $V=W$ or NDR does not hold in the proof of Theorem \ref{ThAB2} of Appendix \ref{AB}} \label{E4}

Theorem \ref{Contmax} is applied, in both the $\max$ and $\min$ versions, to show that in Appendix \ref{AB} we can include in the proof of Theorem \ref{ThAB2} by a continuity argument the case $V=W$ or the Assumption NDR of Appendix \ref{AB} does not hold. Here is how this is done.

Let $(\overline{V},\overline{W})\in (0,1)^2$. We consider the function $g_1:\mathbb{R}^3\times B\rightarrow \mathbb{R}$, where $B$ is an open ball of center $(\overline{V},\overline{W})$ and radius $\varepsilon_0>0$ whose closure is strictly included in $(0,1)^2$, given by
\begin{eqnarray*}
	&&g_1(\alpha,x,\beta,V,W)=\frac{1+V\cos(x+\alpha)}{(1-W\cos\alpha)(1+V\cos(x+\beta))},\ (\alpha,x,\beta,V,W)\in \mathbb{R}^3\times B.
\end{eqnarray*}
Since the partial derivatives
\begin{eqnarray*}
	&&\frac{\partial g_1}{\partial x}(\alpha,x,\beta,V,W)=\frac{V\left(-\sin(x+\alpha)+\sin(x+\beta)-V\sin(\alpha-\beta)\right)}
	{(1-W\cos\alpha)\left(1+V\cos(x+\beta)\right)^2}\\
	&&
	\frac{\partial g_1}{\partial \beta}(\alpha,x,\beta,V,W)=\frac{V\left(1+V\cos(x+\alpha)\right)\sin(x+\beta)}
	{(1-W\cos\alpha)\left(1+V\cos(x+\beta)\right)^2}\\
	&&
	\frac{\partial g_1}{\partial V}(\alpha,x,\beta,V,W)
	=\frac{\cos(x+\alpha)-\cos(x+\beta)}
	{(1-W\cos\alpha)\left(1+V\cos(x+\beta)\right)^2}\\
	&&
	\frac{\partial g_1}{\partial W}(\alpha,x,\beta,V,W)
	=\frac{\left(1+V\cos(x+\alpha)\right)\cos\alpha}
	{\left(1+V\cos(x+\beta)\right)(1-W\cos\alpha)^2}\\
	&&(\alpha,x,\beta,V,W)\in \mathbb{R}^3\times B
\end{eqnarray*}
are bounded in $\mathbb{R}^3\times B$, $g_1$ is Lipschitz continuous with respect to the variables $(x,\beta,V,W)$ and then, by Theorem \ref{Contmax} applied to the function $g_1$, we conclude that the function $g_2:\mathbb{R}^2\times B\rightarrow \mathbb{R}$ given by
$$
g_2(x,\beta,V,W)=\max\limits_{\alpha\in\mathbb{R}}g_1(\alpha,x,\beta,V,W)=H_{VW}(x,\beta),\ (x,\beta,V,W)\in \mathbb{R}^2\times B, 
$$
is Lipschitz continuous. By Theorem \ref{Contmax} applied to the Lipschitz continuous function $g_2$, we conclude that the functions $g_3,g_4:\mathbb{R}\times B\rightarrow \mathbb{R}$ given by
\begin{eqnarray*}
	&&g_3(\beta,V,W)=\max\limits_{x\in\mathbb{R}}g_2(x,\beta,V,W)=H_{VW}^{\max}(\beta)\\
	&&g_4(\beta,V,W)=\min\limits_{x\in\mathbb{R}}g_2(x,\beta,V,W)=H_{VW}^{\min}(\beta)\\
	&&(\beta,V,W)\in \mathbb{R}\times B, 
\end{eqnarray*}
are Lipschitz continuous.

In Appendix \ref{AF}, it is proved the following fact: for any $\varepsilon>0$ with $\varepsilon<\varepsilon_0$, there exist finitely many $(V,W)\in B$ such that 
$$
\Vert (V,W)-(\overline{V},\overline{W})\Vert_2=\varepsilon,
$$
and $V=W$ or NDR does not hold.

Now, consider a sequence $\varepsilon_n$ such that $0<\varepsilon_n<\varepsilon_0$ for any $n$ and $\varepsilon_n\rightarrow 0$, $n\rightarrow \infty$, and a sequence $\left\{\left(V^{(n)},W^{(n)}\right)\right\}$ in $B$ such that, for any $n$, 
$$
\Vert (V^{(n)},W^{(n)})-(\overline{V},\overline{W})\Vert_2=\varepsilon_n,
$$
$V^{(n)}\neq W^{(n)}$ and NDR holds for $\left(V^{(n)},W^{(n)}\right)$.  Since $g_i$, $i\in\{3,4\}$, is a continuous function and
$$
\left(V^{(n)},W^{(n)}\right)\rightarrow \left(\overline{V},\overline{W}\right),\ n\rightarrow \infty,
$$
we have, for any $\beta\in\mathbb{R}$. 
$$
g_i\left(\beta,\overline{V},\overline{W}\right)=\lim\limits_{n\rightarrow \infty}g_i\left(\beta,V^{(n)},W^{(n)}\right)
$$
Therefore, for $\beta_1,\beta_2\in\mathbb{R}$ such that $\beta_1\leq \beta_2$, if
$$
g_i\left(\beta_1,V^{(n)},W^{(n)}\right)\leq g_i\left(\beta_2,V^{(n)},W^{(n)}\right)\text{\ \ for any $n$}
$$
then
$$
g_i\left(\beta_1,\overline{V},\overline{W}\right)\leq g_i\left(\beta_2,\overline{V},\overline{W}\right),
$$
and if
$$
g_i\left(\beta_1,V^{(n)},W^{(n)}\right)\geq g_i\left(\beta_2,V^{(n)},W^{(n)}\right)\text{\ \ for any $n$}
$$
then
$$
g_i\left(\beta_1,\overline{V},\overline{W}\right)\geq g_i\left(\beta_2,\overline{V},\overline{W}\right).
$$

In the proof of Theorem \ref{ThAB2} of Appendix \ref{AB}, it has been proved that, for any $n$, $H^{\max}_{V^{(n)}W^{(n)}}$ is monotonically increasing and  $H^{\min}_{V^{(n)}W^{(n)}}$ is monotonically decreasing in $[0,\pi]$. Therefore, by the previous observations about the functions $g_i$, $i\in\{3,4\}$, we conclude that the same is true for $H^{\max}_{\overline{V}\ \overline{W}}$ and $H^{\min}_{\overline{V}\ \overline{W}}$. Moreover, we have
$$
H_{VW}^{\max}\left(\beta\right)=\lim\limits_{n\rightarrow \infty}H_{V^{(n)}W^{(n)}}^{\max}\left(\beta\right)\text{\ \ and\ \ }H_{VW}^{\min}\left(\beta\right)=\lim\limits_{n\rightarrow \infty}H_{V^{(n)}W^{(n)}}^{\min}\left(\beta\right),\ \beta\in \{0,\pi\},
$$
and the values $H_{V^{(n)}W^{(n)}}^{\max}\left(\beta\right)$ and $H_{V^{(n)}W^{(n)}}^{\min}\left(\beta\right)$ are given by the expressions for the extreme values of $H_{V^{(n)}W^{(n)}}^{\max}$ and $H_{V^{(n)}W^{(n)}}^{\min}$ provided at the end of Appendix \ref{AB} for the case where $V\neq W$ and NDR holds. This show that such expressions remain also valid for the case $V=W$ or NDR does not hold.

\section{} \label{AF}

In this appendix, we show that the equation (\ref{d2Hdx2}) of Appendix \ref{AB} has finitely many solutions in $[0,\pi]$. Recall that we assume $V\neq 0$, $W\neq 0$ and $V\neq W$. Such equation is equivalent to
\begin{equation}
	\cos\left(x_i(\beta)+\beta\right)+V\cos\left(\alpha_{VW}^{\max}\left(x_i(\beta)\right)-\beta\right)=0 \label{eqform}
\end{equation}
(use (\ref{dH}) of Appendix \ref{AB} and (\ref{dfmax}) of Appendix \ref{AA}). The equivalent form (\ref{eqform}) can be written as
\begin{equation}
	\left\{
	\begin{array}{l}
		-V\sin(x+\alpha)-W\sin \alpha+VW\sin x=0\\
		-\sin(x+\alpha)+\sin(x+\beta)-V\sin(\alpha-\beta)=0\\
		\cos(x+\beta)+V\cos(\alpha-\beta)=0	
	\end{array}
	\right.
	\label{systemtwoA}
\end{equation}
where $x=x_i(\beta)$ and $\alpha=\alpha^{\max}_{VW}(x_i(\beta))$. The first equation is the equation for the stationary points $\alpha$ of $f_{VW}(\cdot,x)$ and the second equation is the equation for the stationary points of $H_{VW}(\cdot,\beta)$ (recall (\ref{eq}) in Appendix \ref{AB}): $x=x_i(\beta)$ is a stationary point of $H_{VW}(\cdot,\beta)$.

Set
$$
\gamma=x+\beta\text{\ \ and\ \ }\delta=\alpha-\beta.
$$
The system (\ref{systemtwoA}) becomes
\begin{equation}
	\left\{
	\begin{array}{l}
		-V\sin(\gamma+\delta)-W\sin(\gamma+\delta-x)+VW\sin x=0\\
	-\sin(\gamma+\delta)+\sin\gamma -V\sin \delta=0 \\
	\cos \gamma+V\cos \delta=0.
	\end{array}
	\right.
	\label{lasttwo}
\end{equation}
Expanding $\sin(\gamma+\delta)$ in the second equation of (\ref{lasttwo}), we obtain 
$$
(1-\cos \delta)\sin \gamma-(\cos \gamma+V)\sin \delta=0.
$$
Using the third equation in (\ref{lasttwo}), we obtain
$$
(1-\cos \delta)(\sin \gamma-V\sin \delta)=0.
$$
Thus,
$$
\cos \delta=1\text{\ \ or\ \ }\sin \gamma=V\sin \delta.
$$
Observe that $\sin \gamma=V\sin \delta$ and the third equation in (\ref{lasttwo}) imply
$$
1=\sin^2 \gamma+\cos^2 \gamma=V^2\sin^2 \delta+V^2\cos^2 \delta=V^2,
$$
which is impossible since $V\in(0,1)$. Therefore, any solution $(x,\gamma,\delta)$ of (\ref{lasttwo}) has $\cos \delta=1$, i.e., $\delta=2k\pi$ for some $k\in\mathbb{Z}$.

Therefore, the system (\ref{lasttwo}) reduces to
\begin{equation}
	\left\{
	\begin{array}{l}
		-V\sin\gamma-W\sin(\gamma-x)+VW\sin x=0\\
		\cos \gamma=-V\\
		\delta=2k\pi,\ k\in\mathbb{Z}.
	\end{array}
	\right.
	\label{lasttwo2}
\end{equation}
Expanding $\sin(\gamma-x)$ in the first equation of (\ref{lasttwo2}) and using the second equation, we obtain 
$$
	-(V+W\cos x)\sin \gamma=0.
$$
Since $\sin \gamma\neq0$ (otherwise $\cos\gamma=-V=-1$, and we have $V\in(0,1)$), it follows that
$$
	\cos x=-\frac{V}{W}.
$$
	
We can conclude what follows. If $V>W$, the system (\ref{lasttwo2}) has no solutions. Thus, if $V>W$, the equation (\ref{d2Hdx2}) of Appendix \ref{AB} has no solutions. If $V<W$, the solutions of (\ref{lasttwo2}) are
\begin{equation*}
	\left\{
	\begin{array}{l}
		x=\pm\arccos\!\left(-\frac{V}{W}\right)+2l\pi,\ l\in\mathbb{Z},\\
		\gamma=x+\beta=\pm\arccos(-V)+2m\pi,\ m\in\mathbb{Z}\\
		\delta=\alpha-\beta=2\pi k,\ k\in\mathbb{Z}
	\end{array}
	\right.
\end{equation*}
Thus, if $V<W$, the equation (\ref{d2Hdx2}) of Appendix \ref{AB} has finitely many solutions in $[0,\pi]$.

\section{} \label{A2H}

In this appendix we show that the system
\[
\frac{\partial H_{VW}}{\partial x}(x,\beta)=0
\text{\ \ and\ \ }
\frac{\partial^2 H_{VW}}{\partial x^2}(x,\beta)=0
\]
has finitely many solutions in
\((-\pi,\pi]\times[0,\pi]\). By (\ref{dH}) of Appendix \ref{AB}, this system is equivalent to
\begin{equation}
	\left\{
	\begin{array}{l}
		\left(f^{\max}_{VW}\right)^\prime(x)(1+V\cos(x+\beta))+f_{VW}^{\max}(x)V\sin(x+\beta)=0\\
		\left(f^{\max}_{VW}\right)^{\prime\prime}(x)(1+V\cos(x+\beta))+f_{VW}^{\max}(x)V\cos(x+\beta)=0.
	\end{array}
	\right.  \label{dHd2H=0}
\end{equation}

Let $(x,\beta)$ a solution of (\ref{dHd2H=0}). We cannot have $\left(f^{\max}_{VW}\right)^{\prime\prime}(x)+f_{VW}^{\max}(x)=0$, otherwise $f_{VW}^{\max}(x)=0$ by the second equation in (\ref{dHd2H=0}), which is impossible since $f^{\max}_{VW}$ has positive values (see, in Appendix \ref{AA}, Theorem \ref{maxminf}, or the definition (\ref{fVW}) of $f_{VW}$). Therefore,  by the second equation in (\ref{dHd2H=0}), we have
\begin{equation}
	\cos(x+\beta)=-\frac{\left(f^{\max}_{VW}\right)^{\prime\prime}(x)}{V\left(\left(f^{\max}_{VW}\right)^{\prime\prime}(x)+f_{VW}^{\max}(x)\right)}. \label{cos=}
\end{equation}
Then, by the first equation in (\ref{dHd2H=0}), we obtain
$$
f_{VW}^{\max}(x)^2V^2=\left(\left(f^{\max}_{VW}\right)^\prime(x)^2+\left(f^{\max}_{VW}\right)^{\prime\prime}(x)^2\right)\frac{f^{\max}_{VW}(x)^2}{\left(\left(f^{\max}_{VW}\right)^{\prime\prime}(x)+f_{VW}^{\max}(x)\right)^2}
$$
equivalent to
\begin{equation}
f_{VW}^{\max}(x)^2V^2\left(\left(f^{\max}_{VW}\right)^{\prime\prime}(x)+f_{VW}^{\max}(x)\right)^2-\left(\left(f^{\max}_{VW}\right)^\prime(x)^2+\left(f^{\max}_{VW}\right)^{\prime\prime}(x)^2\right)f^{\max}_{VW}(x)^2=0. \label{3f}
\end{equation}
The right-hand side in (\ref{3f}) is a real analytic function of $x$ (recall that $f_{VW}^{\max}$ is real analytic). Moreover it is not identically zero: set $x=0$ and use table (\ref{tablefmax}). Thus, (\ref{3f}) has finitely many solution in $x\in(-\pi,\pi]$. Now, by (\ref{cos=}), we conclude that there are finitely many solutions $(x,\beta)\in(-pi,pi]\times [0,\pi]$ of (\ref{dHd2H=0}).

\section{}  \label{AH}

In this appendix, we prove  the following fact given in Subsection \ref{E4} of Appendix \ref{ACMAXMIN}: for any $\varepsilon>0$ with $\varepsilon<\varepsilon_0$, there exist finitely many $(V,W)\in B$ such that 
$$
\Vert (V,W)-(\overline{V},\overline{W})\Vert_2=\varepsilon,
$$
and $V=W$ or NDR does not hold.

We look for $V,W\in(0,1)$ with $V\neq W$ such that the Assumption NDR of Appendix \ref{AB} does not hold, i.e., there exists $x_0\in\mathbb{R}$ such that
\begin{equation}
\frac{\partial H_{VW}}{\partial x}(x_0,0)=0\text{\ \ and\ \ }\frac{\partial^2 H_{VW}}{\partial x^2}(x_0,0)=0 \label{==0}.
\end{equation} 
The condition (\ref{==0}) can be written as
\begin{equation}
	\left\{
	\begin{array}{l}
		-V\sin(x_0+\alpha_0)-W\sin \alpha_0+VW\sin x_0=0\\
		-\sin(x_0+\alpha_0)+\sin x_0-V\sin\alpha_0=0\\
			p\left(\cos\left(x_0+\alpha_0\right),\cos x_0,\cos \alpha_0\right)=0,	
	\end{array}
	\right.
	\label{systemtwo}
\end{equation}
where $\alpha_0=\alpha^{\max}_{VW}(x_0)$. The first two equations correspond to $\frac{\partial H_{VW}}{\partial x}(x_0,0)=0$, once the equation for $\alpha_0$ is written explicitly. The third equation, where
$$
p(a,b,c)=V\left(-a+Wb\right)^2(1+Vb)-\left(1-Wc\right)\left(Va+Wc\right)(a-b)
$$
with
$$
a=\cos(x_0+\alpha_0),\ b=\cos x_0\text{\ \ and\ \ }c=\cos \alpha_0,
$$
correspond to $\frac{\partial^2 H_{VW}}{\partial x^2}(x_0,0)=0$.

In fact, by (\ref{dH}) we have
$$\frac{\partial^2 H_{VW}}{\partial x^2}(x_0,0)\notag \\
=\frac{\left(f^{\max}_{VW}\right)^{\prime\prime}(x_0)(1+V\cos x_0)+f_{VW}^{\max}(x_0)V\cos x_0}{(1+V\cos x_0)^2}
$$
Now, by (\ref{dfmax}) , the numerator of this fraction is
\begin{eqnarray*}
&&\left(f^{\max}_{VW}\right)^{\prime\prime}(x_0)(1+V\cos x_0)+f_{VW}^{\max}(x_0)V\cos x_0 \\
&& =\left(\frac{\left(-V\cos\left(x_0+\alpha_0\right)+VW\cos x_0\right)^2}{\left(1-W\cos \alpha_0\right)^2\left(V\cos\left(x_0+\alpha_0\right)+W\cos\alpha_0\right)}-\frac{V\cos\left(x_0+\alpha_0\right)}{1-W\cos \alpha_0}\right)(1+V\cos x_0)\\
&&\quad\quad+\frac{1+V\cos(x_0+\alpha_0)}{1-W\cos\alpha_0}V\cos x_0\\
&& =\left(\frac{\left(-Va+VWb\right)^2}{\left(1-Wc\right)^2\left(Va+Wc\right)}-\frac{Va}{1-Wc}\right)(1+Vb)+\frac{1+Va}{1-Wc}Vb\\
&&=\frac{\left(-Va+VWb\right)^2(1+Vb)}{\left(1-Wc\right)^2\left(Va+Wc\right)}-\frac{V(a-b)}{1-Wc}\\
&&=\frac{Vp(a,b,c)}{\left(1-Wc\right)^2\left(Va+Wc\right)}.
\end{eqnarray*}

In Proposition \ref{propxalpha} of Appendix \ref{AC}, the solutions $(x_0,\alpha_0)$, with $\alpha_0=\alpha^{\max}_{VW}(x_0)$ of the first two equation in (\ref{systemtwo}) have been determined. Such solutions are
\begin{itemize}
	\item $x_0$ and $\alpha_0$ even multiples of $\pi$;
	\item $x_0$ odd multiple of $\pi$ and $\alpha_0$ even multiple of $\pi$ when $V<W$;
	\item $x_0$ odd multiple of $\pi$ and $\alpha_0$ odd multiple of $\pi$ when $V>W$;
	\item the solutions of
	\begin{equation*}
		\cos x_0=-Q\text{\ \ and\ \ }\cos\alpha_0=\frac{Q-K}{L} 
	\end{equation*}
	with $\mathrm{sign}(\sin x_0)=-\mathrm{sign}(\sin \alpha_0)$, when $Q\leq 1$. Here
	$$
	Q=\frac{V(1+W)}{2W},\ L=\frac{V^2-W}{V(1-W)}\text{\ \ and\ \ }K=V-L=\frac{(1-V^2)W}{V(1-W)}.
	$$
\end{itemize}
By recalling Remark \ref{Q=1} of Appendix \ref{AC}, the fourth case with $Q=1$ can be included in the third case.

For $x_0$ and $\alpha_0$ even multiples of $\pi$, we have $a=b=1$, and therefore
$$
p(a,b,c)=V(1-W)^2(1+V)\neq 0.
$$

For $x_0$ odd multiple of $\pi$ and $\alpha_0$ even multiple of $\pi$, when $V<W$, we have $a=b=-1$, and therefore
\[
p(a,b,c)=V(1-W)^2(1-V)\neq 0.
\]

For $x_0$ odd multiple of $\pi$ and $\alpha_0$ an odd multiple of $\pi$ when $V>W$, we have $a=1$, $b=-1$ and $c=-1$, and therefore
\begin{eqnarray*}
	p(a,b,c)&=&	V(1+W)^2(1-V)-2(1+W)(V-W) \\
	&=&	(1+V)(1+W)\left(-V(1+W)+2W\right).
\end{eqnarray*}
Hence
$$
p(a,b,c)=0\Longleftrightarrow V=\frac{2W}{1+W}\Longleftrightarrow Q=1.
$$

It remains to consider the fourth case with $Q<1$, where 
$$
a=\frac{QK-1}{L},\ b=-Q\text{\ \ and\ \ }c=\frac{Q-K}{L}
$$
(for the first one, see (\ref{cosxplusalpha}) of Appendix \ref{AC}). A lengthy but straightforward computation gives
\begin{eqnarray*}
&&p(a,b,c)=\\
&&\frac{V(1-V^2)(1-W)^2((1+W)V-2W)((1+W)V+2W)(-(1+W)V^2+2W)}{8W(V^2-W)^2}.
\end{eqnarray*}
Since $V,W\in(0,1)$, we have
$$
p(a,b,c)=0\ \Longleftrightarrow\  (1+W)V-2W=0\text{\ \ or\ \ }-(1+W)V^2+2W=0.
$$
Since $Q<1$, we have $(1+W)V-2W<0$. Moreover, since $Q<1$ and $V\in(0,1)$, we also have
\[
V^2<V<\frac{2W}{1+W},
\]
and hence $-(1+W)V^2+2W>0$.

We have proved that there exists $x_0\in\mathbb{R}$ such that (\ref{==0}) holds if and only if $Q=1$, i.e., $V=\frac{2W}{1+W}$.

We conclude that for any $\varepsilon>0$ with $\varepsilon<\varepsilon_0$, the points $(V,W)$ in the circle  
$$
\Vert (V,W)-(\overline{V},\overline{W})\Vert_2=\varepsilon,
$$
such that $V=W$ or NDR does not hold are points of intersection with the line $V=W$ or the curve $V=\frac{2W}{1+W}$ and then they are finitely many.

\section{} \label{AC}

In this appendix, we determine the extreme values of the function $H_{VW}(\cdot,0)$, defined by setting $\beta=0$ in (\ref{Halphax}) of Appendix \ref{AB}. We assume $V\neq 0$, $W\neq 0$ and $V\neq W$.

The stationary points of $H_{VW}(\cdot,0)$ are the solutions $x$ of the equation
\begin{equation*}
	-\sin(x+\alpha^{\max}_{VW}(x))+\sin x-V\sin \alpha^{\max}_{VW}(x)=0 
\end{equation*}
(see (\ref{eq}) in Appendix \ref{AB} with $\beta=0$).

Since $ \alpha^{\max}_{VW}(x)$ is a stationary point of the function $f_{VW}(\cdot,x)$, we solve the two equations
\begin{equation*}
	\left\{
	\begin{array}{l}
		-V\sin(x+\alpha)-W\sin \alpha+VW\sin x=0\\
		-\sin(x+\alpha)+\sin x-V\sin\alpha=0
	\end{array}
	\right.
\end{equation*}
in the unknowns $x$ and $\alpha$, where the first equation is the equation for the stationary points $\alpha$ of $f_{VW}(\cdot,x)$ (see (\ref{df}) in Appendix \ref{AA}).

By performing some algebraic manipulations, the previous system of equations is transformed into the equivalent system
\begin{equation}
	\left\{
	\begin{array}{l}\sin x=L\sin \alpha\\
		\sin(x+\alpha)=-K\sin\alpha,
	\end{array}
	\right.  \label{sys}
\end{equation}
where
\begin{equation}
	L:=\frac{V^2-W}{V(1-W)}\text{\  \ and\ \ }K:=V-L=\frac{(1-V^2)W}{V(1-W)}. \label{LK}
\end{equation}
Observe that $L$ can be negative, zero or positive, whereas $K$ is positive.

Regarding the solutions of (\ref{sys}), we have the following result.

\begin{lemma} \label{Lemmasolutions}
	Let
	\begin{equation}
		Q:=\frac{K^2-L^2+1}{2K}=\frac{V(1+W)}{2W}. \label{QQ}
	\end{equation}
	If $Q>1$, then the set of solutions of (\ref{sys}) is given by the pairs $(x,\alpha)$ such that $x$ and $\alpha$ are multiples of $\pi$. If $Q\leq 1$, then the set of solutions also includes the pairs $(x, \alpha)$ such that
	$$
	\cos x=-Q\text{\ \ and\ \ }\cos\alpha=\frac{Q-K}{L}
	$$
	with $\mathrm{sign}(\sin x)=-\mathrm{sign}(\sin \alpha)$.
\end{lemma}

\begin{proof}
	
	To prove the second equality in (\ref{QQ}), observe that
	$$
	K^2-L^2+1=V^2-2VL+1=\frac{(1-V^2)(1+W)}{1-W}.
	$$
	
	Clearly, the pairs $(x,\alpha)$ such that $x$ and $\alpha$ are multiples of $\pi$ are solutions of \eqref{sys}. We now investigate whether any other solutions exist.
	
	Consider the case $L=0$ and then $K=V$. The first equation in (\ref{sys}) is $\sin x=0$. If $\sin x=0$, then $\sin(x+\alpha)=\pm \sin \alpha$ in the second equation, which becomes $(\pm 1+V)\sin \alpha=0$, i.e. $\sin \alpha=0$. Therefore, in the case $L=0$, the solutions of (\ref{sys}) are only the pairs $(x,\alpha)$ such that $x$ and $\alpha$ are multiples of $\pi$. 
	
	Consider the case $L\neq 0$. The solutions of  (\ref{sys}) with $\sin x=0$ have also $\sin \alpha=0$ and then they are the pairs $(x,\alpha)$ such that $x$ and $\alpha$ are multiples of $\pi$.
	
	Now, in the case $L\neq 0$, we look for solutions of (\ref{sys}) with $\sin x\neq 0$. Such solutions have $\sin \alpha\neq 0$. By expanding $\sin(x+\alpha)$ as
	$$
	\sin(x+\alpha)=\sin x\cos\alpha+\cos x\sin\alpha,
	$$
	the second equation in (\ref{sys}) provides
	\begin{equation}
		\sin x \cos\alpha=-(\cos x+K)\sin \alpha. \label{xplusalpha}
	\end{equation}
	By squaring and by using the first equation, we obtain
	$$
	\sin^2 x\left(1-\frac{\sin^2 x}{L^2}\right)=(\cos x+K)^2 \frac{\sin^2 x}{L^2}
	$$
	and then, since $\sin x\neq 0$,
	\begin{equation*}
		1-\frac{\sin^2 x}{L^2}=\frac{(\cos x+K)^2}{L^2} \label{1-}
	\end{equation*}
	equivalent to
	\begin{equation}
		\cos x=\frac{L^2-K^2-1}{2K}=-Q. \label{Q1A}
	\end{equation}
	Concerning the other unknown $\alpha$ in case of $\cos x=-Q$, by using  (\ref{xplusalpha}) and the first equation in (\ref{sys}), we obtain
	\begin{equation*}
		L\sin\alpha \cos\alpha=(Q-K)\sin\alpha. \label{alphaA}
	\end{equation*}
	Since $\sin \alpha\neq 0$, we have
	\begin{equation}
		\cos \alpha=\frac{Q-K}{L}. \label{alphaB}
	\end{equation}
	We have proved that, in the case $L\neq 0$, solutions $(x,\alpha)$ of (\ref{sys}) with $\sin x\neq 0$ satisfy (\ref{Q1A}) and (\ref{alphaB}).

	Therefore, we conclude that if $Q>1$, then the system (\ref{sys}) has only the solutions with $\sin x=0$, i.e., the solutions are only the pair $(x,\alpha)$ such that $x$ and $\alpha$ are multiples of $\pi$. This is true for both cases $L=0$ and $L\neq 0$.
	
	On the other hand, if $Q\leq 1$, and hence $L<0$ (see point 2 of Remark \ref{remApp} below), in addition to the solutions $(x,\alpha)$ such that $x$ and $\alpha$ are multiples of $\pi$, we have the additional solutions $(x,\alpha)$ satisfying (\ref{Q1A}) and (\ref{alphaB})  with $\mathrm{sign}(\sin x)=-\mathrm{sign}(\sin \alpha)$. Observe that the equation (\ref{alphaB}) has solutions: see point 3 of Remark \ref{remApp} below.  
	
	In fact, (\ref{Q1A}) and (\ref{alphaB}) imply
	\begin{equation}
		L^2-\sin^2x=(\cos x+K)^2. \label{L2sin2}
	\end{equation}
	and
	\begin{equation}
		L^2\cos^2\alpha=(\cos x+K)^2  \label{L2cos2}
	\end{equation}
	By comparing (\ref{L2cos2}) and (\ref{L2sin2}), we obtain
	$$
	\sin^2x=L^2\sin^2\alpha,
	$$
	i.e.
	$$
	\sin x=L\sin \alpha\text{\ \ or\ \ }\sin x=-L\sin \alpha.
	$$
	Therefore, by taking $(x,\alpha)$ satisfying (\ref{Q1A}) and (\ref{alphaB})  with $\mathrm{sign}(\sin x)=-\mathrm{sign}(\sin \alpha)$, we obtain the first equation 
	$
	\sin x=L\sin \alpha
	$
	in (\ref{sys}), by recalling that $L<0$. Moreover, (\ref{Q1A}), (\ref{alphaB}) and $\sin x=L\sin \alpha$ imply (\ref{xplusalpha}), equivalent to the second equation in (\ref{sys}): in fact, by (\ref{alphaB}) and (\ref{Q1A}) we obtain
	$$
	L\sin\alpha\cos\alpha=(Q-K)\sin\alpha=-(\cos x+K)\sin\alpha,
	$$
	and by $\sin x=L\sin \alpha$ we obtain (\ref{xplusalpha}).

	There are no other solutions of (\ref{sys}) in the case $Q\leq 1$, since, as we have proved above, each solution $(x,\alpha)$ with $\sin x\neq 0$ satisfies (\ref{Q1A}), (\ref{alphaB})  and $\mathrm{sign}(\sin x)=-\mathrm{sign}(\sin \alpha)$ (recall that the first equation is $\sin x=L\sin \alpha$ with $L<0$).

\end{proof}

\begin{remark} \label{remApp}
	\quad
	\begin{itemize}
		\item [1.]
		$$
		Q\leq 1 \Longleftrightarrow V\leq\frac{2W}{1+W} \Longleftrightarrow W\geq \frac{V}{2-V},
		$$
		with
		$$
		\frac{2W}{1+W}<1\text{\ \ and\ \ }\frac{V}{2-V}<1;
		$$
		\item  [2.] If $Q\leq 1$, then $L<0$. To prove this, use the upper bound of $V$ from point 1 in the numerator $V^2-W$ of the fraction defining $L$.
		\item [3.] If $Q\leq 1$, then
		$$
		\left(\frac{Q-K}{L}\right)^2\leq 1.
		$$
		To prove this, observe that
		\begin{eqnarray*}
			(Q-K)^2&=&Q^2-2QK+K^2=Q^2-(K^2-L^2+1)+K^2\\
			&=&Q^2+L^2-1.
		\end{eqnarray*}
		\item [4.] 
		$$
		Q-K\geq 0 \Longleftrightarrow V^2\geq \frac{2W^2}{1+W^2},
		$$
		with
		$$
		\frac{2W^2}{1+W^2}<1.
		$$
		To prove the equivalence $\Longleftrightarrow$, use
		$$
		Q=\frac{V(1+W)}{2W}\text{\ \ and\ \ }K=\frac{(1-V^2)W}{V(1-W)}.
		$$
	\end{itemize}
\end{remark}

Starting from the parameters $V$ and $W$, we have defined the quantities $L$ and $K$ in (\ref{LK}) and $Q$ in (\ref{QQ}). The next lemma shows how $K$, $L$ and $W$ can be expressed in terms of $V$ and $Q$.
\begin{lemma} \label{LemmaLKW}
	We have
	\begin{equation*}
		L=-\frac{V^2-2QV+1}{2(Q-V)},\ K=\frac{1-V^2}{2(Q-V)}\text{\ \ and\ \ }W=\frac{V}{2Q-V}
		\label{formula}
	\end{equation*}
	with
	\begin{equation*}
		Q-V=\frac{V(1-W)}{2W}>0.
	\end{equation*}
\end{lemma}
\begin{proof}
	By using the definitions
	\begin{equation*}
		Q=\frac{K^2-L^2+1}{2K}\text{\ \ and\ \ }K=V-L,
	\end{equation*}
	we obtain
	$$
	Q=\frac{V^2-2VL+1}{2(V-L)}
	$$
	and then
	$$
	2(V-Q)L=V^2-2QV+1.
	$$
	The expression for $L$ follows. Then, the expression for $K$ follows:
	$$
	K=V-L=V+\frac{V^2-2QV+1}{2(Q-V)}=\frac{1-V^2}{2(Q-V)}.
	$$
	The expression for $W$ is obtained by (\ref{QQ}). Finally, we have
	$$
	Q-V=\frac{V(1+W)}{2W}-V=\frac{V(1-W)}{2W}>0.
	$$
\end{proof}

The solutions of (\ref{sys}) are given in Lemma \ref{Lemmasolutions}. It is of interest to determine which solutions $(x,\alpha)$ of (\ref{sys}) correspond to $\alpha=\alpha^{\max}_{VW}(x)$, i.e., $\alpha$ is a maximum point of $f_{VW}(\cdot,x)$. The first components $x$ of such solutions $(x,\alpha)$ are the stationary points of $H_{VW}(\cdot,0)$. 

\begin{proposition}\label{propxalpha}
	The solutions $(x,\alpha)$ of (\ref{sys}) such that $\alpha$ is a maximum point of $f_{VW}(\cdot,x)$ are:
	\begin{itemize}
		\item $x$ and $\alpha$ even multiples of $\pi$;
		\item $x$ odd multiple of $\pi$ and $\alpha$ even multiple of $\pi$ when $V<W$;
		\item $x$ odd multiple of $\pi$ and $\alpha$ odd multiple of $\pi$ when $V>W$;
		\item the solutions of
		\begin{equation}
			\cos x=-Q\text{\ \ and\ \ }\cos\alpha=\frac{Q-K}{L} \label{coscos}
		\end{equation}
		with $\mathrm{sign}(\sin x)=-\mathrm{sign}(\sin \alpha)$, when $Q\leq 1$. 
	\end{itemize}
\end{proposition}

\begin{proof}
	The condition under which $\alpha$ is a maximum point of $f_{VW}(\cdot,x)$ is
	\begin{equation}
		V\cos(x+\alpha)+W\cos \alpha>0.  \label{V1W1geq0}
	\end{equation}
	This condition is obtained by imposing a negative value to the second derivative $\frac{\partial^2 f_{VW}}{\partial \alpha^2}(\alpha,x)$, given in (\ref{df}) of Appendix \ref{AA} for $\alpha$ stationary point of $f_{VW}(\cdot,x)$. Indeed, \eqref{V1W1geq0} is the condition for local maximum point, but Theorem \ref{LemmafvW} of Appendix \ref{AA} shows that the local maximum point are global maximum point.
	
	For a solution $(x,\alpha)$ of (\ref{sys}) such that $x$ and $\alpha$ are multiples of $\pi$, the values of the left-hand side in (\ref{V1W1geq0}) are given in the next table
	\begin{equation*}
		\begin{tabular}{|c|c|c|}
			\hline
			& $x$ even  multiple of $\pi$             & $x$ odd multiple of $\pi$\\
			\hline
			$\alpha$ even multiple of $\pi$ & $V+W$  &  $-V+W$\\
			\hline
			$\alpha$ odd  multiple of $\pi$ & $-V-W$ & $V-W$ \\
			\hline
		\end{tabular} 
	\end{equation*}
	Positive values, and therefore maximum points $\alpha$, are obtained in the following cases:
	\begin{itemize}
		\item [a)] $x$ and $\alpha$ even multiples of $\pi$;
		\item [b)] $x$ odd multiple of $\pi$ and $\alpha$ even multiple of $\pi$, when $V<W$;
		\item [c)] $x$ odd multiple of $\pi$ and $\alpha$ odd multiple of $\pi$, when $V>W$.
	\end{itemize}
	
	Now, we consider the additional solutions (\ref{coscos}) in the case of $Q\leq 1$. We expand  $\cos(x+\alpha)$ in (\ref{V1W1geq0}) as
	\begin{eqnarray}
		\cos(x+\alpha)&=&\cos x\cos \alpha-\sin x\sin \alpha  \notag\\
		&=&-Q\cos\alpha-L\sin^2\alpha\notag \\
		&&\ \ \ \text{($\cos x=-Q$ and $\sin x=L\sin\alpha$: see (\ref{coscos}) and (\ref{sys}))}\notag\\
		&=&L\cos^2\alpha-Q\cos\alpha-L \notag\\
		&=&L\left(\frac{Q-K}{L}\right)^2-Q\frac{Q-K}{L}-L \notag\\
			&&\ \ \ \text{($\cos \alpha=\frac{Q-K}{L}$: see (\ref{coscos}))}\notag\\
		&=&\frac{K^2-QK-L^2}{L}\notag \\
		&=&\frac{QK-1}{L}, \notag \\
		\label{cosxplusalpha}
	\end{eqnarray}
	where the last equality follows by
	$$
	L^2=K^2-2QK+1
	$$
	(see (\ref{QQ})). Since $L<0$ (recall point 2 in Remark \ref{remApp}), we obtain that the condition (\ref{V1W1geq0}) is equivalent to
	\begin{equation}
		V(QK-1)+W(Q-K)<0. \label{cV1W1}
	\end{equation}
	By using the expressions for $K$ and $W$ given in Lemma \ref{LemmaLKW}, we obtain
	$$
	V(QK-1)=\frac{V}{2(Q-V)}\cdot\left(-Q(1+V^2)+2V\right)
	$$
	and
	$$
	W(Q-K)=\frac{V}{2(2Q-V)(Q-V)}\cdot\left(2Q^2-2QV-1+V^2\right)
	$$
	and then
	\begin{eqnarray*}
		&&V(QK-1)+W(Q-K)\\
		&&=\frac{V}{2(2Q-V)(Q-V)}\\
		&&\ \ \ \ \cdot\left(-2V^2Q^2+(V^3+3V)Q-(1+V^2)\right).
	\end{eqnarray*}
	The trinomial 
	$$
	-2V^2Q^2+(V^3+3V)Q-(1+V^2),
	$$
	in $Q$, which has the same sign as the left-hand side of (\ref{cV1W1}) (recall Lemma \ref{LemmaLKW}: we have $Q-V>0$), is a downward parabola with vertex at
	$$
	-\frac{V^3+3V}{2(-2V^2)}=\frac{V^2+3}{4V}>1.
	$$
	Therefore, by confining to $Q\leq 1$, the maximum value of the trinomial is at $Q=1$ and such maximum value is
	$$
	(V-1)^3<0.
	$$
	We conclude that, for $Q\leq 1$, the trinomial has negative sign and then (\ref{cV1W1}) holds.
	
	Consequently, the solutions $(x,\alpha)$ given by (\ref{coscos}) have $\alpha$ as a maximum point of $f_{VW}(\cdot,x)$.
\end{proof}

\begin{remark}\label{Q=1}
	In the previous Proposition \ref{propxalpha}, the solutions of the fourth case with $Q=1$ are included in the third case $x$ odd multiple of $\pi$ and $\alpha$ odd multiple of $\pi$ when $V>W$. In fact, $Q=\frac{V(1+W)}{2W}=1$ implies $V>W$. Moreover (\ref{coscos}) implies $x$ odd multiple of $\pi$ and $\alpha$ odd multiple of $\pi$: for the latter, observe that if $Q=1$, then
	$$
	L=-\frac{1-V}{2}\text{\ \ and\ \ }K=\frac{1+V}{2}.
	$$	 
\end{remark}

The following proposition, which follows from the previous Proposition \ref{propxalpha}, specifies the stationary points of the function $H_{VW}(\cdot,0)$ and their corresponding values.
\begin{proposition}\label{lastprop}
	The stationary points of the function $H_{VW}(\cdot,0)$ are:
	\begin{itemize}
		\item The even multiples $x$ of $\pi$ with value
		\begin{equation}
			H_{VW}(x,0)=\frac{1}{1-W}. \label{DeltaH}
		\end{equation}
		\item The odd multiples $x$ of $\pi$ with value
		\begin{equation}
			H_{VW}(x,0)=\left\{
			\begin{array}{l}
				\frac{1}{1-W}	\text{\ if\ }V<W\\
				\\
				\frac{1+V}{(1-V)(1+W)}	\text{\ if\ }V>W.
			\end{array}
			\right.  \label{DeltaDeltaH}
		\end{equation}
		\item If $Q\leq 1$, the solutions $x$ of $\cos x=-Q$ with value
		\begin{equation}
			H_{VW}(x,0)=\frac{1-V^2}{(1-QV)(1-W)}. \label{Hx0}
		\end{equation}
	\end{itemize}
\end{proposition}
\begin{proof}
	We have that $x\in\mathbb{R}$ is a stationary point of $H_{VW}(\cdot,0)$ if and only if $x$ is the first component of a pair $(x,\alpha)$ in Proposition \ref{propxalpha}. Regarding the values at the stationary points, we utilize  
	\begin{equation}
		H_{VW}(x,0)=\frac{f_{VW}^{\max}(x)}{1+V\cos x}=\frac{1+V\cos(x+\alpha)}{(1+V\cos x)(1-W\cos \alpha)}  \label{H0}
	\end{equation}
	with $\alpha=\alpha_{VW}^{\max}(x)$ (recall (\ref{Halphax}) of Appendix \ref{AB} and (\ref{dfmax}) of Appendix \ref{AA}).
	
	For the cases $x$ even multiple of $\pi$ and $x$ odd multiple of $\pi$, the computations in (\ref{H0}) are straightforward and we obtain (\ref{DeltaH}) and (\ref{DeltaDeltaH}).
	
	For the case $Q\leq 1$, use in (\ref{H0}):
	$$
	\cos x=-Q,\ \cos\alpha=\frac{Q-K}{L}\text{\ \ and\ \ }\cos(x+\alpha)=\frac{QK-1}{L}
	$$
	(refer to (\ref{coscos}) and (\ref{cosxplusalpha})). We obtain
	\begin{eqnarray*}
		H_{VW}(x,0)&=&\frac{L+VQK-V}{(1-VQ)(L-W(Q-K))}\\
		&=&\frac{VQK-K}{(1-VQ)(L-W(Q-K))}\text{\ \ (use $L-V=-K$ in the numerator: see (\ref{LK}))}\\
		&=&\frac{-K}{L-W(Q-K)}\\
		&=&\frac{-K}{L-W(Q-V+L)}\text{\ \ (use $-K=-V+L$ in the denominator)}\\
		&=&\frac{-K}{(1-W)L-W(Q-V)}\\
		&=&\frac{K}{\left(-L+\frac{W}{1-W}(Q-V)\right)(1-W)}.
	\end{eqnarray*}
	By using the expression for $W$ in Lemma \ref{LemmaLKW}, we have
	$$
	\frac{W}{1-W}=\frac{V}{2(Q-V)}
	$$
	and then
	\begin{eqnarray*}
		H_{VW}(x,0)&=&\frac{K}{\left(-L+\frac{V}{2}\right)(1-W)}\\
		&=&\frac{K}{\left(K-\frac{V}{2}\right)(1-W)}\text{\ \ (use $-L=K-V$ in the denominator)}.
	\end{eqnarray*}
	By using the expression for $K$ in Lemma \ref{LemmaLKW}, we obtain (\ref{Hx0}).
	
\end{proof}

Here is our final result concerning the values $H_{VW}^{\max}(0)$ and $H_{VW}^{\min}(0)$. It follows from the previous proposition regarding the stationary points of $H_{VW}(\cdot, 0)$. 
\begin{theorem}	\label{lastThAC}
	We have
	\begin{equation*}
		H_{VW}^{\max}(0)=\left\{
		\begin{array}{l}
			\frac{1-V^2}{(1-QV)(1-W)}\text{\ if\ }Q\leq1\\
			\\
			\frac{1+V}{(1-V)(1+W)}\text{\ if\ }Q>1,
		\end{array}
		\right.
	\end{equation*}
	where $Q$ is defined in (\ref{QQ}), and
	\begin{equation*}
		H_{VW}^{\min}(0)=\frac{1}{1-W}. 
	\end{equation*}
	Moreover, for $x_0\in\mathbb{R}$, we have $H_{VW}^{\max}(0)=H_{VW}(x_0,0)$ if and only if
	$$
	\left\{
	\begin{array}{l}
		\cos x_0=-Q\text{\ if\ }Q\leq 1\\
		\\
		x_0\text{\ is an odd multiple of $\pi$\ if\ }Q>1,
	\end{array}
	\right.
	$$
	and we have $H_{VW}^{\min}(0)=H_{VW}(x_0,0)$ if and only if
	$$
	\left\{
	\begin{array}{l}
		x_0\text{\ is a multiple of $\pi$\ if\ }V<W\\
		\\
		x_0\text{\ is an even multiple of $\pi$\ if\ }V>W.
	\end{array}
	\right.
	$$
\end{theorem}

\begin{proof}
	We distinguish the three cases:
	\begin{equation}
		a):V<W,\ b):W< V\leq\frac{2W}{1+W}\text{\ \ and\ \ }c):\frac{2W}{1+W}< V. \label{abc}
	\end{equation}
	
	In point 1 of Remark \ref{remApp}, we observed that $Q\leq 1$ is equivalent to
	$$
	V\leq\frac{2W}{1+W}.
	$$
	Therefore, in the cases a) and b), we have $Q\leq 1$ and, in the case c), we have $Q>1$. Moreover, the case b) is well-defined, since 
	$$
	\frac{2W}{1+W}>W
	$$
	holds.
	
	In all three cases a), b) and c), $H^{\max}_{VW}(0)$ and $H^{\min}_{VW}(0)$ are determined by comparing the values $H_{VW}(x,0)$ when $x$ is a stationary point of $H_{VW}(\cdot,0)$, since $H^{\max}_{VW}(0)$ and $H^{\min}_{VW}(0)$ are values at stationary points.

	Case a). For $V<W$, we have the following stationary points of $H_{VW}(\cdot,0)$ (see Proposition \ref{lastprop}). 
	\begin{itemize}
		\item The multiples $x$ of $\pi$ with value
		$$
		H_{VW}(x,0)=\frac{1}{1-W}.
		$$
		\item  Since $Q\leq 1$, the solutions $x$ of $\cos x=-Q$ with value
		$$
		H_{VW}(x,0)=\frac{1-V^2}{(1-QV)(1-W)}>\frac{1}{1-W},
		$$
		where the inequality $>$ holds by $Q>V$ (see Lemma \ref{LemmaLKW}).
	\end{itemize}
	We conclude that
	$$
	H_{VW}^{\max}(0)=\frac{1-V^2}{(1-QV)(1-W)} \text{\ \ and\ \ }H_{VW}^{\min}(0)=\frac{1}{1-W}.
	$$
	Moreover, for $x_0\in\mathbb{R}$, we have $H_{VW}^{\max}(0)=H_{VW}(x_0,0)$ if and only if $\cos x_0=-Q$ and we have $H_{VW}^{\min}(0)=H_{VW}(x_0,0)$ if and only if  
	$x_0$ is a multiple of $\pi$.
	
	Case b). For $W<V\leq\frac{2W}{1+W}$, we have the following stationary points of $H_{VW}(\cdot,0)$ (see Proposition \ref{lastprop}).
	\begin{itemize}
		\item The even multiples $x$ of $\pi$ with value
		$$
		H_{VW}(x,0)=\frac{1}{1-W}.
		$$
		\item The odd multiples $x$ of $\pi$ with value
		$$
		H_{VW}(x,0)=\frac{1+V}{(1-V)(1+W)}> \frac{1}{1-W},
		$$
		where the inequality $>$ follows by $V>W$.
		\item Since $Q\leq 1$, the solutions $x$ of $\cos x=-Q$ with value
		$$
		H_{VW}(x,0)=\frac{1-V^2}{(1-QV)(1-W)}\geq \frac{1+V}{(1-V)(1+W)},
		$$
		where the inequality $\geq$ is now proved. We rewrite the inequality as
		\begin{equation}
			\frac{(1-V)^2}{1-QV}\geq \frac{1-W}{1+W}.\label{ineqrew}
		\end{equation}
		Using the expression for $W$ in Lemma \ref{LemmaLKW}, we obtain
		$$
		\frac{1-W}{1+W}=\frac{Q-V}{Q}
		$$
		and then (\ref{ineqrew}) is equivalent to
		$$
		(1-V)^2Q\geq(1-QV)(Q-V),
		$$
		that reduces to
		$$
		Q^2-2Q+1=(Q-1)^2\geq 0.
		$$
		Observe that the inequality is strict if and only if $Q<1$. When $Q=1$, equality holds, and the condition $\cos x=-Q$ is equivalent to $x$ being an odd multiple of $\pi$.
	\end{itemize}
	We conclude that
	$$
	H_{VW}^{\max}(0)=\frac{1-V^2}{(1-QV)(1-W)}\text{\ \ and\ \ }H_{VW}^{\min}(0)=\frac{1}{1-W}.
	$$
	Moreover, for $x_0\in\mathbb{R}$, we have $H_{VW}^{\max}(0)=H_{VW}(x_0,0)$ if and only if $\cos x_0=-Q$, and we have $H_{VW}^{\min}(0)=H_{VW}(x_0,0)$ if and only if $x_0$ is an even multiple of $\pi$.

	Case c). For $\frac{2W}{1+W}<V$, since $Q>1$ and $V>W$, we have the following stationary points of $H(\cdot,0)$ (see Proposition \ref{lastprop}).
	\begin{itemize}
		\item 
		The even multiples $x$ of $\pi$ with value
		$$
		H_{VW}(x,0)=\frac{1}{1-W}.
		$$
		\item The odd multiples $x$ of $\pi$ with value
		$$
		H_{VW}(x,0)=\frac{1+V}{(1-V)(1+W)}>\frac{1}{1-W},
		$$
		where the inequality follows by $V>W$.
	\end{itemize}
	We conclude that
	$$
	H_{VW}^{\max}(0)=\frac{1+V}{(1-V)(1+W)}\text{\ \ and\ \ }H_{VW}^{\min}(0)=\frac{1}{1-W}.
	$$
    Moreover, for $x_0\in\mathbb{R}$, we have  $H_{VW}^{\max}(0)=H(x_0,0)$ if and only if $x_0$ is an odd multiple of $\pi$, and we have $H^{\min}(0)=H(x_0,0)$ if and only if $x_0$ is even multiple of $\pi$.
	
	Observe that in the transition case between b) and c), i.e., when $V=\frac{2W}{1+W}$ and $Q=1$, formulas in b) and c) for $H_{VW}^{\max}(0)$ give the same value.
\end{proof}

\begin{figure}
	\centering
	\subfloat[Case a): $V=0.45$ and $W=0.5$]{\includegraphics[width=1\textwidth]{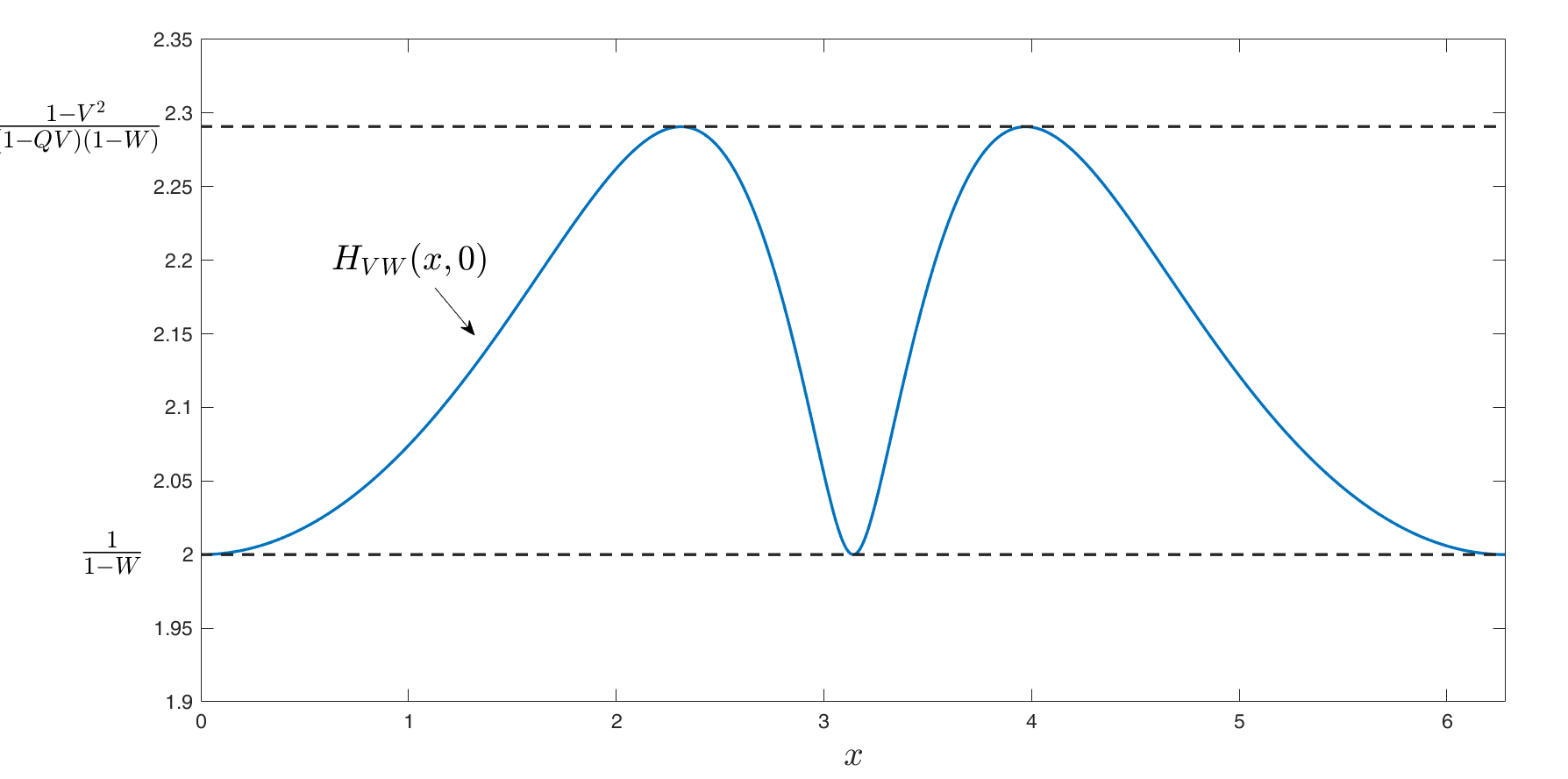}}

	\vspace{0.5cm} 
	
	\subfloat[Case b): $V=0.55$ and $W=0.5$]{\includegraphics[width=1\textwidth]{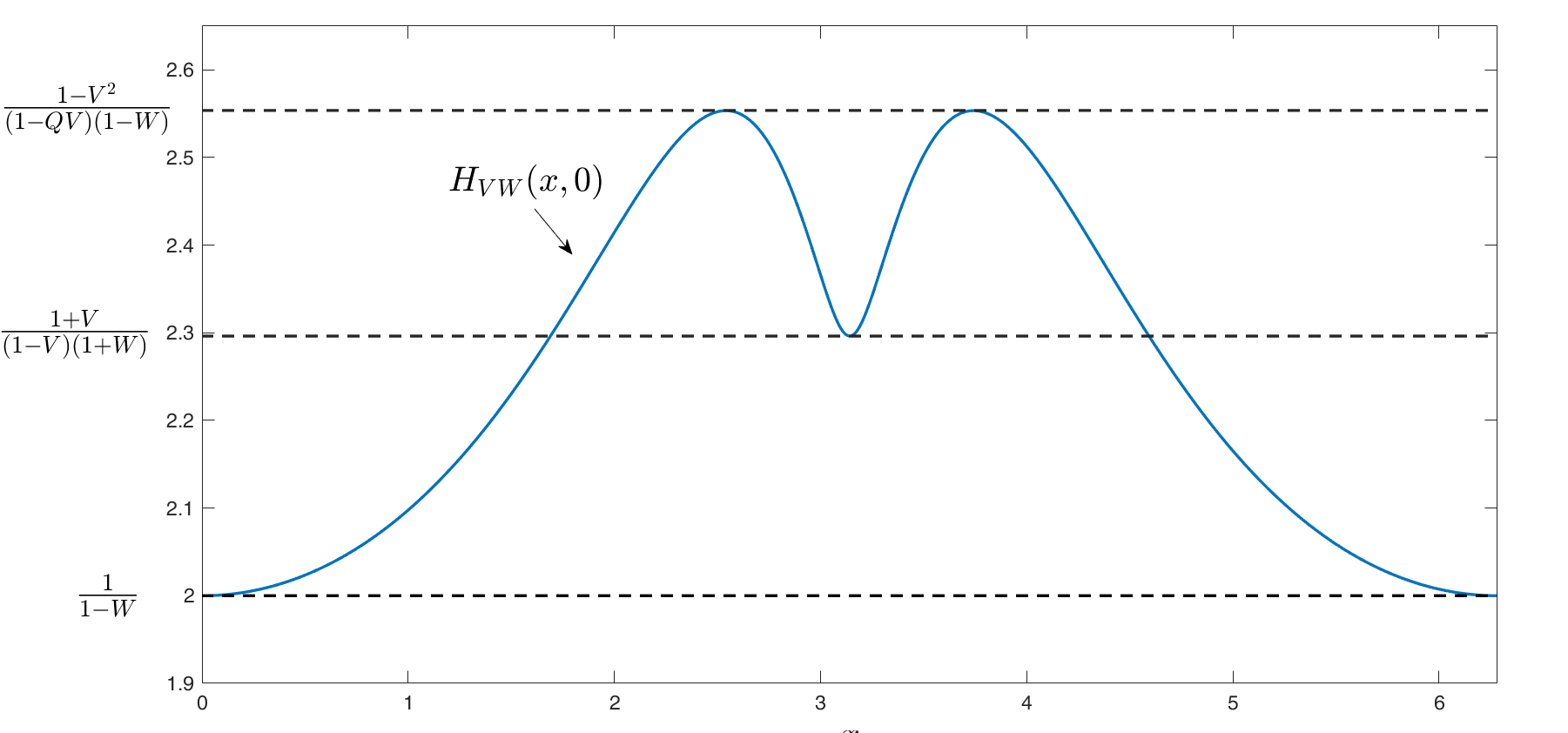}}
	
	\vspace{0.5cm} 
	
	\subfloat[Case c): $V=0.7$ and $W=0.5$]{\includegraphics[width=1\textwidth]{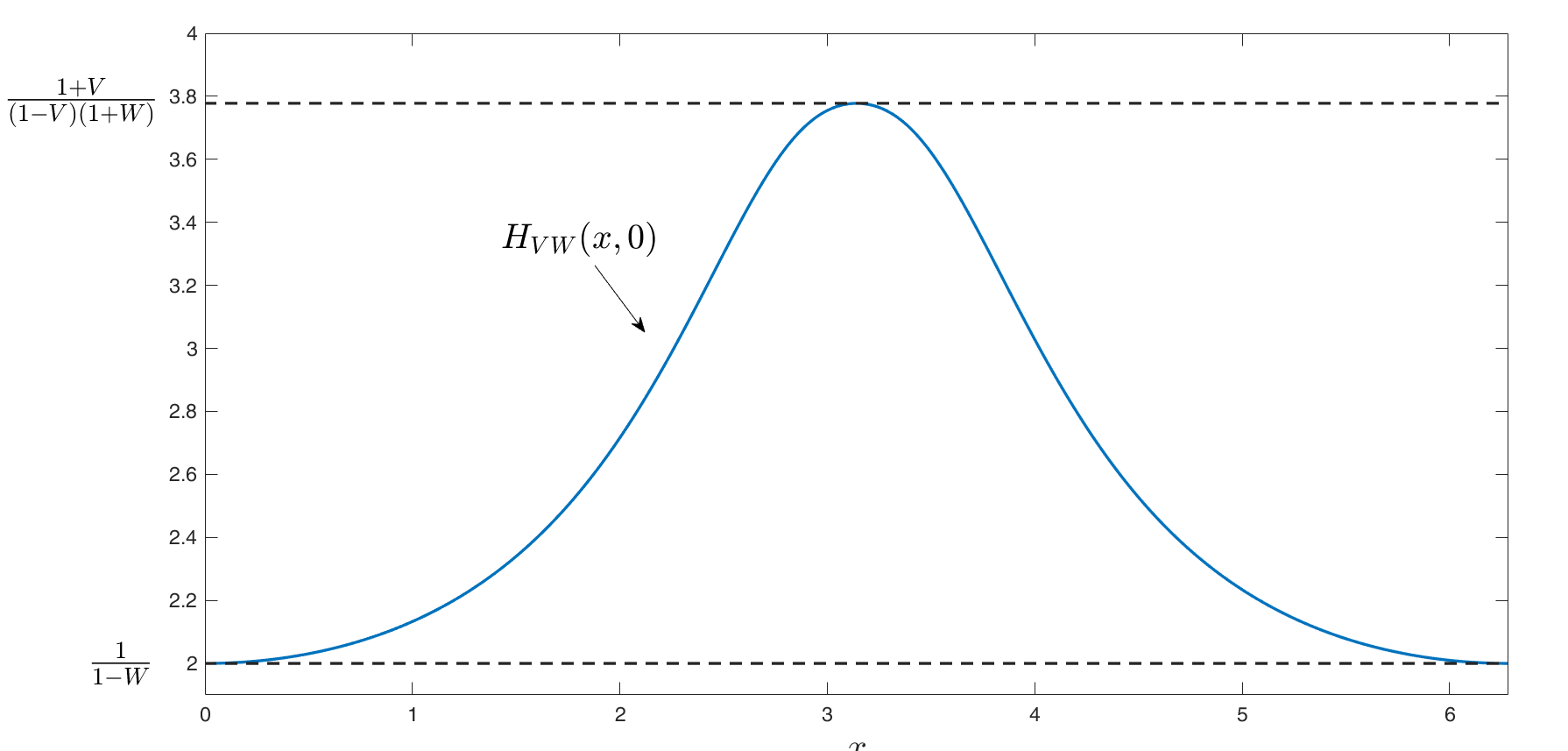}}
	
	\caption{$H_{VW}(x,0)$ for $x\in[0,2\pi]$.}
	
	\label{lastfigure}
\end{figure}

In Figure \ref{lastfigure}, we see $H_{VW}(x,0)$ for $x \in [0, 2\pi]$ in the three cases a), b), and c) in (\ref{abc}) (the values of $V$ and $W$ in this figure are the same as in Figures \ref{branches} and \ref{branchespartial} of Appendix \ref{AB}). Observe that these curves also represent a scaled version of $\mathrm{OT}(t,y_0)$ as $t$ varies (see (\ref{OTy}) in the main paper) for
$$
\Delta_1(\widehat{y}_0) =2(\gamma_1(\widehat{y}_0)-\theta_1)=0,
$$
i.e., $y_0$ is orthogonal to the second right singular vector of the matrix $R_1$ of Proposition \ref{euclidean} of the main paper (see Subsection \ref{sssection} there). In the case b), there is a stationary point that is a local minimum but not the global minimum. This fact can be also seen in Figure  \ref{FigureV1} (B) of Subsection \ref{test} in the main paper. 

\begin{thebibliography}{99}
	
	
	
	
	
	
	
	\bibitem{Mohy2008} A.~Al-Mohy and N.~Higham. Computing the Frech\'{e}t
	derivative of the matrix exponential, with an application to condition
	number estimation. SIAM Journal on Matrix Analysis and Applications, 30
	(2008/2009) no. 4, 1639-1657.
	
	
	\bibitem{Burgisser2013} P. B{\"{u}}rgisser and F.~Cucker. Condition: the
	geometry of numerical algorithms. Springer, 2013.
	
	
	
	
	\bibitem{Kagstrom1977} B.~K{\aa }gstr{\"{o}}m. Bounds and perturbation
	bounds for the matrix exponential. BIT 17 (1977) no. 1, 39-57.
	
	\bibitem{Levis1969} A.~Levis. Some computational aspects of the matrix
	exponential. IEEE Transactions on Automatic Control, AC-14 (1969) no. 4,
	410-411.
	
	\bibitem{Van1977} C.~Van Loan. The sensitivity of the matrix exponential.
	SIAM Journal on Numerical Analysis, 14 (1977) no. 6, 971-981.
	
	
	
	
	
	
	
	
	\bibitem{M0} S.~Maset. Asymptotic condition numbers for linear ordinary differential equations. arXiv (2026).
	
	\bibitem{Maset2022} S.~Maset. A long-time relative error analysis for linear ordinary differential equations with perturbed initial value. arXiv (2026). 
	
	\bibitem{moler2003} C.~Moler and C.~Van~Loan. Nineteen dubious ways to
	compute the exponential of a matrix, twenty-five years later. SIAM Review,
	45 (2003) no. 1, 3-49.
	
	\bibitem{van2006} C. Van~Loan. A study of the matrix exponential. Numerical
	Analysis Report No. 10, University of Manchester, Manchester, UK, August
	1975. Reprinted November 2006.
	
	\bibitem{Zhu2008} W.~Zhu, J.~Xue, and W.~Gao. The sensitivity of the
	exponential of an essentially nonnegative matrix. Journal of Computational
	Mathematics, 26 (2008) no. 2, 250-258.
	
	
	
	
	
\end{thebibliography}
\end{document}